\documentclass[mnsc,nonblindrev]{informs3}

\OneAndAHalfSpacedXI % current default line spacing
\usepackage{natbib}
 \bibpunct[, ]{(}{)}{,}{a}{}{,}%
 \def\bibfont{\small}%
\usepackage{setspace}
\usepackage{multirow}
\usepackage{algorithm}
\usepackage{algpseudocode}
\usepackage{url}

\TheoremsNumberedThrough     % Preferred (Theorem 1, Lemma 1, Theorem 2)
\ECRepeatTheorems

\EquationsNumberedThrough    % Default: (1), (2), ...
\MANUSCRIPTNO{MS-0001-1922.65}

\begin{document}
%%%%%%%%%%%%%%%%

%\setcounter{page}{1}
%\DOI{10.1287/xxxx.0000.0000}%


\RUNTITLE{Block Scheduling in Two-stage Outpatient Clinics: Appointment Template Design}

% \TITLE{Bundling Information Goods of Decreasing Value}
\TITLE{Block Scheduling in Two-stage Outpatient Clinics: Appointment Template Design}

\ARTICLEAUTHORS{%
\AUTHOR{Pelin Ke\c{s}rit}
\AFF{Department of Industrial \& Systems Engineering, Texas A\&M University, \\ \EMAIL{pelin@tamu.edu}} %, \URL{}}
\AUTHOR{Chelliah Sriskandarajah}
\AFF{Department of Information \& Operations Management, Mays Business School, Texas A\&M University, \EMAIL{chelliah@mays.tamu.edu}}
\AUTHOR{Jon M. Stauffer}
\AFF{Department of Information \& Operations Management, Mays Business School, Texas A\&M University, \EMAIL{jstauffer@mays.tamu.edu}}
} % end of the block

\ABSTRACT{%
Increasing the efficiency and effectiveness of the healthcare system is a challenge faced worldwide. Many outpatient clinics have implemented two-stage service systems, with both a physician and physician assistant, to enhance capacity and reduce costs. Some patients only visit a physician assistant while some patients visit both providers depending on their patient type. However, minimizing provider idle time and overtime while reducing patient waiting time is challenging in two-stage service systems. Thus, our objective is to find daily appointment templates, based on block scheduling, that minimize a weighted sum of these metrics. A block schedule divides the overall schedule into several time blocks and assigns patients of different types into each block in proportion to their daily demand to balance the workload throughout the day. Since the problem is shown to be strongly $\mathcal{NP}$-Hard, we develop a heuristic algorithm that provides a no-idle time appointment template that is easily implementable. We expand our study to include stochastic service times and show that our algorithm yields an efficient block schedule under practically relevant conditions. The algorithm is able to provide a solution similar in cost to the stochastic model when patient wait time costs are low by maintaining lower physician idle times with at most a 16 minute/patient increase in patient wait times. Comparing our heuristic to a First Come, First Appointment scheduling rule, we show that our heuristic is able to better minimize provider idle time, which mimics many real-life settings where clinics prioritize the efficiency of the healthcare providers.
% We also extend our block scheduling approach to situations with positive patient no-show probabilities.
}%

%\KEYWORDS{deterministic inventory theory; infinite linear programming duality;
%  existence of optimal policies; semi-Markov decision process; cyclic schedule}

% \HISTORY{This paper was
% first submitted on April 12, 1922 and has been with the authors for
% 83 years for 65 revisions.}

\maketitle
%%%%%%%%%%%%%%%%%%%%%%%%%%%%%%%%%%%%%%%%%%%%%%%%%%%%%%%%%%%%%%%%%%%%%%

% Samples of sectioning (and labeling) in MNSC
% NOTE: (1) \section and \subsection do NOT end with a period
%       (2) \subsubsection and lower need end punctuation
%       (3) capitalization is as shown (title style).
%
%\section{Introduction.}\label{intro} %%1.
%\subsection{Duality and the Classical EOQ Problem.}\label{class-EOQ} %% 1.1.
%\subsection{Outline.}\label{outline1} %% 1.2.
%\subsubsection{Cyclic Schedules for the General Deterministic SMDP.}
%  \label{cyclic-schedules} %% 1.2.1
%\section{Problem Description.}\label{problemdescription} %% 2.

% Text of your paper here

\section{Introduction}\label{intro}
Healthcare providers all around the world face a common challenge, that is increasing the effectiveness and efficiency of the healthcare system. This is an especially important issue in outpatient clinics due to increasing demand. Outpatient services consist of low-acuity medical services that do not require an overnight stay \citep{deceuninck2018outpatient} as opposed to inpatient services, which are more complex and require longer hospital stays. Outpatient services include but are not limited to imaging services, ambulatory surgery centers, specialized clinics, and primary care clinics. Lower costs \citep{davis1972substitution,stgeorge2021}, improved patient experience, and the use of newer technologies \citep{vitikainen2010substituting, uncwo2021} are some of the factors that attract patients to outpatient services. From 2012 to 2015, revenue obtained from outpatient services grew faster than inpatient services in the United States \citep{outpatientcare2018}. By 2016, outpatient services made up nearly 47$\%$ of total hospital revenue \citep{outpatientcare2018}. The COVID-19 pandemic contributed to the increase in demand for outpatient services as well, with nearly 40\% of physicians reporting that they are more likely to refer their patients to non-hospital locations for procedures and surgeries \citep{growth2020}. Therefore, it is important for outpatient service capacity to match this increase in demand.

Physician shortages in the United States are also projected to increase in the coming years \citep{malayala2021primary}, providing additional importance for minimizing provider idle time and overtime. Patient scheduling is a major factor that impacts efficiency in healthcare systems \citep{sun2021stochastic}. In particular, outpatient clinics need to develop appointment templates to balance minimizing provider idle time and overtime with reductions in patient waiting time \citep{zacharias2017managing}. However, satisfying rapidly increasing demand without increasing costs and decreasing patient satisfaction is a challenge for healthcare providers. Despite efforts in the literature to design an optimal appointment template that minimizes total costs, the solutions provided are often not adopted by outpatient clinics due to difficulties in implementation \citep{huang2012alternative}. Therefore, developing an implementable schedule utilizing healthcare resources efficiently (low provider idle time and overtime) while maintaining a high level of patient satisfaction is essential for outpatient clinics.

% Healthcare services are typically divided into inpatient and outpatient services based on the procedures and the length of stay of the patients. 

\subsection{Problem Definition} \label{pd}
\vspace{-1mm}

Outpatient clinics can be categorized according to the number of stages that patients experience throughout their visit. Single-stage outpatient clinics consist of a single healthcare provider, typically the physician. \cite{klassen2019appointment} show that medical clinics can decrease patient waiting time with the addition of mid-level service providers (MLSPs). MLSPs can treat moderately ill patients since they are better trained than entry level healthcare providers, but not as well trained as physicians. MLSPs generally include physician assistants (PA), nurse practitioners (NP), clinical nurse specialists (CNS), and certified nurse-midwives (CNM) \citep{bishop2012advanced}. Without loss of generality, MLSPs will be referred to as $PA$s throughout this paper. In addition, many clinic managers have included $PA$ clinic roles to reduce costs and maximize physician efficiency, as the cost of healthcare in the United States has significantly increased in the last decade \citep{kumar2011examining,nhefacts2021}. With the inclusion of $PA$s, the outpatient clinic becomes a two-stage system since $PA$s see some patients before $P$ or act as a substitute for $P$ \citep{white2017ice} depending on the patient type. However, two-stage systems complicate appointment scheduling and make it more challenging to effectively utilize this enhanced capacity. The potential for patient waiting time at each stage \citep{srinivas2020designing} and the performance of a stage being affected by the performance of an earlier stage \citep{klassen2019appointment} are some of these challenges.

% With the inclusion of MLSPs, the outpatient clinic then becomes a two-stage system. As the cost of healthcare in the US has significantly increased in the last decade \citep{kumar2011examining,nhefacts2021}, many clinic managers are seeking ways to maximize efficiency while maintaining a high level of patient care.  Many clinics have implemented MLSPs to enhance system capacity. However, their addition complicates the system and makes it challenging to effectively utilize this enhanced capacity. 

\begin{figure}[h!] %t h b
	\centering
	\includegraphics[scale=.6]{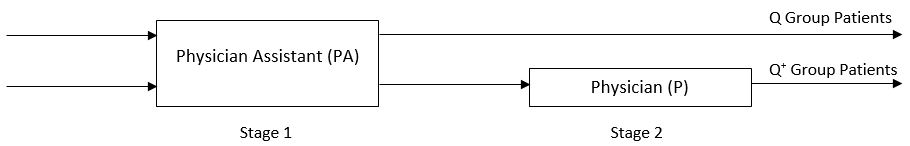}
    \caption{Patient Flow Chart}
	\label{process1}
	\vspace{-3mm}
\end{figure}

This study tackles these challenges by considering a two-stage outpatient medical clinic under patient heterogeneity (multiple patient types) with one physician assistant ($PA$) in stage~1 and one physician ($P$) in stage~2 (We discuss strategies for extending this problem setting to multiple provider systems in Appendix~\ref{general_setting}). While $PA$ serves all patients, $P$ only sees some of the patients after they exit stage~1, depending on the patient type. The system design is illustrated in Figure~\ref{process1}. All patients receive treatment first at stage~1, then they may continue to stage~2 if they need further treatment. However, if the patients are fully treated by the $PA$ in stage~1, they can exit the system without seeing the physician ($P$). Thus, upon entering the system, patients are assumed to follow one of two paths: (i) see the $PA$ and leave ($Q$ Group Patients), (ii) see $PA$ first, then see $P$ and then leave the clinic ($Q^+$ Group Patients). Without loss of generality, average service times for $P$ and $PA$ are assumed to be different from each other and dependent on the patient type. 

\vspace{-4mm}

%This study assumes that different types of patients exist in the system and can be grouped according to the similarities in the medical service that they require and other characteristics. These characteristics may include but are not limited to having scheduled specific medical procedures and being a new or recurring patient. The problem of creating a two-stage schedule consisting of both $PA$s and $P$s that maximizes provider utilization while improving patient flow can be a challenging task. By considering all of the factors above, we aim to create a daily appointment template for a two-stage ($PA$ and $P$) system in order to minimize  $P$ and $PA$ idle time and overtime, and reduce patient waiting time.

Our research is motivated by a practice-based scenario in a Women’s Health Clinic studied in \cite{huang2015appointment}, where patients are classified into types according to their service times, and the clinic needs to assign patients to a daily schedule. %By considering patient heterogeneity in appointment scheduling, provider idle time and patient waiting time can be minimized \citep{rohleder2000using} and the cost of a schedule can be reduced \citep{deceuninck2018outpatient}. 
We build on the block scheduling concept proposed by \cite{lee2018outpatient} in order to minimize provider idle time and overtime over the course of a planning horizon. In each two-stage appointment block, the clinic schedules patients by type into the appropriate slot. Then, the clinic repeats the same block several times to satisfy the demand and establish a planning horizon schedule (appointment template).

%Our research is motivated by a practice-based scenario in a Women’s Health Clinic studied in \cite{huang2015appointment}, where patients are classified into types according to their service times, and the clinic needs to assign patients to a daily schedule. By considering patient heterogeneity in appointment scheduling, provider idle time and patient waiting time can be minimized \citep{rohleder2000using} and the cost of a schedule can be reduced \citep{deceuninck2018outpatient}. We build on the block scheduling concept proposed by \cite{lee2018outpatient} in order to balance provider workload throughout the day and design a two-stage daily appointment template. In each block, the clinic schedules patients of each type according to their demand ratios, which are the number of each patient type seen in a block. Then, the clinic repeats the same block schedule several times to satisfy the daily demand and establishes a one-day schedule (appointment template) that minimizes the provider idle cost, provider overtime cost, and patient waiting time cost.

\subsection{Research Questions and Contributions}\label{prd}
\vspace{-1mm}

This study proposes a solution to the cost minimization problem of appointment template design for two-stage outpatient clinics that increases the efficiency of healthcare provider utilization while improving patient flow. Thus, the study aims to answer the main research question: \textit{How can a two-stage outpatient clinic design its appointment templates under patient heterogeneity and deterministic service times in such a way that minimizes the total system cost?} We define the total system cost as a combination of patient waiting time plus idle time and overtime of the $P$ and the $PA$. We address this research question by adopting a block scheduling policy for the appointment template and formulating a mathematical model based on block scheduling. We show that the problem of minimizing patient waiting time in a no-idle time block schedule is strongly $\mathcal{NP}$-Hard, and thus develop a heuristic algorithm.

% By appropriately overbooking patients in the block schedule, we also show that our heuristic algorithm can handle patient no-shows and still minimize provider idles times.

% Our heuristic algorithm performs well under practically relevant problem instances with deterministic service times.

Additionally, we extend our main problem by asking: \textit{How can a two-stage outpatient clinic design its appointment templates under patient heterogeneity and stochastic service times in such a way that minimizes the total expected cost?} As an extension to the primary research question, we develop a stochastic programming model solved using Sample Average Approximation (SAA) and develop a heuristic algorithm that provides a no-idle block schedule when the deviation from the mean service time is less than 46\% on average. Additionally, we compare our implementable heuristic to a first come, first appointment scheduling rule from literature and the stochastic model SAA solution. We find that our heuristic can achieve lower physician idle times than both of these approaches, while only increasing wait time per patient by 9.7 and 16 minutes to the first come, first appointment scheduling rule and SAA solution, respectively. 

The rest of the paper is organized as follows. Section~\ref{lr} reviews related literature. We introduce the block scheduling concept in two-stage outpatient clinic scheduling with mathematical modeling and heuristic approaches in Section~\ref{twosec}. Section~\ref{msmp} provides appointment templates using the mathematical model and heuristic algorithm for the planning horizon. Section~\ref{ucd} analyzes scenarios with stochastic service times. Section~\ref{compstu} presents the computational results.  Section~\ref{con} concludes and discusses future research directions.

% Section~\ref{robu} discusses the strategies for handling no-show uncertainty with block scheduling. 

\section{Literature Review} \label{lr}
\vspace{-1mm}

This research builds on both single-stage and multiple-stage healthcare scheduling literature, which we discuss in more detail in Section~\ref{lr-s} and Section~\ref{lr-m}, respectively.

\subsection{Single-Stage Appointment Scheduling} \label{lr-s}
\vspace{-1mm}

\looseness -1 Single-stage appointment scheduling literature considers clinical settings that consist of one healthcare provider. Common objectives in single-stage appointment scheduling problems are minimizing the patient waiting time, $P$ idle time \citep{cayirli2003outpatient} and overtime. Aiming to reduce these three costs, \cite{shuang2019two} first assign $P$s to days of the week according to capacity constraints, then schedule patients into daily appointments under uncertain service times. \cite{huang2015appointment} determine appropriate appointment lengths using a reclassification technique subject to clinical constraints such as patient waiting time, $P$ overtime, and idle time. \cite{anvaryazdi2020appointment} determine the number of each patient type in each appointment slot by utilizing a stochastic program to minimize the patient waiting time. They find that the greater the capacity of the provider, the better their proposed solution performs. \cite{shehadeh2021using} propose an adaptive rescheduling policy under patient heterogeneity. They use a stochastic mixed integer program to propose a solution that is near optimal when patients arrive on time or early. It is also possible to encounter studies that focus on revenue maximization based on similar performance metrics. \cite{yu2020appointment} develop an index policy to maximize the revenue of a healthcare facility with series patients, that require multiple visits to complete a patient's treatment. \cite{taiwo2023patient} show that delaying the choice of healthcare provider ($P$ or $PA$) in a single-stage system can increase the profit.

 \begin{table}[b]
 \caption{Overview of the Literature}
    \centering
     \resizebox{\textwidth}{!}{\begin{tabular}{|c|c|c|c|c|c|}
    \hline
        \multirow{2}{*}{Paper} &
        \multirow{2}{*}{Number of Stages} &
        \multirow{2}{*}{Method} &  \multicolumn{3}{c|}{Objective}
         \\
        \cline{4-6}
        & & & Patient Waiting Time & Physician Overtime & Physician Idle Time\\
        \hline
       \cite{anvaryazdi2020appointment}& 1&	Two-stage stochastic programming&	X&	&	\\
        \cite{fan2020appointment}&2&	Gradient search&	X&	X&		\\
        \cite{huang2015appointment} &1&	Reclassification technique&	&	&	\\
        \cite{klassen2019appointment}&2&	Simulation optimization&	X&	X&	X	\\
        \cite{lee2018outpatient}&1&	Block scheduling&	X&	X&	X	\\
        \shortcites{oh2013guidelines}\citet{oh2013guidelines}&2&	Stochastic integer programming&	X&	&	X	\\
        \cite{shehadeh2021using} &1&	Stochastic mixed-integer&	X&	X&	X	\\
        \cite{shuang2019two}&1&	Two-stage stochastic optimization&	X&	X&	X\\
       \cite{srinivas2020designing}&2&	Stochastic optimization&	X&	X&	X\\
        \cite{taiwo2020outpatient}&1&	Two-stage optimization&	X&	X&	X	\\
        \cite{taiwo2023patient} & 1 & Heuristic algorithm & X & X & X\\
        \cite{wang2020managing}&1&	Data-analytic approach&	X&	X&	X	\\
        \cite{white2017ice}& 2& Discrete-event simulation& X& &
        X\\
       \cite{yu2020appointment}&1 (multiple visits)&	One-step policy improvement&	&	X&	\\
        \cite{zhou2019appointment}&multiple&	Two-stage program&	X&	&	X	\\
        \cite{zhou2021sequencing}&multiple&	Standard Benders decomposition&	X&	&	X	\\
        
        Our Paper & 2 & Block scheduling & X & X & X
        \\
        \hline
    \end{tabular}}
    % \vspace{0.01in}
    \label{litrew}
    \vspace{-1mm}
\end{table}

Handling walk-in patients is also a commonly addressed challenge in single-stage appointment scheduling. \cite{taiwo2020outpatient} show that stipulating the optimal time window for patient walk-ins performs better than an open walk-in policy when the objective is to minimize the cost of patient waiting time, $P$ idle time, and $P$ overtime. \cite{wang2020managing} advise purposely leaving some appointment slots empty to anticipate walk-ins. \cite{youn2022planning} indicate that the focus on healthcare scheduling should not be limited to a single method or a small area of application but rather be comprehensive through a multidisciplinary perspective.

Among these studies, our objectives in a two-stage setting mostly align with \cite{lee2018outpatient} that minimize the total cost of patient waiting time, healthcare provider idle time, and overtime using the block scheduling approach in a single-stage setting. While we use a similar block scheduling approach, our paper differs from \cite{lee2018outpatient} by considering a two-stage healthcare system. This greatly increases the complexity of the required block schedule.

\subsection{Multi-stage Appointment Scheduling} \label{lr-m}
\vspace{-1mm}

 This paper is primarily related to the literature on multi-stage outpatient clinics and appointment scheduling problems. These studies mostly agree on the benefit of introducing $PA$s into the system for common metrics such as a reduction in patient waiting time. \cite{fan2020appointment} compare single-stage with two-stage appointment scheduling under no-shows and state the trade-offs between $P$ overtime and system congestion. The addition of $PA$s is also considered to be an approach to increase patient satisfaction without increasing the total cost. \cite{klassen2019appointment} observe that the more tasks a $PA$ can carry out, the more patient waiting time will decrease. \cite{white2017ice} show that the clinic can benefit most from $PA$s when they see the patient before $P$ if they consist of lower-cost staff and as a substitute for $P$s when they consist of higher-cost staff.

\looseness -1 Like single-stage problems, multi-stage appointment scheduling also considers a similar set of objectives. \cite{srinivas2020designing} develop a stochastic optimization model to schedule pre-booking and same-day appointments in order to minimize patient waiting time, $P$ idle and overtime as well as denied appointment requests. \cite{zhou2019appointment} study the cost minimization problem in a two-stage sequential service system where only the first stage allows appointments and additional interactions operate on a first-come, first-served basis. \cite{zhou2021sequencing} extend this study by adapting a SAA approach. 

% \textcolor{blue}{Although this study considers an outpatient clinical setting, it is possible to encounter multi-stage scheduling problems in healthcare literature such as ambulatory surgery centers where the goal is to minimize the total cost with the common metric of healthcare provider overtime cost\citep{youn2022adaptive}.}

Popular methods in solving multi-stage scheduling problems are metaheuristics and multi-agent methods \citep{marynissen2019literature}. However, it is also possible to encounter studies that adopt exact methods. Formulating the problem as a stochastic integer program, \cite{oh2013guidelines} demonstrate the necessity of different amounts of slack time in a schedule depending on the patient type under a fixed appointment times assumption. Case studies, queuing theory, and mathematical models are some of the other commonly used methods both in single and multiple-stage scheduling problems \citep{cayirli2003outpatient}. An overview of the comparable scheduling literature can be found in Table~\ref{litrew}. Despite having some common cost minimization objectives such as overtime, multi-stage scheduling problems from other areas of the healthcare industry (such as ambulatory surgery centers) focus mainly on capacity planning  \citep{youn2022adaptive}. Thus, this scheduling research does not align with the characteristics of multi-stage scheduling in an outpatient clinical setting and therefore is excluded from the scope of this review. %Despite having some common cost minimization objectives such as overtime, multi-stage scheduling problems from other areas of the healthcare industry (such as ambulatory surgery centers) focus mainly on capacity planning  \citep{youn2022adaptive}. Thus, this scheduling research does not align with the characteristics of multi-stage scheduling in an outpatient clinical setting and therefore is excluded from the scope of this review.

%Despite having some common cost minimization objectives such as overtime, multi-stage scheduling problems from other areas of the healthcare industry (such as ambulatory surgery centers) focus mainly on capacity planning \citep{youn2022adaptive}, which do not align with the characteristics of multi-stage scheduling in an outpatient clinical setting and therefore are excluded from the scope of this review.

% Even though exact methods are regarded as less popular due to complexity, it is possible to encounter studies that adopt this approach in scheduling. \cite{zhou2019appointment} consider a two-stage sequential service system and formulate it first as a stochastic program then as a two stage program by establishing linear relationships between waiting times. \cite{zhou2021sequencing} then extend this work by adapting sample average approximation approach. 

%  The current literature streams include research regarding two-stage outpatient appointment scheduling problems, albeit the definition of such systems differ in various contexts. \cite{fan2020appointment} study two-stage outpatient clinic appointment scheduling problem under patient no-shows. The existence of a second stage varies with respect to the probability that the patient will be required to have further diagnostic examination. \cite{oh2013guidelines} study the two-stage outpatient appointment scheduling under patient heterogeneity with stochastic nurse and physician service times. They demonstrate the necessity of different amounts of slack time in the schedule depending on the patient type. 

Although there is an extensive range of literature on multi-stage appointment scheduling, most of the studies either lack one of the performance metrics considered to be important for cost calculation (e.g. overtime) or their appointment systems only consider scheduling for some of the stages even though they are multi-stage systems, which reduce the complexity of the problem that they consider. Our paper aims to fill the gap in the outpatient appointment scheduling literature by extending the single-stage block scheduling concept proposed by \cite{lee2018outpatient} to two-stage outpatient clinics. \cite{cayirli2003outpatient} indicate that most of the studies on outpatient scheduling literature assume that patients are homogeneous with respect to service times. Current literature lacks the focus on two-stage outpatient clinics that consist of $PA$s and $P$s, operating under patient heterogeneity, that this research provides.

\section{Model and Analysis for Two-stage Scheduling of a Block} \label{twosec} 
\vspace{-1mm}

This section introduces the block scheduling concept in two-stage outpatient clinics. Section~\ref{bs} defines a block schedule for the problem setting and discusses the assumptions of the system. Section~\ref{msm} introduces the block schedule mathematical model and related notation, while Section~\ref{algsintro} provides appointment templates based on the heuristic algorithms developed in this section.

\subsection{Two-stage Block Scheduling} \label{bs}
\vspace{-1mm}

Block scheduling was initially introduced in the Toyota Production System to smooth manufacturing workloads by evenly distributing the production of different product types across the planning horizon on an automotive assembly line. In a similar spirit, block scheduling in healthcare refers to a cyclic scheduling approach in which a sequence of different patient types is assigned within a single block based on their smallest demand ratios. These blocks are then repeated until the overall demand for the planning horizon is fulfilled. More specifically, a block is defined as a sequence of patient types ordered by their smallest demand ratios, representing the aggregate demand distribution over the planning period. A two-stage block schedule refers to the scheduling of patients across multiple such blocks, where the blocks are repeated consecutively according to their predefined sequence until the planning horizon's demand is satisfied. Building on the work of \cite{lee2018outpatient}, we extend the concept of block scheduling to a two-stage outpatient clinic setting. The parameters relevant to this environment are summarized in Table~\ref{not0}.

%Block scheduling was first introduced in the Toyota Production System to smooth manufacturing workloads by distributing the production of different product types evenly over the planning horizon on an automotive assembly line. In a similar spirit, the block schedule is a cyclic schedule where we first assign a sequence of different patient types within a single block according to their smallest demand ratios and then multiple such blocks are repetitively scheduled until the demand for the planning horizon is met. More precisely, a single block is defined as a sequence of different patient types according to their smallest demand ratios representing the total demand during the planning horizon. Two-stage Block Schedule is the schedule of patients in multiple such blocks that are repetitively and consecutively scheduled according to the sequence given in each block until the demand for the planning horizon is met. Building on \cite{lee2018outpatient}, we extend the block scheduling concept to two-stage outpatient clinics. Parameters for this two-stage environment are listed in Table~\ref{not0}.

\begin{table}[h]
\caption{Block Schedule Notation}
    \centering
    % \resizebox{\textwidth}{!}{
    \begin{tabular}{|l l|}
    \hline
    \multicolumn{2}{|l|}{\textbf{Parameters}}\\
         \hline
         $m$ & Number of patient types.\\
         $r_i$ & Demand ratio of type $i$ patients during the planning horizon\\
          & which is also the number of type $i$ patients, where $i=1,2, \ldots, m$.\\
         $r$ & Total number of patients scheduled within a single block, $r=\sum_{i=1}^{m}r_i$.\\
        $v$ & Total number of $Q$ group patients scheduled within a single block, $v=\sum_{i=1}^{q}r_i$.\\
         $\lambda_{i}$ & The mean service time of type $i$ patient at stage~1.\\
         $\mu_{i}$ & The mean service time of type $i$ patient at stage~2.\\
         $k$ & Number of blocks to be scheduled during the planning horizon (a day).\\
         $n_i$ & Number of type $i$ patients scheduled during the planning horizon, $n_i=kr_i$.\\
         $n$ & Total number of patients scheduled during the planning horizon, $n=\sum_{i=1}^{m}n_i$.\\
         $R$ & The length of regular time in day.\\
\hline
    \end{tabular}
    % }
    % \vspace{0.01in}
    \label{not0}
    \vspace{-3mm}
\end{table}

 The common assumptions for two-stage block scheduling are summarized as follows:
\begin{itemize}
%\item $n-T$, the total number of overbooked patients scheduled.
\item There are $m$ types of patients. All patients call in advance to make appointments, which is a common practice in many health care systems \citep{cayirli2003outpatient}. The clinic knows patient types while making appointments. Factors such as medical record, age of the patient, and the number of previous visits, as well as the required treatment can help identify the patient types \citep{deceuninck2018outpatient}.
\item There are neither unscheduled walk-in patients nor scheduled patient no-shows. We consider positive no-show probabilities and examine overbooking strategies in Appendix~\ref{robu}.
\item Mean service times of the $m$ types of patients are deterministic and known but different for different stages. We later relax this assumption and introduce stochastic service times in Section~\ref{ucd} and present the computational results in Section~\ref{compstu}.
\item Patients, $P$, and $PA$ are punctual. The punctuality assumption eases the tractability of the analysis. \citep{taiwo2023patient}.
\item The clinic can flexibly assign patients to the clinic's preferred time according to its appointment template. In compliance with the literature, the patients are assigned to possible slots according to their types rather than their preferences \citep{cayirli2003outpatient}. 
\item The duration of regular time $R$, is exogenously determined by the clinic. An overtime cost will be incurred if the schedule length goes beyond time $R$.
\item  We assume that the clinic can successfully schedule $n$ patients to a given day due to the large demand for services at the clinic.
\end{itemize}

%\vspace{-0.2in}

\begin{table}[h]
\caption{Patient related parameters}
    \centering
    \begin{tabular}{|c|c|c|c|}
    \hline
        \multirow{2}{*}{Patient}
         &  \multicolumn{2}{c|}{Service Time} & \multirow{2}{*}{Ratio ($r_i$)}
         \\
        \cline{2-3}
        & $PA$ ($\lambda_i$) & $P$  ($\mu_i$) & \\
        \hline
        Type 1 & $\lambda_1$ & $\mu_1 = 0$ &  $r_1$\\
        Type 2 & $\lambda_2$ & $\mu_2 = 0$ & $r_2$\\
        \vdots & \vdots & \vdots & \vdots\\
        Type $q$ & $\lambda_q$ & $\mu_q = 0$ & $r_q$\\
        Type $q+1$ & $\lambda_{q+1}$ & $\mu_{q+1} > 0$ & $r_{q+1}$\\
        \vdots & \vdots & \vdots & \vdots\\
        Type $m$ & $\lambda_m$ & $\mu_m > 0$ & $r_m$\\
        \hline
    \end{tabular}
    % \vspace{0.05in}
    
    \label{notation1}
    \vspace{-3 mm}
\end{table}

%\vspace{-0.2in}

% Some of the assumptions are relaxed later in the sections. Those include the study of stochastic service time and patient with positive no-show probabilities.  
Table~\ref{notation1} provides the patient type details. Here the patients of type $i=1,2,\ldots,q$, collectively referred to as $Q$ Group Patients, have a positive service time at stage~1, $\lambda_i>0$, and zero service time at stage~2, $\mu_i=0$, which exemplifies the case where the patient is only served by $PA$. All other types have $\mu_i\geq\lambda_i>0$, where $i=q+1,q+2, \ldots, m$. 
These patients are collectively referred to as $Q^+$ Group Patients, and demonstrate the case where the patients receive medical service from both $PA$ and $P$ consecutively. We assume $\lambda_i \leq \mu_i$ for these patients as they require more treatment from $P$ \citep{klassen2019appointment}.  \cite{oh2013guidelines} show that a similar patient classification scheme, where the mean service time of $P$ is greater than or equal to $PA$, is broadly applicable in primary care settings. \cite{srinivas2020designing} exemplify a similar service time structure for outpatient clinics.
The smallest demand ratios of each patient type, $r_i$, $i=1,2, \ldots, m$, are assumed to be known from past data. We define workloads on $PA$ and $P$ in a block as $L_{a} =\sum_{i=1}^{m} r_i\lambda_{i}$ and $L_{p} =\sum_{i=1}^{m} r_i\mu_{i}$, respectively.

\noindent{\bf Definition of a Block Schedule for a Two-stage system:} In this environment, the block schedule is formed as follows: (i) First assign $r_i$ patient type $i$, $i=1, 2, \ldots, m$, within a single block; (ii) If $L_{a}\leq L_{p}$, the desired block is formed and this block is repeated desired number of times (say, $k$) until the demand is met for the planning horizon. The number of blocks in a planning horizon ($k$) depends on the length of the planning horizon and single block, so it is known based on clinic scheduling preferences. If $L_{a} > L_{p}$, we apply the {\it Balance-Workload Procedure}, discussed in Appendix~\ref{app0}, to obtain $L_{a}\leq L_{p}$ and then use the same approach to form a block schedule.

\noindent{\bf Illustration of Block Schedule, Example 1:}
Table~\ref{table-Pr1} presents the patient types, service times (in minutes), and demand ratios for~Example 1, where $m=4$, $q=2$, and $r=9$. Thus, $L_{a} =\sum_{i=1}^{m} r_i\lambda_{i} = 125$ and $L_{p} =\sum_{i=1}^{m} r_i\mu_{i}= 130$.  

\begin{table}
\caption{Example~1, where $m=4$, $r=\sum_{i=1}^{m} r_i= 9$, and $q=2$, for one block}
	\begin{center}
		%\large
		\begin{tabular}{|c|c|c|c|} 
			\hline
			Patient type $Ti$	& Stage 1 time,  $\lambda_{i}$ & Stage 2 time, $\mu_{i}$ & Demand ratio, $r_i$ \\
			\hline
			Type $T1$ & 10 & 0  &  3 \\
			Type $T2$ & 15 & 0  &  2 \\
			Type $T3$ & 20 & 25 & 1\\
			Type $T4$ & 15 & 35  & 3\\
			\hline
		\end{tabular}
	\end{center}
%	\vspace{0.1in}
	\label{table-Pr1}
    \vspace{-3mm}
\end{table}

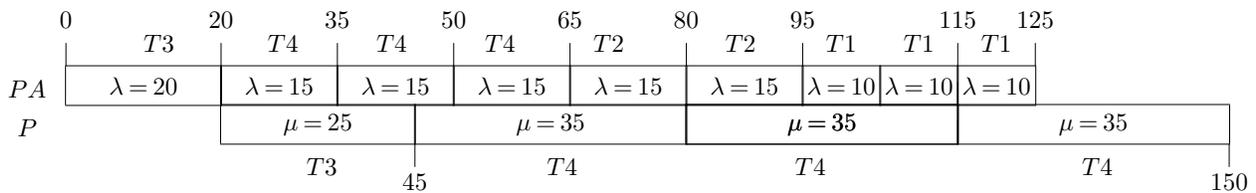
\begin{figure}[htbp]
\resizebox{\textwidth}{!}{
	{\normalsize
		\setlength{\unitlength}{1.2mm}
		\begin{picture}(158,24)(12,4) 
		%\thicklines
		
		\put(27,23){\makebox(0,0){$T3$}}
		
		\put(43,23){\makebox(0,0){$T4$}}
		
		\put(57,23){\makebox(0,0){$T4$}}
		
		\put(71,23){\makebox(0,0){$T4$}}
		
		\put(85,23){\makebox(0,0){$T2$}}
		
		\put(102,23){\makebox(0,0){$T2$}}
		
		\put(115,23){\makebox(0,0){$T1$}}
		\put(125,23){\makebox(0,0){$T1$}}
		\put(135,23){\makebox(0,0){$T1$}}
		
		\put(48,7){\makebox(0,0){$T3$}}
		\put(79,7){\makebox(0,0){$T4$}}
		\put(111,7){\makebox(0,0){$T4$}}
		\put(148,7){\makebox(0,0){$T4$}}

		\multiput(15,20)(27,0){1}{\line(0,1){3}}
		\put(15,26){\makebox(0,0){0}}
		
		\multiput(35,20)(27,0){1}{\line(0,1){3}}
		\put(35,26){\makebox(0,0){20}}
		
		\multiput(50,20)(27,0){1}{\line(0,1){3}}
		\put(50,26){\makebox(0,0){35}}
		
		\multiput(65,20)(27,0){1}{\line(0,1){3}}
		\put(65,26){\makebox(0,0){50}}
		
		\multiput(80,20)(27,0){1}{\line(0,1){3}}
		\put(80,26){\makebox(0,0){65}}
		
		\multiput(95,20)(27,0){1}{\line(0,1){3}}
		\put(95,26){\makebox(0,0){80}}
		
		\multiput(110,20)(27,0){1}{\line(0,1){3}}
		\put(110,26){\makebox(0,0){95}}

		\multiput(130,20)(27,0){1}{\line(0,1){3}}
		\put(130,26){\makebox(0,0){115}}

		\multiput(140,20)(27,0){1}{\line(0,1){3}}
		\put(140,26){\makebox(0,0){125}}

		\multiput(60,7)(27,0){1}{\line(0,1){3}}
		\put(60,5){\makebox(0,0){45}}
		
		\multiput(165,7)(27,0){1}{\line(0,1){3}}
		\put(165,5){\makebox(0,0){150}}

		\multiput(130,15)(27,0){1}{\line(0,1){3}}
		\multiput(130,10)(27,0){1}{\line(0,1){3}}
		%\multiput(130,5)(27,0){1}{\line(0,1){3}}
		
		%======PA ====
		
		\put(15,15){\framebox(20,5){$\lambda=20$}}
		\put(35,15){\framebox(15,5){$\lambda=15$}}
		\put(50,15){\framebox(15,5){$\lambda=15$}}
		\put(65,15){\framebox(15,5){$\lambda=15$}}
		\put(80,15){\framebox(15,5){$\lambda=15$}}
		\put(95,15){\framebox(15,5){$\lambda=15$}}
		
		\put(110,15){\framebox(10,5){$\lambda=10$}}
		\put(120,15){\framebox(10,5){$\lambda=10$}}
		\put(130,15){\framebox(10,5){$\lambda=10$}}

		\put(10,17){\makebox(0,0){$PA$}}
		
		%======P==
		\put(35,10){\framebox(25,5){$\mu=25$}}
		\put(60,10){\framebox(35,5){$\mu=35$}}
		\put(95,10){\framebox(35,5){$\mu=35$}}
		
		\put(130,10){\framebox(35,5){$\mu=35$}}
		
		\put(95,10){\framebox(35,5){$\mu=35$}}

		\put(10,12){\makebox(0,0){$P$}}
		
		\end{picture}}}	
	%	\vspace{0.1in}	
	\caption{A No-Idle Block Schedule $\pi_1$ for Example~1. $\sigma=\{T_3, T_4, T_4, T_4, T_2,T_2,T_1,T_1,T_1\}$ is the patient sequence in $\pi_1$. The total patient waiting time in $\pi_1$ is 90.}
	\label{figx:block1}	
    % \vspace{-4mm}
\end{figure}

A schedule is created in which there is no idle time for $P$ (respectively, $PA$) from the start time of the block schedule until the finish time of $P$ (respectively, $PA$) in the block. More precisely, the idle time is measured with respect to the start and the end times of the first and the last appointments within the block for both $P$ and $PA$. A No-Idle Block Schedule for Example~1 is shown in Figure~\ref{figx:block1}.

\begin{figure}[htbp]
	{\normalsize
		\setlength{\unitlength}{1.1mm}
		\begin{picture}(100,13)(5,7) 
		%\thicklines
		
		\put(15,17){\makebox(0,0){$PA$}}
		
		\put(15,12){\makebox(0,0){$P$}}
				
		%======Head1 ====
		
		\put(20,15){\framebox(15,5){$|H_1|=20$}}
				
		%===========Body1=============
		
		\put(35,10){\framebox(25,10){$|B_1|=105$}}
		
		%===========Tail1=============
		
		\put(60,10){\framebox(20,5){$|T_1|=25$}}		
		
		\end{picture}}	
	\vspace{-0.2in}	
	\caption{The Block Schedule $\pi_1$  with $H_1$, $B_1$, $T_1$ for Example~1}
	\label{figx:block2}	
    % \vspace{-8mm}
\end{figure}
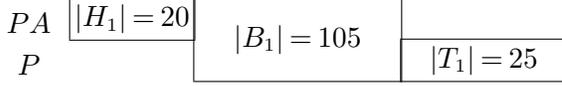

In a planning horizon, there are $k$ blocks, each block $\pi_c$, $c=1,2, \ldots, k$, can be divided into three distinct sections, the head ($H_c$), body ($B_c$), and tail ($T_c$) (see Figure~\ref{figx:block2}).

\noindent{\bf Head Section:} It is the beginning portion of $\pi_c$, where $P$ is idle throughout. We call this section the ``head" of block schedule $\pi_c$. We define $|H_c|$ as the length of the head section. 

%\medskip

\noindent{\bf Body Section:} It is the middle portion of $\pi_c$ where both $P$ and $PA$ are working. We call this section the ``body" of block schedule $\pi_c$. We define $|B_c|$ as the length of the body section. 

%\medskip

\noindent{\bf Tail Section:} It is the last portion of $\pi_c$, where $PA$ is idle throughout. We call this section the ``tail" of block schedule $\pi_c$. The length of the tail section is measured by $|T_c| = C_c - (|H_c| + |B_c|)$, that is the difference between the completion time of the schedule ($C_c$) and the sum of the duration of the head and the body sections.

%Consider a multiple block schedule as shown in Figure~\ref{figx:block00}: the schedule consists of $k$ jobs (corresponding to $k$ blocks), each job having a head $H$, a body $B$, and a tail $T$, where $k=2$,

\vspace{-5mm}

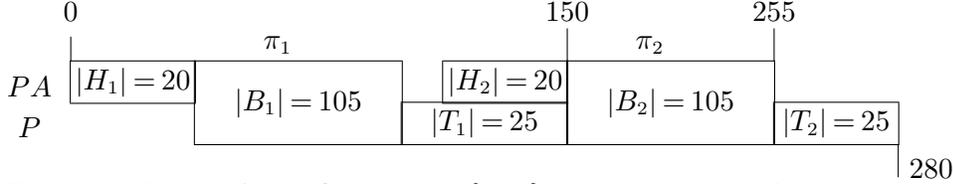
\begin{figure}[htbp]
	{\normalsize
		\setlength{\unitlength}{1.1mm}
		\begin{picture}(100,15)(10,15) 
		%\thicklines

		\put(15,17){\makebox(0,0){$PA$}}
		
		\put(15,12){\makebox(0,0){$P$}}
		
			\multiput(20,20)(27,0){1}{\line(0,1){3}}
		\put(20,26){\makebox(0,0){0}}

		%======Head1 ====
		
		\put(20,15){\framebox(15,5){$|H_1|=20$}}

		%===========Body1=============
		
		\put(35,10){\framebox(25,10){$|B_1|=105$}}
		
		%===========Tail1=============
		
		\put(60,10){\framebox(20,5){$|T_1|=25$}}

		%======Head2 ====
		
		\put(65,15){\framebox(15,5){$|H_2|=20$}}

		%===========Body2=============
		
		\put(80,10){\framebox(25,10){$|B_2|=105$}}
		
		%===========Tail2=============
		
		\put(105,10){\framebox(15,5){$|T_2|=25$}}

		%======Head3 ====

		\multiput(105,20)(27,5){1}{\line(0,1){4}}
		\put(105,26){\makebox(0,0){$255$}}
		
			\multiput(80,20)(27,5){1}{\line(0,1){4}}
		\put(80,26){\makebox(0,0){$150$}}
		
			\multiput(120,6)(27,5){1}{\line(0,1){4}}
		\put(124,7){\makebox(0,0){$280$}}

		\put(45,22){\makebox(0,0){$\pi_1$}}
		\put(90,22){\makebox(0,0){$\pi_2$}}
		
		\end{picture}}	
	\vspace{0.3in}	
	\caption{Example~1: The Schedule $\Gamma= \{\pi_1\pi_2\}$ is a concatenation of two blocks. $\sigma$ is patient sequence in $\pi_1$ and $\pi_2$.}
	\label{figx:block00}	
    % \vspace{-4mm}
\end{figure}

% \text{

% }
 \vspace{0.25in}
% \vspace{-0.2in}

Let $\Gamma$ denote the schedule consisting of $k$ blocks. If $k=2$, this means repeating the patients described in Table~\ref{table-Pr1} twice to satisfy the daily demand. Figure~\ref{figx:block00} shows the no-wait block schedule $\Gamma = \{\pi_1\pi_2\}$ of two blocks, where the no-wait constraint is maintained between the blocks in scheduling terminology. This is useful in the sense that the system idle time and patient waiting time remain the same in each block even after forming the no-wait schedule. The reader may refer to \cite{sriskandarajah1986some} for the definition of no-wait schedules. 

% Thus, the total idle time of $P$, $D_p(\sigma) = 0$, the total idle time of $PA$, $D_a(\sigma)=5$, and the patient waiting time at stages~1 and 2, respectively, $W_p(\sigma) = 180$, $W_a(\sigma)=0$. 

% As illustrated in Example~1, this study aims to find a no-idle time block schedule for two-stage outpatient appointments under patient heterogeneity with deterministic service times and zero no-show probabilities. For the problem setting described in Section~\ref{pd}, the study conjectures the existence of a no-idle time block schedule for practically relevant conditions ($\lambda>0,\, \mu=0 \text{ for $Q$ Group Patients; } \lambda\leq\mu,\text{ for $Q^+$ Group Patients}$), as discussed in Section~\ref{noshow}.

\begin{table}[b]
\caption{Notation for Single Block Scheduling Model}
    \centering
    \resizebox{\textwidth}{!}{\begin{tabular}{|l l|}
         \hline
         \multicolumn{2}{|l|}{\textbf{Parameters}}\\
         \hline
         $\alpha$ & Unit cost for patient waiting time.\\
         $\beta_{p}$ & Unit cost for $P$ idle time.\\
         $\beta_{a}$ & Unit cost for $PA$ idle time.\\
         $\pi$ & The schedule of a block.\\
          $\sigma$ & The sequence pf patients in $\pi$.\\
         $M$ & Sufficiently large number\\
         \hline
         \multicolumn{2}{|l|}{\textbf{Decision Variables}}\\
         \hline
         $x_{lj}$ & Binary variable that takes the value 1 if patient $l$ is assigned to the appointment slot $j$ and 0\\
         & otherwise, where $j=1,2,\cdots,r$ and the set of all patients in a block, $l=1,2,\cdots,r$.\\ 
         $y_j$ & Binary variable that takes the value 1 if the slot $j$ is assigned to a patient that has a nonzero\\
         & service time at $P$ and 0 otherwise, where $j=1,2,\cdots,r$. \\
         \hline
         \multicolumn{2}{|l|}{\textbf{Internal Variables}}\\
         \hline
         $e^j_a$ & The start time of patient in slot $j$ at stage 1 in $\pi$, where $j=1,2, \ldots, r$. \\
         $e^j_p$ & The start time of patient in  slot $j$ at stage 2 in $\pi$,  where $j=1,2, \ldots, r$. \\
         $f^j_a$ & The finish time of  patient in  slot $j$ at stage 1 in $\pi$, where $j=1,2, \ldots, r$. \\
         $f^j_p$ & The finish time of patient in  slot $j$ at stage 2 in $\pi$,  where $j=1,2, \ldots, r$. \\
         $w^j_a$ & The waiting time of patient in  slot $j$ at stage 1 in $\pi$, where $j=1,2, \ldots, r$. \\
         $\bar{w}^j_p$ & The waiting time of patient in  slot $j$ at stage 2 in $\pi$ (denoted by $w^j_p$), if the patient has nonzero\\
         & service time at stage 2 (denoted by $\bar{w}^j_p$), where $j=1,2, \ldots, r$.\\
         $d^j_a$ & The $PA$ idle time just before serving the patient in  slot $j$ at stage 1.\\
         $d^j_p$ & The $P$ idle time just before serving the patient in  slot $j$ at stage 2.\\
         \hline
    \end{tabular}
    }
    % \vspace{0.02in}
    \label{not1}
    \vspace{-3mm}
\end{table}

% \vspace{-0.05in}
\subsection{Two-stage Scheduling Model for a Block} \label{msm}
% Appointment Templates for Two-stage Scheduling of a Block
\looseness -1 In order to address the first research question discussed in Section~\ref{prd}, we develop a two-stage scheduling model. Table~\ref{not1} describes the parameters and variables in this model. This is in addition to the block parameters in Table~\ref{not0}. The variables that determine the patient and appointment slot assignments are the main decision variables, whereas the internal variables are calculated based on these slot assignments.

Since we are considering a single block in this model, overtime is ignored until we consider the model for a planning horizon (or full day) in Section~\ref{msmp}. Thus, we formulate the mathematical model that minimizes the weighted sum of the cost of patient waiting times and the idle times of the $P$ and $PA$, for a single block schedule.
% in Section~\ref{app1.0} of E-Companion.

\textbf{Single Block Model}
\begin{align}
    \min\; & \alpha \sum_{j=1}^r \bar{w}^j_p+\beta_a \sum_{j=1}^r d^j_a+\beta_p \sum_{j=1}^r d^j_p\label{m1}\\
    \text{subject to } & \sum_{j=1}^r x_{lj}=1\hspace{3mm} l=1,2,\cdots,r\label{m2}\\
    & \sum_{l=1}^{r}  x_{lj}=1\hspace{3mm}  j=1,2,\cdots,r\label{m3}\\
    & \sum_{s\in\{1,2,\cdots,q\}} x_{s1}=0\label{m35}\\
    & e_a^{j+1}\geq e_a^j + \sum_{l=1}^r x_{lj} \lambda_l \hspace{3mm}  j=1,2,\cdots,r-1\label{m4}\\
    & e_p^j \geq e_a^j + \sum_{l=1}^r x_{lj} \lambda_l \hspace{3mm} j=1,2,\cdots,r\label{m5}\\
    & e_p^{j+1} \geq e_p^j + \sum_{l=1}^r x_{lj} \mu_l \hspace{3mm} j=1,2,\cdots,r-1\label{m6}\\
    & d_a^{j+1} \geq e_a^{j+1}-e_a^j-\sum_{l=1}^r x_{lj} \lambda_l \hspace{3mm} j=1,2,\cdots,r-1\label{m7}\\
    & d_p^{j+1} \geq e_p^{j+1}-e_p^j-\sum_{l=1}^r x_{lj}\mu_l \hspace{3mm}  j=1,2,\cdots,r-1\label{m8}\\
    & e_a^1=0\label{m9}\\
    & \bar{w}_p^j \geq e_p^j-e_a^j-\sum_{l=1}^r x_{lj} \lambda_l -M(1-y_j)\hspace{3mm} j=1,2,\cdots,r\label{m10}\\
    & \frac{y_j}{M} \leq \sum_{l=1}^r \mu_l x_{lj} \hspace{3mm} j=1,2,\cdots,r\label{m11}\\
    & My_j \geq \sum_{l=1}^r \mu_l x_{lj}\hspace{3mm} j=1,2,\cdots,r\label{m12}\\
    & x_{lj}\in \{0,1\}\hspace{3mm} l=1,2,\cdots,r;\hspace{3mm} j=1,2,\cdots,r\label{m13}\\
    & y_j\in \{0,1\}\hspace{3mm} j=1,2,\cdots,r\label{m14}\\
    & d_a^j, d_p^j, e_a^j, e_p^j, \bar{w}_p^j \geq 0\hspace{3mm}   j=1,2,\cdots,r\label{m15}
    %& M\text{: sufficiently large number}\label{m16} \nonumber
\end{align}

The objective function (\ref{m1}) is the weighted sum of the total patient waiting time, $PA$ idle time, and $P$ idle time costs. As per the problem setting with punctual patients, the patient waiting time at stage~1 is assumed to be zero and thus, is not included in the objective function. Constraints (\ref{m2}) and (\ref{m3}) ensure that each patient is assigned to one appointment slot and each appointment slot has one patient assigned to it respectively. Constraint (\ref{m35}) ensures that no patient with zero service time at stage 2 will be scheduled at the first position. Constraint (\ref{m4}) indicates that a patient cannot start his or her appointment at stage~1 before the previous patient leaves stage~1. Constraints (\ref{m5}) and (\ref{m6}) state that a patient can start his or her appointment at stage~2 either immediately after completing his or her appointment at stage~1 or after the previous patient leaves stage~2, whichever is greatest. Constraint (\ref{m7}) indicates that the idle time of $PA$ at each appointment slot is greater than or equal to the difference between the start time of the patient at stage~1 and the finish time of the previous patient at stage~1. Similarly, constraint (\ref{m8}) indicates that the idle time of $P$ at each appointment slot is greater than or equal to the difference between the start time of the patient at stage~2 and the finish time of the previous patient at stage~2. Constraint (\ref{m9}) denotes that the first patient starts his or her appointment at stage~1 at time zero. Constraint (\ref{m10}) defines the waiting time of each patient at stage~2 as the difference between the start time at stage~2 and the finish time at stage~1 whenever the patient has a nonzero service time at stage~2. Constraints (\ref{m11}) and (\ref{m12}) identify the patients that have nonzero service times at stage~2 by using a binary variable. Constraints (\ref{m13}) and (\ref{m14}) define the binary variables whereas the constraint (\ref{m15}) ensures non-negativity of all of the remaining variables.

\vspace{-1mm}
\subsection{Single Block Schedule Minimizing System Idle Time} \label{algsintro}
\vspace{-1mm}

We first examine the environment where patients always show up as scheduled with deterministic service times. We study this class of schedules because a clinic may be interested in finding a class of schedules which eliminate the idle time and overtime of $PA$ and $P$, while patients waiting time is also reduced. Later in Section~\ref{ucd}, we introduce stochastic service times for each patient type. 

\noindent{\bf Patient Types:} Recall that patients are classified according to their service requirements in two stages. $Q$ Group Patients have $\lambda_i>0$  and $\mu_i=0$, and $Q^+$ Group Patients have $0<\lambda_i\leq\mu_i$. In each block, $r_i$ patients of type $i$ should be scheduled. 

 Among the block schedules with zero $PA$ and $P$ idle time, we are interested in those that minimize the patient waiting time. Thus, we describe Problem~$P_1$ below. %Since Problem~$P_1$ is shown to be $\mathcal{NP}$-Hard in general (Theorem~\ref{thm-block1}), we develop heuristic Algorithm~\ref{alg-MA} to solve it. 

\noindent {\bf Problem $P_1$}: Find a block schedule $\pi$  with zero total idle time such that the total patient waiting time $W_p(\pi)+W_a(\pi)$ is minimized, where $W_a(\pi)$, $W_p(\pi)$ denote patient waiting time at stages 1, 2, respectively in $\pi$.

Note that $r$ patients are scheduled in $\pi$. Let $\sigma$ be a sequence of patients in $\pi$ given appointments at the clinic in a block, where $\sigma=\{\sigma_{1}, \sigma_{2}, \ldots, \sigma_{r}\}$ and $\sigma_{j}$ denotes the patient sequenced in $j^{th}$ position (or slot $j$) in $\pi$. The appointment times are assigned to the patients in $\pi$ as follows. Let $\tau^{\sigma_j}_a$ denote the appointment time of the patient sequenced in the $j^{th}$ position.

\noindent \textbf{Appointment times:} $\tau^{\sigma_j}_a=e^j_a$, where $e^j_a$ is the start time of $j^{th}$ patient at stage 1 in $\pi$, where $j=1,2, \ldots, r$. Note that $\tau^{\sigma_1}_a=0$ and $W_a(\pi)=0$ for any sequence $\sigma$ in $\pi$ due to our assumption that patients, $P$, and $PA$ are punctual.

\noindent Minimizing patient waiting time, we start with block schedule $\pi$ having zero system idle time and propose Theorem~\ref{thm-block1}. The reader may refer to Appendix~\ref{app1} for proofs of theorems and lemmas in this section.

\begin{theorem} \label{thm-block1}
Minimizing patient waiting time in a block schedule $\pi$ having zero system idle time (Problem $P_1$) is strongly $\mathcal{NP}$-Hard.
\end{theorem} 

% \begin{algorithm}
% \begin{algorithmic}
% \caption{No-idle Block Schedule Heuristic}
% \label{alg-A}
%  \Statex \textbf{Step 1}: Arrange patients such that $\lambda_{q+1} \geq \lambda_{q+2} \geq \ldots \geq \lambda_m$, for $i=q+1,q+2, \ldots, m$, ($Q^+$ Group Patients).
%  \Statex \textbf{Step 2:} Schedule patients at earliest available time starting with all type $q+1$, followed by all type $q+2$, and so on until scheduling all type $m$, and finally all $Q$ Group Patients at the end in any order.
% \end{algorithmic}
% \end{algorithm}

Note that Problem $P_1$ is the mathematical representation of the problem. Since this problem is $\mathcal{NP}$-Hard (Theorem 1), we develop heuristic algorithms to solve Problem $P_1$. Algorithm~\ref{alg-MA} obtains a no-idle block schedule. The complexity of Algorithm~\ref{alg-MA} is in the order of $O(r \,log \, r)$ time as it involves sorting of $r$ patients in Step 1 and Step 2 can be done in $O(r)$ time.

\begin{algorithm}
\begin{algorithmic}
% \begin{singlespace}
\caption{No-idle Block Schedule Heuristic for a Single Block}
\label{alg-MA}

    \Statex \textbf{Step 1}: Arrange patients such that $\lambda_{q+1} \geq \lambda_{q+2} \geq \ldots \geq \lambda_m$,
     \Statex \hskip3.5em where $r_i=1$, $i=1,\cdots,m$.
     \Statex \hskip3.5em If there is tie $\lambda_{i} = \lambda_{i+1}$ the tie is broken in favor of scheduling the patient
     \Statex \hskip3.5em with the smallest value of $(\mu_i - \lambda_{i})$ {\bf (or equivalently the smallest, $\mu_i$)}.
     \Statex \hskip3.5em Front sequence, $\sigma^f =\{\sigma_{1}, \sigma_{2}, \ldots, \sigma_{m-q}\}= \{q+1, q+2, \ldots, q+m\}$.
     \Statex \textbf{Step 2:} Arrange patients $i$, $i=1,2, \ldots, q$ ($Q$ Group Patients) in any order to form 
     \Statex \hskip3.5em the back sequence.
     \Statex \hskip3.5em Back sequence, $\sigma^b=\{\sigma_{m-q+1}, \sigma_{m-q+2}, \ldots, \sigma_{m}\}= \{1, 2, \ldots, q\}$.
      \Statex \textbf{Step 3:} Sequence, $\sigma= \sigma^f \cup \sigma^b$.  
% \end{singlespace}
\end{algorithmic}
\end{algorithm}

For simplicity in exposition, we assume that $r_i=1$ for all $i=1,\cdots,m$. The results in the Lemmas below are valid for any general $r_i \geq1$ for all $i=1,\cdots,m$ for the following reason: if $r_i > 1$, replicate those patients with $r_i$ distinguishing patients having each $r_i=1$, where $m = \sum_{i=1}^{} r_i$. This converted problem is equivalent to that of the problem with $r_i=1$ for all $i$. Let Algorithm~\ref{alg-MA} provide a sequence, $\sigma$ which is a sequence of patients in the appointment block schedule $\pi$, where $\sigma=\{\sigma_{1}, \sigma_{2}, \ldots, \sigma_{m-q}, \sigma_{m-q+1}, \ldots, \sigma_{m}\}$ and $\sigma_{j}$ denotes the patient type sequenced in the $j^{th}$ position in $\sigma$. Consider two partial sequences, $\sigma^f$ and $\sigma^b$ of $\sigma$, where $\sigma= \sigma^f \cup \sigma^b$:\\ (i) Front sequence, $\sigma^f =\{\sigma_{1}, \sigma_{2}, \ldots, \sigma_{m-q}\}= \{q+1, q+2, \ldots, m\}$, \\(ii) Back sequence, $\sigma^b=\{\sigma_{m-q+1}, \sigma_{m-q+2}, \ldots, \sigma_{m}\}= \{1, 2, \ldots, q\}$.

In Lemma~\ref{L-block1}, we show that there exists a block schedule $\pi$ with zero $PA$ and $P$ idle time.

\begin{lemma}\label{L-block1}
Algorithm~\ref{alg-MA} provides a No-Idle Block Schedule.
\end{lemma}

% \begin{algorithm}
% \begin{algorithmic}
% \caption{Algorithm $A$}
%  \Statex \textbf{Step 1}: Arrange patients such that $\lambda_{q+1} \geq \lambda_{q+2} \geq \ldots \geq \lambda_m$, for $i=q+1,q+2, \ldots, m$, 
%  \Statex \hskip3.5em ($Q^+$ Group Patients).
%  \Statex \textbf{Step 2:} Schedule patients at earliest available time starting with all type $q+1$, 
%  \Statex \hskip3.5em followed by all type $q+2$, and so on until
%  \Statex \hskip3.5em scheduling all type $m$, and finally all $Q$ Group Patients at the end in any order.
%   \Statex \hskip3.5em Without loss of generality it is assumed that $r_i=1$ for all $i=1, \cdots, m$.
% \end{algorithmic}
% \label{alg-A}
% \end{algorithm}

\begin{lemma}\label{l-block2}
The total patient waiting time in a no-idle time block schedule $\pi$ obtained by Algorithm~\ref{alg-MA}, is $W_a(\pi)+W_p(\pi)$, where $W_a(\pi)+W_p(\pi) = \sum_{j=1}^{\gamma-1} (\gamma-j)(\mu_{v+j} - \lambda_{v+j+1}), \,\, \hbox{where} \,\, \gamma= r-v$.
\end{lemma}

As minimizing patient waiting time in a no-idle time block schedule is $\mathcal{NP}$-Hard, we aim to reduce the patient waiting time by Algorithm~\ref{alg-MA}. This yields a no-idle block schedule with a unique front schedule that has the minimum patient waiting time among all possible solutions. The algorithm schedules $Q^+$ Group patients in the front sequence and $Q$ Group patients in the back sequence. Front and back sequences are defined according to their position in the block and each consists of patients from the same group. Together they form the block schedule. The no-idle block schedule obtained by Algorithm~\ref{alg-MA} is illustrated in Figures~\ref{figx:block1} and \ref{figx:block2}.

In Algorithm~\ref{alg-MA}, we first arrange patients in $Q^+$ Group such that $\lambda_{q+1} \geq \lambda_{q+2} \geq \ldots \geq \lambda_m$. Let $S$ be the set of all possible schedules such that patients are ordered: $\lambda_{q+1} \geq \lambda_{q+2} \geq \ldots \geq \lambda_m$. $S_1$ is the set of all possible schedules that may be obtained by Algorithm~\ref{alg-MA}. Note that $S_1 \subseteq S$.
% Algorithm~\ref{alg-MA} provides one schedule in $S_1$. 

\begin{lemma}\label{modif}
The schedule $s' \in S_1$ given by Algorithm~\ref{alg-MA} yields a front sequence that provides the smallest total patient waiting time of any sequence $\bar{s} \in \{S-S_1\}$.
\end{lemma}

As a consequence of Step~1 of Algorithm~\ref{alg-MA}, the following condition is satisfied for the front sequence, $\sigma^f =\{\sigma_{1}, \sigma_{2}, \ldots, \sigma_{m-q}\}= \{q+1, q+2, \ldots, m\}$: If $\lambda_{\sigma(i)} = \lambda_{\sigma(i+1)}$ then $(\mu_{\sigma(i)} - \lambda_{\sigma(i)}) \leq (\mu_{\sigma(i+1)} - \lambda_{\sigma(i+1)})$  (or equivalently, $\mu_{\sigma(i)}  \leq \mu_{\sigma(i+1)}$), where $q+1 \leq i \leq m-1$. Recall that $r=\sum_{i=1}^{m} r_i$ and  $v$ is the total number of $Q$ group patients scheduled within one block, where $v=\sum_{i=1}^{q}r_i$. Algorithm~\ref{alg-MA} reorders $(r-v)$ patients of $Q^+$ Group Patients according to their mean service times such that $\lambda_{v+1} \geq \lambda_{v+2} \geq \ldots \geq \lambda_r$. We now develop an upper bound on the total patient waiting time in a no-idle  block schedule $\pi$ obtained by Algorithm~\ref{alg-MA}.

%\bigskip

\begin{lemma}\label{Perf-bound-block1}
Total patient waiting time in a no-idle  block schedule $\pi$ obtained by Algorithm~\ref{alg-MA} is bounded above by $\frac{(r-v)(r-v-1)}{2}[\gamma_1 - \gamma_2]$, where
	$\gamma_1 =\argmax_{1 \leq j \leq r-v-1}\{\mu_j\}$ and $\gamma_2 =\argmin_{2 \leq j \leq r-v}\{\lambda_j\}$. This bound is tight.
\end{lemma}

% Total patient waiting time in a no-idle  block schedule $\sigma$ given by Algorithm~\ref{alg-MA} is bounded above by $\frac{(r-v)(r-v-1)}{2}[\mu_{max} - \lambda_{min}]$, where $\mu_{max} = \max \{\mu_{\sigma_{k}}, \, k=1, 2, \ldots, r-v \}$  and $\mu_{min} = \min \{\lambda_{\sigma_{k}}\}$, $k=1, 2, \ldots, r-v $. This bound is tight.

% $a_1=\underset{1\leq i \leq r-v}{\argmax} \{\mu_{\sigma_i}\}, \hspace{8mm} a_2=\underset{2\leq i \leq r-v+1}{\argmax} \{\lambda_{\sigma_i}\}-1$

% The reader may verify that  $W_p(\sigma^f) =  {(r-v)(r-v-1) \over 2}[\mu_{max} - \lambda_{min}]=120$ for the example given in Table~\ref{table-Pr1}, where $r=9$, $v=5$, and $\mu_{max} - \lambda_{min} =20$.

\begin{figure}[htbp]
\resizebox{\textwidth}{!}{
	{\normalsize
		\setlength{\unitlength}{1.15mm}
		\begin{picture}(160,24)(12,4) 
		%\thicklines
		
		\put(27,23){\makebox(0,0){$T3$}}
        \put(41,23){\makebox(0,0){$T1$}}
		\put(53,23){\makebox(0,0){$T4$}}
        \put(66,23){\makebox(0,0){$T1$}}
        \put(75,23){\makebox(0,0){$T1$}}
		\put(86,23){\makebox(0,0){$T4$}}
		\put(107,23){\makebox(0,0){$T2$}}
		\put(122,23){\makebox(0,0){$T4$}}

		\put(48,7){\makebox(0,0){$T3$}}
		\put(79,7){\makebox(0,0){$T4$}}
		\put(111,7){\makebox(0,0){$T4$}}
		\put(148,7){\makebox(0,0){$T4$}}

		\multiput(15,20)(27,0){1}{\line(0,1){3}}
		\put(15,26){\makebox(0,0){0}}
		
		\multiput(35,20)(27,0){1}{\line(0,1){3}}
		\put(35,26){\makebox(0,0){20}}
		
		\multiput(60,20)(27,0){1}{\line(0,1){3}}
		\put(60,26){\makebox(0,0){45}}
		
		\multiput(95,20)(27,0){1}{\line(0,1){3}}
		\put(95,26){\makebox(0,0){80}}
		
		\multiput(130,20)(27,0){1}{\line(0,1){3}}
		\put(130,26){\makebox(0,0){115}}

		\multiput(165,7)(27,0){1}{\line(0,1){3}}
		\put(165,5){\makebox(0,0){150}}

		\multiput(130,15)(27,0){1}{\line(0,1){3}}
		\multiput(130,10)(27,0){1}{\line(0,1){3}}
		%\multiput(130,5)(27,0){1}{\line(0,1){3}}
		
		%======PA ====
		
		\put(15,15){\framebox(20,5){$\lambda=20$}}
        \put(35,15){\framebox(10,5){$\lambda=10$}}
  
		\put(45,15){\framebox(15,5){$\lambda=15$}}
        \put(60,15){\framebox(10,5){$\lambda=10$}}
        \put(70,15){\framebox(10,5){$\lambda=10$}}
  
		\put(80,15){\framebox(15,5){$\lambda=15$}}
        \put(100,15){\framebox(15,5){$\lambda=15$}}
		\put(115,15){\framebox(15,5){$\lambda=15$}}

		\put(10,17){\makebox(0,0){$PA$}}
		
		%======P==
		\put(35,10){\framebox(25,5){$\mu=25$}}
		\put(60,10){\framebox(35,5){$\mu=35$}}
		\put(95,10){\framebox(35,5){$\mu=35$}}
		
		\put(130,10){\framebox(35,5){$\mu=35$}}
		
		\put(95,10){\framebox(35,5){$\mu=35$}}

		\put(10,12){\makebox(0,0){$P$}}
		
		\end{picture}}}	
		\vspace{-0.3in}	
	\caption{An illustration of Steps 4 and 5 of Algorithm~\ref{alg-B} for Example~1}
	\label{algB_ex2}	
    \vspace{-5mm}
\end{figure}
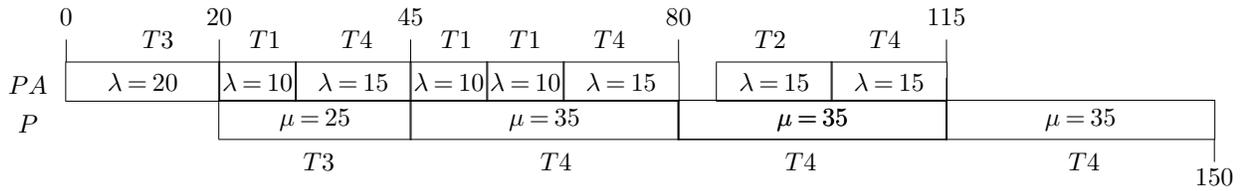

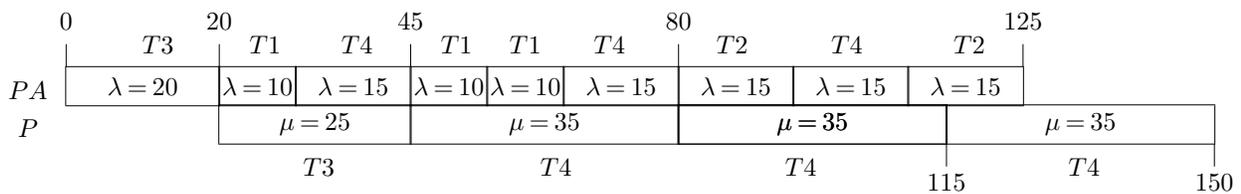
\begin{figure}[htbp]
\resizebox{\textwidth}{!}{
	{\normalsize
		\setlength{\unitlength}{1.15mm}
		\begin{picture}(160,24)(12,4) 
		%\thicklines
		
		\put(27,23){\makebox(0,0){$T3$}}
        \put(41,23){\makebox(0,0){$T1$}}
		\put(53,23){\makebox(0,0){$T4$}}
        \put(66,23){\makebox(0,0){$T1$}}
        \put(75,23){\makebox(0,0){$T1$}}
		\put(86,23){\makebox(0,0){$T4$}}
		\put(102,23){\makebox(0,0){$T2$}}
		\put(117,23){\makebox(0,0){$T4$}}
		\put(133,23){\makebox(0,0){$T2$}}
  
		\put(48,7){\makebox(0,0){$T3$}}
		\put(79,7){\makebox(0,0){$T4$}}
		\put(111,7){\makebox(0,0){$T4$}}
		\put(148,7){\makebox(0,0){$T4$}}

		\multiput(15,20)(27,0){1}{\line(0,1){3}}
		\put(15,26){\makebox(0,0){0}}
		
		\multiput(35,20)(27,0){1}{\line(0,1){3}}
		\put(35,26){\makebox(0,0){20}}
		
		\multiput(60,20)(27,0){1}{\line(0,1){3}}
		\put(60,26){\makebox(0,0){45}}
		
		\multiput(95,20)(27,0){1}{\line(0,1){3}}
		\put(95,26){\makebox(0,0){80}}
		
		\multiput(140,20)(27,0){1}{\line(0,1){3}}
		\put(140,26){\makebox(0,0){125}}
		
		\multiput(130,7)(27,0){1}{\line(0,1){3}}
		\put(130,5){\makebox(0,0){115}}
		\multiput(165,7)(27,0){1}{\line(0,1){3}}
		\put(165,5){\makebox(0,0){150}}

		%\multiput(130,15)(27,0){1}{\line(0,1){3}}
		\multiput(130,10)(27,0){1}{\line(0,1){3}}
		%\multiput(130,5)(27,0){1}{\line(0,1){3}}
		
		%======PA ====
		
		\put(15,15){\framebox(20,5){$\lambda=20$}}
        \put(35,15){\framebox(10,5){$\lambda=10$}}
  
		\put(45,15){\framebox(15,5){$\lambda=15$}}
        \put(60,15){\framebox(10,5){$\lambda=10$}}
        \put(70,15){\framebox(10,5){$\lambda=10$}}
  
		\put(80,15){\framebox(15,5){$\lambda=15$}}
        \put(95,15){\framebox(15,5){$\lambda=15$}}
		\put(110,15){\framebox(15,5){$\lambda=15$}}
        \put(125,15){\framebox(15,5){$\lambda=15$}}

		\put(10,17){\makebox(0,0){$PA$}}
		
		%======P==
		\put(35,10){\framebox(25,5){$\mu=25$}}
		\put(60,10){\framebox(35,5){$\mu=35$}}
		\put(95,10){\framebox(35,5){$\mu=35$}}
		
		\put(130,10){\framebox(35,5){$\mu=35$}}
		
		\put(95,10){\framebox(35,5){$\mu=35$}}

		\put(10,12){\makebox(0,0){$P$}}
		
		\end{picture}}}
		\vspace{-0.3in}	
	\caption{An illustration of Steps 6 and 7 of Algorithm~\ref{alg-B} for Example~1}
	\label{algB_ex3}	
    \vspace{-5mm}
\end{figure}

Algorithm~\ref{alg-MA} sorts $Q^+$ Group patients in a non-increasing order of their service times at stage~1, which creates waiting time due to greater service times at stage~2. In order to further reduce the patient waiting time while maintaining our no-idle time goal, one can alternate between $Q^+$ and $Q$ Group patients in the schedule. This approach reduces the waiting time by scheduling a $Q$ Group patient in stage~1 while a $Q^+$ Group patient is still busy at stage~2, instead of having another $Q^+$ Group patient wait for $P$ to complete their appointment with the previous patient. Motivated by this approach and aiming to reduce the patient waiting time while preserving a no-idle block schedule, we develop heuristic Algorithm~\ref{alg-B}. We later validate this reasoning during computational studies, based on the comparisons to the optimal solution. Figures~\ref{algB_ex1}, \ref{algB_ex2} and \ref{algB_ex3} illustrate the steps of Algorithm~\ref{alg-B} on Example~1. 

% In order to decrease patient waiting time without having any idle time for both $P$ and $PA$, we propose a (heuristic) algorithm as follows. 
\begin{algorithm}
\begin{algorithmic}
\caption{Improved No-idle Block Schedule Heuristic for a Single Block}
\label{alg-B}
 \Statex \textbf{Step 1}: Arrange $Q^+$ Group Patients such that $\lambda_{q+1} \geq \lambda_{q+2} \geq \ldots \geq \lambda_m$, where $r_i=1$, $i=1,\cdots,m$.
 \Statex \textbf{Step 2}: Schedule patients as no-wait in manner at earliest available time starting with all type ${q+1}$ followed by all type ${q+2}$, and so on until all $Q^+$ Group are scheduled (See Figure~\ref{algB_ex1} of Example 1). In this no-wait schedule, the finish time of each patient service at stage 1 ($PA$) is equal to the start time of this patient at stage 2 ($P$).
 \Statex \textbf{Step 3}: The partial schedule in Step 2 creates $(m-q-1)$ nonnegative idle time slots at $PA$  (See Figure~\ref{algB_ex1} of Example 1).
 \Statex \textbf{Step 4:} Arrange $Q$ Group Patients such that $\lambda_{1} \leq \lambda_{2} \leq \ldots \leq \lambda_q$.
 \Statex \textbf{Step 5:}  Schedule $Q$ Group Patients into $(m-q-1)$ idle time slots at stage 1  as follows: schedule each $Q$ Group patient in the given order in Step 4 in the first available idle time slot at the earliest time in that slot if the patient can be fit in that slot. If the patient is scheduled in that slot then update the remaining idle time of that slot. Continue this process until either all idle time slots at stage 1 are filled or there are no $Q$ Group Patients left with service time less than or equal to any of the remaining duration of idle time slots. (See Figure~\ref{algB_ex2} of Example 1.)   
 \Statex \textbf{Step 6:} After completing Step 5, positive idle time slots left at stage 1 (if any) are removed as follows: Number the idle time slots in PA from the earliest to the latest as slots $1, 2, \ldots, \ell$. To remove the idle time slot 1, advance the start times of all patients schedule in PA between slot 1 and 2 by the duration of the slot 1. Continue the process for the remaining slots $2,3 \ldots, \ell$ until all the slots are removed.
\Statex \textbf{Step 7:} If any $Q$ Group Patients are left after completing Step 5, these patients are scheduled at any order at the end of the partial schedule obtained at Step 6 (at earliest time at PA). (See Figure~\ref{algB_ex3} of Example 1.)
\end{algorithmic}
\end{algorithm}

The complexity of Algorithm~\ref{alg-B} can be determined as follows: Each of Steps 1, 2, and 3 can be determined in the order of at most $O(m \,log \, m)$ time as those steps involve sorting of patients to achieve a sequence. Step 4 can be determined in the order of at most  $O(m \,log \, m)$ time as it involves sorting of patients to achieve a sequence. Step 5 can be determined in the order of at most $O(m^2)$ time as it involves inserting $Q$ patients in at most $m$ idle time slots. Step 6 can be determined in the order of at most $O(m)$ time as it involves moving at most $r$ patients to close idle time gaps at $PA$. Thus the complexity of Algorithm~\ref{alg-B} is $O(m^2)$.
The following observation is made without proof as the final Step 6 of Algorithm~\ref{alg-B} eliminates the idle time at $PA$.

\textbf{Observation:} Algorithm~\ref{alg-B} provides a no-idle block schedule.

% \begin{lemma} \label{L-algB}
% Algorithm~\ref{alg-B} provides a no-idle block schedule.
% \end{lemma}

\begin{figure}[htbp]
\resizebox{\textwidth}{!}{
	{\normalsize
		\setlength{\unitlength}{1.15mm}
		\begin{picture}(160,24)(12,4) 
		%\thicklines
		
		\put(27,23){\makebox(0,0){$T3$}}
		
		\put(53,23){\makebox(0,0){$T4$}}
		
		\put(86,23){\makebox(0,0){$T4$}}
		
		\put(122,23){\makebox(0,0){$T4$}}

		\put(48,7){\makebox(0,0){$T3$}}
		\put(79,7){\makebox(0,0){$T4$}}
		\put(111,7){\makebox(0,0){$T4$}}
		\put(148,7){\makebox(0,0){$T4$}}

		\multiput(15,20)(27,0){1}{\line(0,1){3}}
		\put(15,26){\makebox(0,0){0}}
		
		\multiput(35,20)(27,0){1}{\line(0,1){3}}
		\put(35,26){\makebox(0,0){20}}
		
		\multiput(60,20)(27,0){1}{\line(0,1){3}}
		\put(60,26){\makebox(0,0){45}}
		
		\multiput(95,20)(27,0){1}{\line(0,1){3}}
		\put(95,26){\makebox(0,0){80}}
		
		\multiput(130,20)(27,0){1}{\line(0,1){3}}
		\put(130,26){\makebox(0,0){115}}

		\multiput(165,7)(27,0){1}{\line(0,1){3}}
		\put(165,5){\makebox(0,0){150}}

		\multiput(130,15)(27,0){1}{\line(0,1){3}}
		\multiput(130,10)(27,0){1}{\line(0,1){3}}
		%\multiput(130,5)(27,0){1}{\line(0,1){3}}
		
		%======PA ====
		
		\put(15,15){\framebox(20,5){$\lambda=20$}}
		\put(45,15){\framebox(15,5){$\lambda=15$}}
		\put(80,15){\framebox(15,5){$\lambda=15$}}
		\put(115,15){\framebox(15,5){$\lambda=15$}}

		\put(10,17){\makebox(0,0){$PA$}}
		
		%======P==
		\put(35,10){\framebox(25,5){$\mu=25$}}
		\put(60,10){\framebox(35,5){$\mu=35$}}
		\put(95,10){\framebox(35,5){$\mu=35$}}
		
		\put(130,10){\framebox(35,5){$\mu=35$}}
		
		\put(95,10){\framebox(35,5){$\mu=35$}}

		\put(10,12){\makebox(0,0){$P$}}
		
		\end{picture}}}	
		\vspace{-0.3in}	
	\caption{An illustration of Steps 1-3 of Algorithm~\ref{alg-B} for Example~1}
	\label{algB_ex1}
    \vspace{-5mm}
\end{figure}
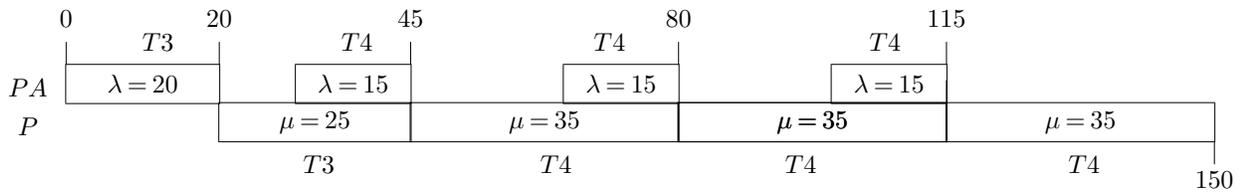

% \begin{algorithm}
% \begin{algorithmic}
% \caption{Algorithm $B$ ($AB$)}
%  \Statex \textbf{Step 1}: Arrange $Q^+$ Group Patients such that $\lambda_{q+1} \geq \lambda_{q+2} \geq \ldots \geq \lambda_m$ 
%  \Statex \hskip3.5em Consecutively schedule these patients in stage 2 such that $e_p^j=f_p^{j-1}$,
%  \Statex \hskip3.5em where there is no idle time for $P$, and
%  \Statex \hskip3.5em $f_a^{j}=e_p^j$, where the there is no wait time for patients in stage 2
%  \Statex \textbf{Step 2:} Arrange $Q$ Group Patients such that $\lambda_{1} \geq \lambda_{2} \geq \ldots \geq \lambda_q$
%  \Statex \hskip3.5em Until either the idle time for $PA$ before serving patient $j$, $d_a^j$ is filled or
%   \Statex \hskip3.5em there is no $Q$ Group patient left on the arranged list with $\lambda \geq d_a^j$
%   \Statex \hskip3.5em schedule $Q$ Group Patients into the idle times at stage 1,  
%   \Statex \hskip3.5em in the arranged order for the entire block
%   \Statex \textbf{Step 3:} If there are idle time for $PA$ at stage 1 after completing Step 2,
%   \Statex \hskip3.5em Move each patient ahead in the schedule until $e_a^j=f_a^{j-1}$ is achieved 
%   \Statex \hskip3.5em for every patient in the schedule
% \Statex \textbf{Step 4:} If there are any $Q$ Group Patients left after completing Step 2 for the entire block,  
% \Statex \hskip3.5em schedule these patients in any order at the end to form the back schedule.
% \end{algorithmic}
% \label{alg-B}
% \end{algorithm}

\section{Model and Analysis for Two-stage Scheduling for the Planning Horizon} \label{msmp}
\vspace{-1mm}

In this section, we extend the two-stage scheduling model introduced in Section~\ref{msm} from one block to the entire planning horizon to obtain a morning, afternoon, or entire day schedule. In addition to the parameters and variables used in the previous model, planning horizon specific notation is defined in Table~\ref{not2}. 

\vspace{-2mm}
\begin{table}[h]
\caption{ Notation for Planning Horizon Schedule}
\centering
    \resizebox{\textwidth}{!}{\begin{tabular}{|l l|}
         \hline
\multicolumn{2}{|l|}{\textbf{Parameters}}\\
\hline
$o_{p}$ & Unit cost for $P$ overtime.\\
$o_{a}$ & Unit cost for $PA$ overtime.\\
$\Gamma$ &  The block schedule of the planning horizon. $\Gamma =\{\pi_1, \pi_2, \ldots, \pi_k\}$ is defined as the concatenation of $k$ blocks.\\
\hline
\multicolumn{2}{|l|}{\textbf{Decision Variables}}\\
\hline
$x_{lct}$  & Binary variable that takes the value 1 if patient $l$ of block $c$ is assigned to the appointment slot\\ 
& $t$ and 0 otherwise, where $t=1,2,\cdots,n$ and the set of all patients in a block $l=1,2,\cdots,r$. \\
$y_t$ & Binary variable that takes the value 1 if the slot $t$ is assigned to a patient that has a nonzero\\
& service time at $P$ and 0 otherwise, where $t=1,2,\cdots,n$.\\
$\tau^{t}_a$ & The appointment time of the patient scheduled in position $t$.\\
\hline
\multicolumn{2}{|l|}{\textbf{Internal Variables}}\\
\hline
$e^t_a$ & The start time of patient in slot $t$ at stage 1 in $\Gamma$, where $t=1,2, \ldots, n$.\\
$e^t_p$ & The start time of patient in slot $t$ at stage 2 in $\Gamma$   where $t=1,2, \ldots, n$.\\
$f^t_a$ & The finish time of patient in slot $t$ at stage 1 in $\Gamma$  where $t=1,2, \ldots, n$.\\
$f^t_p$ & The finish time of patient in slot $t$ at stage 2 in $\Gamma$,  where $t=1,2, \ldots, n$.\\
$w^t_a$ & The waiting time of patient in slot $t$ at stage 1 in $\Gamma$, where $t=1,2, \ldots, n$.\\
$\bar{w}^t_p$ & The waiting time of patient in slot $t$ at stage 2 in $\Gamma$  (denoted by $w_p^t$), if the patient has nonzero\\
& service time at stage 2 (denoted by $\bar{w}^t_p$), where $t=1,2, \ldots, n$.\\
$d^t_a$ & The $PA$ idle time just before serving the patient in slot $t$ in $\Gamma$  at stage 1.\\
$d^t_p$ & The $P$ idle time just before serving the patient in slot $t$ in $\Gamma$  at stage 2.\\
$b_a$ & The $PA$ overtime\\
$b_p$ & The $P$ overtime\\
\hline
    \end{tabular}}

    \label{not2}
\end{table}

We formulate the Planning Horizon Model that minimizes the weighted sum of the cost of patient waiting times, the idle times, and overtimes of $P$ and $PA$, for the planning horizon schedule.

% \vspace{-3mm}

%t=(c-1)r+1,\cdots,cr\hspace{3mm}\forall c
%\newpage
\noindent\textbf{Planning Horizon Model}
%\alpha \sum_{t=1}^n(w^t_a + \bar{w}^t_p)+\beta_a \sum_{t=1}^n d^t_a+\beta_p \sum_{t=1}^n d^t_p+ o_a b_a + o_p b_p
\begin{align}
    \min\; & \alpha \sum_{t=1}^n w^t_a + \alpha \sum_{t=1}^n \bar{w}^t_p+\beta_a \sum_{t=1}^n d^t_a+\beta_p \sum_{t=1}^n d^t_p + o_ab_a+ o_pb_p \label{mm1}\\
    % \min\; &  U_p(\sigma)+ I_a(\sigma) + I_p(\sigma)+O_a(\sigma)+O_p(\sigma) \label{mm1}\\
    \text{subject to } & \sum_{t=(c-1)r+1}^{cr} x_{lct}=1\hspace{3mm} l=1,2,\cdots,r\hspace{3mm} c=1,2, \cdots,k \label{mm2}\\
    & \sum_{c=1}^{k}\sum_{l=1}^{r}  x_{lct}=1\hspace{3mm}   t=1,2,\cdots,n \label{mm3}\\
    & \sum_{s\in\{1,2,\cdots,q\}} x_{s11}=0\label{mm4}\\
    & e_a^{t+1}\geq e_a^t + \sum_{c=1}^{k}\sum_{l=1}^{r} x_{lct} \lambda_l \hspace{3mm}  t=1,2,\cdots,n-1\label{mm5}\\
    & e_p^t \geq e_a^t + \sum_{c=1}^{k}\sum_{l=1}^{r} x_{lct} \lambda_l \hspace{3mm} t=1,2,\cdots,n\label{mm6}\\
    & e_p^{t+1} \geq e_p^t + \sum_{c=1}^{k}\sum_{l=1}^{r} x_{lct} \mu_l \hspace{3mm} t=1,2,\cdots,n-1\label{mm7}\\
    & d_a^{t+1} \geq e_a^{t+1}-e_a^t- \sum_{c=1}^{k} \sum_{l=1}^r x_{lct} \lambda_l \hspace{3mm} t=1,2,\cdots,n-1\label{mm8}\\
    & d_p^{t+1} \geq e_p^{t+1}-e_p^t- \sum_{c=1}^{k} \sum_{l=1}^r x_{lct}\mu_l \hspace{3mm}  t=1,2,\cdots,r-1\label{mm9}\\
    & e_a^1=0\label{mm10}\\
    & \bar{w}_p^t \geq e_p^t-e_a^t-\sum_{c=1}^k\sum_{l=1}^r x_{lct} \lambda_l -M(1-y_t)\hspace{3mm} t=1,2,\cdots,n\label{mm12}\\
    & \frac{y_t}{M} \leq \sum_{c=1}^k \sum_{l=1}^r \mu_l x_{lct} \hspace{3mm} t=1,2,\cdots,n\label{mm13}\\
    & My_t \geq \sum_{c=1}^k \sum_{l=1}^r \mu_l x_{lct}\hspace{3mm} t=1,2,\cdots,n\label{mm14}\\
    & \tau_a^t \leq e_a^{t}\hspace{3mm}t=1,2,\cdots,n \label{mmtau1} \\
    & w_a^{t} \geq e_a^{t}-\tau_a^t \hspace{3mm}t=1,2,\cdots,n \label{mmtau2}\\
    & b_a \geq e_a^1 + \sum_{c=1}^k \sum_{t=(c-1)r+1}^{cr} \sum_{l=1}^r \lambda_l x_{lct} + \sum_{t=1}^n d_a^t - R\label{mm15}\\
    & b_p \geq e_p^1 + \sum_{c=1}^k \sum_{t=(c-1)r+1}^{cr} \sum_{l=1}^r \mu_l x_{lct} + \sum_{t=1}^n d_p^t - R\label{mm16}\\
    % & W_p(\sigma)= \sum_{t=1}^n \bar{w}^t_p\label{mm23}\\
    % & U_p(\sigma)=\alpha W_p(\sigma)\label{mm24}\\
    % & D_a(\sigma)= \sum_{t=1}^n d^t_a\label{mm25}\\
    % & I_a(\sigma)=\beta_a D_a(\sigma)\label{mm26}\\
    % & D_p(\sigma)= \sum_{t=1}^n d^t_p\label{mm27}\\
    % & I_p(\sigma)=\beta_p D_p(\sigma)\label{mm28}\\
    % & O_a(\sigma)=o_a b_a\label{mm29}\\
    % &  O_p(\sigma)=o_p b_p\label{mm30}\\
    & x_{lct}\in \{0,1\}\hspace{3mm} l=1,2,\cdots,r;\hspace{3mm} c=1,2,\cdots,k;\hspace{3mm} t=1,2,\cdots,n\label{mm17}\\
    & y_t\in \{0,1\}\hspace{3mm} t=1,2,\cdots,n\label{mm18}\\
    & d_a^t, d_p^t, e_a^t, e_p^t, w_a^t, \bar{w}_p^t,\tau_a^t, b_a, b_p \geq 0\hspace{3mm}   t=1,2,\cdots,n\label{mm19}
    %& M\text{: sufficiently large number}\label{mm20} \nonumber
\end{align}

% Different than the Single Block Model, the patient waiting time at stage~1 is not assumed to be zero as combining different blocks could introduce gaps due to their arrangement.

The objective function (\ref{mm1}) is the weighted sum of the total patient waiting time, idle time, and overtime costs of the $P$ and $PA$, respectively. Constraints (\ref{mm2}) through (\ref{mm14}) of the Planning Horizon Model are similar to constraints (\ref{m2}) through (\ref{m12}) of the Single Block Model. They ensure the same block structure, only extended to include the entire planning horizon with the notation provided in Table~\ref{not2}.

Constraint (\ref{mmtau1}) indicates that a patient cannot enter stage 1 before his or her appointment time and constraint (\ref{mmtau2}) defines the waiting time of each patient at stage 1 as the difference between the start and the appointment time. Constraint (\ref{mm15}) defines the overtime at stage 1 as the maximum of zero and the difference between the finish time of the last appointment and the regular time available for a day. Similarly, constraint (\ref{mm16}) defines the overtime at stage~2 as the maximum of zero and the difference between the finish time of the last appointment and the regular time available for a day. Constraints (\ref{mm17}) and (\ref{mm18}) define the binary variables and constraint (\ref{mm19}) ensures non-negativity of all of the remaining variables.

% Constraint (\ref{mm23}) defines the total patient waiting time at stage~2. Constraint (\ref{mm24}) denotes the total patient waiting cost over the planning horizon at stage~2. Constraints (\ref{mm25}) and (\ref{mm27}) define the total idle time of the $PA$ and $P$, respectively. Constraints (\ref{mm26}) and (\ref{mm28}) denote the total idle time cost of the $PA$ and $P$, respectively. Constraints (\ref{mm29}) and (\ref{mm30}) denote the overtime cost of the $PA$ and $P$, respectively.  

\subsection{Scheduling Multiple Blocks for a Planning Horizon}
\vspace{-1mm}

Repeating the blocks obtained by Algorithms~\ref{alg-MA} and \ref{alg-B} multiple ($k$) times over the planning horizon yields an appointment template that is simple while preserving the block structure. Thus, in order to design an appointment template that minimizes the weighted sum of patients waiting time, and the idle and overtimes of the $PA$ and $P$ we propose heuristic Algorithm~\ref{alg-RMA} for the entire planning horizon, that is based on repetition of Algorithm~\ref{alg-MA}. To further reduce the patient waiting time (as discussed in Section~\ref{msm}), we propose heuristic Algorithm~\ref{alg-RB} for the entire planning horizon, based on Algorithm~\ref{alg-B}. %\textcolor{red}{Different than the mathematical models developed in Sections~\ref{msm} and \ref{msmp}, the Algorithms~\ref{alg-RMA} and \ref{alg-RB} implicitly consider the cost parameters while designing appointment templates.} 

% This is achieved by preserving the priority structure of the provider related and patient related metrics imposed by the costs discussed in Section~\ref{cdp}.

\begin{algorithm}
\begin{algorithmic}
\caption{Two Stage Schedule for Planning Horizon}
\label{alg-RMA}
 \Statex \textbf{Step 1}: Run Algorithm~\ref{alg-MA} for one block and obtain the schedule $\pi$. The planning horizon has multiple blocks $\pi_c$,  $c=1,\cdots,k$, where $\pi_c=\pi$. $\Gamma$ is the planning horizon schedule, where $\Gamma=(\pi_1, \pi_2, \ldots, \pi_k)$.
 \Statex \textbf{Step 2:} Schedule $\pi_c$, $c=1,\cdots,k$, consecutively in $\Gamma$ such that $(e_p^1)_{c+1}=(f_p^r)_c$, where $c=1,\cdots,k-1$.
  \Statex \textbf{Step 3:} Repeat until each block is scheduled in $\Gamma$.
\end{algorithmic}
\end{algorithm}

\vspace{-5mm}

\begin{algorithm}
\begin{algorithmic}
\caption{Improved Two Stage Schedule for Planning Horizon}
\label{alg-RB}
 \Statex \textbf{Step 1}: Run Algorithm \ref{alg-B} for one block and obtain the schedule $\pi$. The planning horizon schedule is $\Gamma$, where $\Gamma=(\pi_1, \pi_2, \ldots, \pi_k)$, and $\pi_c=\pi$, $c=1,\cdots,k$.
 \Statex \textbf{Step 2:} Schedule $\pi_c$, $c=1,\cdots,k$,  consecutively in $\Gamma$ such that the start time of $P$ at the beginning of each block is equal to the finish time of $P$ at the end of the previous block. 
  \Statex \textbf{Step 3:} Repeat until each block is scheduled in $\Gamma$.
\end{algorithmic}
\end{algorithm}

Algorithms~\ref{alg-RMA} and \ref{alg-RB} yield appointment templates where each identical block obtained by Algorithms~\ref{alg-MA} and \ref{alg-B} are consecutively repeated $k$ times, respectively. Assuming $L_a \leq L_p$ the start time of $P$ in block $c+1$ is set equal to the finish time of $P$ in block $c$, where $c=1,\cdots,k-1$. The appointment templates obtained by Algorithm~\ref{alg-RMA} and \ref{alg-RB} provide no-idle time schedules for $P$. 
Figure~\ref{figx:block00} illustrates the Algorithm~\ref{alg-RMA} appointment template for Example~1 in Section~\ref{bs}, when $k=2$. We now develop an upper bound on the total patient waiting time in planning horizon schedule $\Gamma$ obtained by Algorithm~\ref{alg-RMA}.
The reader may refer to Appendix~\ref{app1.1a} for the proof of Lemma~\ref{Perf-bound-block-muti}.

\begin{lemma}\label{Perf-bound-block-muti}
Total patient waiting time $W$ in a no-idle $k$ repetitive block schedule $\Gamma$ obtained by Algorithm~\ref{alg-RMA} is bounded above by $\frac{k(r-v)(r-v-1)}{2}[\gamma_1 - \gamma_2] + {\frac{k(k-1)(r-v)}{2} }\theta$, where
	$k$ is the number blocks, $\theta =  \max \{0, \, [\,\sum_{j=1}^{r-v} \mu_{\sigma_{j}} - \sum_{j=1}^{r} \lambda_{\sigma_{j}}]\}$,  $\gamma_1 =\argmax_{1 \leq j \leq r-v}\{\mu_j\}$ and $\gamma_2 =\argmin_{2 \leq j \leq r-v}\{\lambda_j\}$,  This bound is tight.	
\end{lemma}

\noindent{\bf Remark:} Since Algorithm~\ref{alg-RB} provides a schedule by improving the wait times of the patients in the schedule obtained by Algorithm~\ref{alg-RMA}, the bound in Lemma~\ref{Perf-bound-block-muti} holds for a schedule obtained by Algorithm~\ref{alg-RB} also.

\looseness -1 The reader may refer to Appendix~\ref{app2} for the computational results of the heuristic algorithms and the Planning Horizon Model with deterministic service times. We find that Algorithm~\ref{alg-RMA} maintains $P$ zero idle time, with only small increases in $PA$ idle time and the overtime for both $P$ and $PA$. However, it greatly increases the patient wait time from roughly 30 minutes in the optimal Planning Horizon Model to 200 minutes. More importantly, we find that Algorithm~\ref{alg-RB} can maintain the overtime and idle time values in Algorithm~\ref{alg-RMA}, but reduce the patient wait time to only 62 minutes. These insights are similar to the ones presented in Section~\ref{compstu} for stochastic service times.
We also present different block patterns provided by the optimal solution and Algorithms~\ref{alg-RMA} and \ref{alg-RB} for deterministic service times in Table~\ref{res1.3} in Appendix~\ref{app2}. 
One can observe that most of the optimal block patterns alternate between $Q$ and $Q^+$ patient types in consecutive positions in the block. This shows that the optimal approach and the appointment template provided by Algorithm~\ref{alg-RB} use a similar approach. 

% In order to design an appointment template that further reduce the weighted sum of patients waiting time, and the idle and overtimes of the $PA$ and $P$ we propose a (heuristic) algorithm.

 % The parameters used in this computational study are discussed in Section~\ref{cdp}.

% \Statex \textbf{Step 1}: Run Algorithm \ref{alg-B} for one block and obtain the schedule $\pi$
%  \Statex \textbf{Step 2:} Schedule each $\pi_c$ consecutively in the appointment template $\sigma$
%  \Statex \hskip3.5em such that $(e_p^1)_{c+1}=(f_p^r)_c$, where $c=1,\cdots,k-1$.
%   \Statex \textbf{Step 3:} Repeat until each block is scheduled in the appointment template.

% Algorithm~\ref{alg-RB}, yields an appointment template where each identical block obtained by Algorithm~\ref{alg-B} is consecutively repeated $k$ times. Assuming $L_a \leq L_p$ the start time of $P$ in block $c+1$ is set equal to the finish time of $P$ in block $c$, where $c=1,\cdots,k-1$. The appointment template $\sigma$ obtained by Algorithm~\ref{alg-RB} provides a no-idle time schedule for $P$.

\section{Service Time Uncertainty} \label{ucd}
\vspace{-1mm}

In this section, we introduce service time uncertainty into two-stage scheduling by developing a stochastic programming model in Section~\ref{sto}. We computationally demonstrate that our heuristic algorithm is able to provide no-idle block schedules under practically relevant conditions in Section~\ref{uni}. In order to solve the stochastic programming model, we adopt the SAA approach in Section~\ref{saamodel}. 

% We study stochastic service times for $PA$ and $P$, .  

%  and discuss the expansion of the second stage from one $P$ to two $P$s

\subsection{Two-stage Scheduling Model for a Block with Stochastic Service Times} \label{sto}
\vspace{-1mm}

% In this subsection, we follow the approach proposed by \cite{keefer1983three}, where only three demand points are known with their probability of occurrences as shown in Table~\ref{3point}. The percentiles are the fractiles from the cumulative distribution function (CDF) of the random demand. The limited knowledge about the underlying distribution restricts us to parameterize the mean demand and demand variance. The distribution is discretized based on a few points. However, as \cite{keefer1983three} indicated, the underlying assumption in this method is that the points required for the approximation are ``elicited from the experts with perfect accuracy." 

% \begin{table}[]
%     \centering
%     \begin{tabular}{|c|c|c|}
%     \hline
%          Demand point $(d_i)$ & Fractile & Probability $(p_i)$  \\
%          \hline
%          $d_{1l}$ & 0.05 & 0.185\\
%          $d_{1}$ & 0.5 & 0.630\\
%          $d_{1u}$ & 0.95 & 0.185\\
%          \hline
%     \end{tabular}
%     \caption{Three-point approximation of continuous distribution}
%     \label{3point}
%     \vspace{-0.03in}
% \end{table}

For each block schedule, there exists $r$ time slots and $r$ patients. Thus, service times (at stage 1 and at stage 2) will have $N$ sample paths (scenarios), where $N$ is very large. 
We denote the objective function of each scenario $s$ with respect to the total patient waiting time, $PA$ idle time and $P$ idle time costs as $\Pi_s= \alpha\sum_{j=1}^r w^{j,s}_a+ \alpha \sum_{j=1}^r \bar{w}^{j,s}_p+\beta_a \sum_{j=1}^r d^{j,s}_a+\beta_p \sum_{j=1}^r d^{j,s}_p$. We define the occurrence probability of each scenario $s$ as $\rho_s$. We let $\lambda_{i,s}$ and $\mu_{i,s}$ be the value of the service times of patient $i$ at stages 1 and 2, respectively for a given scenario $s$, where $s=1,2,\cdots,N$. The formulation provides the expected cost of the objective function, which is $\Psi_s=\sum_{s=1}^N \rho_s \Pi_s(w_a^{j,s},\bar{w}_p^{j,s},d_a^{j,s},d_p^{j,s})$. The constraints of the Stochastic Programming Model are similar to those of the constraints of the Planning Horizon Model developed in Section~\ref{msmp}. The Stochastic Programming Model for a Single Block and the Planning Horizon can be found in Appendix~\ref{app3}. Note that the constraints of the Stochastic Programming Model are formed by extending the constraints of the Planning Horizon Model developed in Section~\ref{msmp} to account for each scenario. 

% below yields an optimal feasible solution with an objective value that does not deviate from a $\delta$ fraction of the optimal objective value $\Pi_s^*$ for any scenario.

% We then extend SAA approach for the LP model to the planning horizon (with multiple blocks). The reader may refer to E-Companion for SAA Model for the Planning Horizon Model. 

% Constraint (\ref{m30}) ensures that the objective value does not deviate from a $\delta$ fraction of the optimal objective value $\Pi_s^*$ for any scenario.
\vspace{-1mm}
\subsection{Scheduling with Stochastic Service Times} \label{uni}
\vspace{-1mm}

Suppose each patient type $i$ faces a uniform distribution of service time (as in similar studies such as \cite{denton2003sequential}, \cite{erdogan2015online}), $\tilde \lambda_{i}$ (respectively, $\tilde \mu_{i}$) drawn from the interval $[(1-w/2)\lambda_{i}, (1+w/2)\lambda_{i}]$ (respectively, $[(1-w/2)\mu_{i}, (1+w/2)\mu_{i}]$). Then, the resulting coefficient of variation is $\frac{w}{2 \sqrt{3}}$ for each service time, so that effectively $w$ represents the level of service time uncertainty faced by every patient type. We solve the resulting optimal block scheduling problem for three different values of demand uncertainty with $w=0.1$, $w=0.2$, $w=0.4$, i.e., respectively, 10\%; 20\%; 40\% deviations from the mean values.

Recall that the following attributes describe the problem instance $I$ encountered in our research:

\begin{itemize}
	
	\item Type $i$ patients scheduled, where $i=1, 2, \ldots, m$. 
	
	\item $r_i$ number of patients of type $i$ to be scheduled in each block. 
	
	\item $k$ number of blocks to be scheduled repetitively. 
	
	\item $Q$ Group Patients have service times $\lambda_{i}>0$ and $\mu_{i}= 0$, where $i=1,\cdots,q$.
	
	\item $Q^+$ Group Patients have $\lambda_{i}>0$, $\mu_{i} > 0$, where $\mu_{i} \geq \lambda_{i} >0$, and $i=q+1,\cdots,m$.
		
\end{itemize}

\medskip

% Recall that $\sum_{i=1}^{m} r_i =r$. Let $\sum_{i=1}^{q} r_i = v$. Algorithm~\ref{alg-MA} reorders $(r-v)$ patients of $Q^+$ Group Patients according to their mean service times such that $\lambda_{v+1} \geq \lambda_{v+2} \geq \ldots \geq \lambda_r$. Then it schedules patient $(v+1)$ at earliest available time, followed by patient $(v+2)$, and so on until it schedules patient $r$ (If there is tie $\lambda_{i} = \lambda_{i+1}$ the tie is broken in favor of scheduling patient with the smallest $\mu_i$), and finally, it schedules $v$ patients of $Q$ Group Patients at earliest available time at any order in 

Recall that Algorithm~\ref{alg-MA} provides a patient sequence $\sigma$ in block schedule $\pi$, where $\sigma=\{\sigma_{1}, \sigma_{2},$ $ \ldots,$ $ \sigma_{r-v}, \sigma_{r-v+1}, \ldots, \sigma_{r}\}$ and $\sigma_{j}$ denote the patient sequenced in the $j^{th}$ position in $\sigma$. Consider two partial sequences, $\sigma^f$ and $\sigma^b$ of $\sigma$, where $\sigma= \sigma^f \cup \sigma^b$:\\ (i) Front sequence, $\sigma^f =\{\sigma_{1}, \sigma_{2}, \ldots, \sigma_{r-v}\}= \{v+1, v+2, \ldots, r\}$, \\(ii) Back sequence, $\sigma^b=\{\sigma_{r-v+1}, \sigma_{r-v+2}, \ldots, \sigma_{r}\}= \{1, 2, \ldots, v\}$.

% \vspace{-0.1in}

The following result specifies the maximum uncertainty $w$ in the service times by characterizing the appointment times of patients $\tau$ under which Algorithm~\ref{alg-MA} always provides a block schedule with no idle time for $P$ in a stochastic service time environment. The reader may refer to Appendix~\ref{app1.1} for the proof of Lemma~\ref{L-block1_st}.

%\textcolor{blue}{Note that $r$ patients are scheduled in block $\pi$. Let $\sigma$ be a sequence of patients in $\pi$ given appointment at the clinic in a day, where $\sigma=\{\sigma_{1}, \sigma_{2}, \ldots, \sigma_{r}\}$  and $\sigma_{j}$ denote the patient scheduled in $j^{th}$ position (or slot) in $\sigma$. Recall that $\tau^{\sigma_j}_a$ denotes the appointment time of patient scheduled in $j^{th}$ position.}

\looseness -1 Note that under a deterministic time environment, Algorithm~\ref{alg-MA} provides a no idle schedule $\pi$ with a sequence of patients, $\sigma=\{\sigma_{1}, \sigma_{2}, \ldots, \sigma_{r}\}$ in $\pi$ for any problem instance $I$ (Lemma~\ref{L-block1}). Recall that $\tau^{\sigma_j}_a$ denotes the appointment time of patient scheduled in $j^{th}$ position in $\pi$. Here, we set $\tau^{\sigma_j}_a=e^j_a$, $\tau^{\sigma_1}_a=0$, $e^j_a=\sum_{\ell=1}^{j-1} \lambda_{\sigma_\ell}$, $j=2, 3, \ldots, r$,  and $e^j_a$ is the start time of $j^{th}$ patient at stage 1 in $\pi$.

\begin{lemma}\label{L-block1_st}
In service time uncertainty defined above by $w$ for a problem instance $I$, under the condition $w \leq \min_{1 \leq j \leq r-v-1} \left[\frac{2\sum_{\ell=1}^{j}(\mu_{\sigma_\ell} - \lambda_{\sigma_{\ell+1}})}{\sum_{\ell=1}^{i}(\mu_{\sigma_\ell} + \lambda_{\sigma_{\ell+1}})}\right]$, we obtain no a idle schedule $\pi$ by setting  new appointment times of patients $\tau^{\sigma_1}_a=0$ and $\tau^{\sigma_{j}} = 
	(1-\frac{w}{2}) \sum_{\ell=1}^{j-1} \lambda_{\sigma_\ell}$, $j=2, 3, \ldots, r$, in $\pi$.
\end{lemma}

In Example~1 (see Table~\ref{table-Pr1}), $\mu_{v+1}=25$ and $\lambda_{v+2}=15$. Thus,  $w \leq \frac{2(\mu_{v+1} - \lambda_{v+2})}{\mu_{v+1} + \lambda_{v+2}} =\frac{2 \times 10}{40}= 0.5$, i.e, 50\% deviation from the mean which is given by $i=1$ in the expression, $w \leq \min_{1 \leq i \leq r-v-1} \left[\frac{2\sum_{k=1}^{i}(\mu_{\sigma_k} - \lambda_{\sigma_{k+1}})}{\sum_{k=1}^{i}(\mu_{\sigma_k} + \lambda_{\sigma_{k+1}})}\right]$. Note that Algorithm~\ref{alg-MA} provides a block schedule with no $P$ idle time for Example~1 for all sample paths of random occurrence, if $w \leq 0.5$.

% \begin{corollary} \label{L-cor}
% Under service time uncertainty, Algorithm~\ref{alg-MA} always provides a block schedule with no idle time at $P$ for a problem instance $I$ with new appointment times of patients $\tau^{\sigma_1}_a=0$ and $\tau^{\sigma_{i}}_a = 
% (1-{w \over 2}) \sum_{k=1}^{i-1} \lambda_{\sigma_k}$, $i=2, 3, \ldots, r$, if $w \leq [{2(\mu_{\sigma_1} - \lambda_{\sigma_{2}}) \over \sum_{k=1}^{r-v-1}(\mu_{\sigma_k} + \lambda_{\sigma_{k+1}})}]$, where service times are specific to the problem instance $I$ under consideration.
% \end{corollary}

% For Example~1, $w \leq [{2(\mu_{\sigma_1} - \lambda_{\sigma_{2}}) \over \sum_{k=1}^{r-v}(\mu_{\sigma_k} + \lambda_{\sigma_{k}}) - \lambda_{\sigma_1} - \mu_{\sigma_{r-v}}}] = 1/7 = 0.1429$.
% That is block schedule with no idle time at $P$ for Example~1 for all sample paths of random occurrence, if $w \leq 0.1429$, i.e., 14.29\% deviation from the mean.

\begin{figure}[htbp]
\resizebox{\textwidth}{!}{
	{\normalsize
		\setlength{\unitlength}{1.2mm}
		\begin{picture}(160,24)(12,4) 
		%\thicklines
		
		\put(27,23){\makebox(0,0){$T3$}}
		
		\put(43,23){\makebox(0,0){$T4$}}
		
		\put(57,23){\makebox(0,0){$T4$}}
		
		\put(71,23){\makebox(0,0){$T4$}}
		
		\put(85,23){\makebox(0,0){$T2$}}
		
		\put(102,23){\makebox(0,0){$T2$}}
		
		\put(115,23){\makebox(0,0){$T1$}}
		\put(125,23){\makebox(0,0){$T1$}}
		\put(135,23){\makebox(0,0){$T1$}}
		
		\put(48,7){\makebox(0,0){$T3$}}
		\put(79,7){\makebox(0,0){$T4$}}
		\put(111,7){\makebox(0,0){$T4$}}
		\put(148,7){\makebox(0,0){$T4$}}

		\multiput(15,20)(27,0){1}{\line(0,1){3}}
		\put(15,26){\makebox(0,0){0}}
		
		\multiput(35,20)(27,0){1}{\line(0,1){3}}
		\put(35,26){\makebox(0,0){20}}
		
		\multiput(50,20)(27,0){1}{\line(0,1){3}}
		\put(50,26){\makebox(0,0){35}}
		
		\multiput(65,20)(27,0){1}{\line(0,1){3}}
		\put(65,26){\makebox(0,0){50}}
		
		\multiput(80,20)(27,0){1}{\line(0,1){3}}
		\put(80,26){\makebox(0,0){65}}
		
		\multiput(95,20)(27,0){1}{\line(0,1){3}}
		\put(95,26){\makebox(0,0){80}}
		
		\multiput(110,20)(27,0){1}{\line(0,1){3}}
		\put(110,26){\makebox(0,0){95}}

		\multiput(130,20)(27,0){1}{\line(0,1){3}}
		\put(130,26){\makebox(0,0){115}}

		\multiput(140,20)(27,0){1}{\line(0,1){3}}
		\put(140,26){\makebox(0,0){125}}

		\multiput(60,7)(27,0){1}{\line(0,1){3}}
		\put(60,5){\makebox(0,0){45}}
		
		\multiput(165,7)(27,0){1}{\line(0,1){3}}
		\put(165,5){\makebox(0,0){150}}

		\multiput(130,15)(27,0){1}{\line(0,1){3}}
		\multiput(130,10)(27,0){1}{\line(0,1){3}}
		%\multiput(130,5)(27,0){1}{\line(0,1){3}}
		
		%======PA ====
		
		\put(15,15){\framebox(20,5){$\lambda_{6}=20$}}
		\put(35,15){\framebox(15,5){$\lambda_{7}=15$}}
		\put(50,15){\framebox(15,5){$\lambda_{8}=15$}}
		\put(65,15){\framebox(15,5){$\lambda_{9}=15$}}
		\put(80,15){\framebox(15,5){$\lambda_{1}=15$}}
		\put(95,15){\framebox(15,5){$\lambda_{2}=15$}}
		
		\put(110,15){\framebox(10,5){\small $\lambda_{3}=10$}}
		\put(120,15){\framebox(10,5){\small $\lambda_{4}=10$}}
		\put(130,15){\framebox(10,5){\small $\lambda_{5}=10$}}

		\put(10,17){\makebox(0,0){$PA$}}
		
		%======P==
		\put(35,10){\framebox(25,5){$\mu_{6}=25$}}
		\put(60,10){\framebox(35,5){$\mu_{7}=35$}}
		\put(95,10){\framebox(35,5){$\mu_{8}=35$}}
		
		\put(130,10){\framebox(35,5){$\mu_{9}=35$}}

		\put(10,12){\makebox(0,0){$P$}}
		
		\end{picture}}}	
	%	\vspace{0.1in}	
	\caption{A Block Schedule $\pi_1$ obtained by Algorithm~\ref{alg-MA} for Example~1 for a sample path with $w=0$. The total patient waiting time is 90.}
	\label{figy:block1}	
    % \vspace{-8mm}
\end{figure}
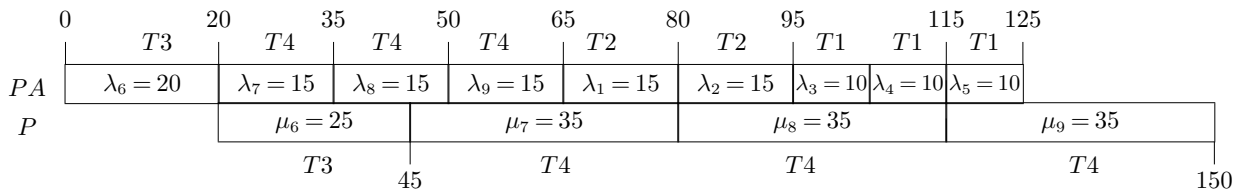

\noindent \textbf{Algorithm~\ref{alg-MA} Performance on Stochastic Service Times:} 
We computationally estimate $w$ for each of the practically relevant problem instances, based on Lemma~\ref{L-block1_st}. The estimated value of $w$ indicates the efficiency of Algorithm~\ref{alg-MA} (which is the single block version of Algorithm~\ref{alg-RMA}) to uncertain service times and shows the conditions under which it performs the best.

\looseness -1 We calculate $w$ for each of the 100 service time datasets randomly generated for each patient type in Section~\ref{cdp}. According to Lemma~\ref{L-block1_st}, Algorithm~\ref{alg-MA} provides a block schedule with no-idle time at $P$ if $w \leq 0.46$ based on the average calculated using 100 random instances. The minimum and maximum values of $w$ were calculated as 0.17 and 0.69, respectively. This means that in a real-world scenario where the deviation from the mean service time is less than 46$\%$ on average, Algorithm~\ref{alg-MA} accounts for service time uncertainty without compromising the no-idle time goal for $P$. 

\vspace{-2mm}
\subsection{Solution Procedure for Two-stage Scheduling Model for a Block with Stochastic Service Times} \label{saamodel}
\vspace{-1mm}

% This section introduces the structure of computational studies to evaluate the performance of Algorithm~\ref{alg-MA} under stochastic service times. 

% \subsubsection{Comparison with SAA Method} \label{SAA}

Solving the Stochastic Programming Model for stochastic service times yields exponential sample paths. Thus, following \cite{stauffer2021elasticity}, we reformulate the Stochastic Programming Model using a Sample Average Approximation (SAA) approach. The SAA Model assumes equal occurrence probability for each scenario and aims to minimize the average of a selected number of scenarios ($K$), that is the sample size, replicating each scenario set $K$, $\nu$ times. The number of replications or the sample size can be increased until achieving a desired precision level. By taking the average of all of the objective values provided by the SAA Model, we obtain an estimate of the (mean) optimal objective value of the Stochastic Programming Model. SAA Models and the detailed procedure can be found in Appendix~\ref{app3}.

% \begin{algorithm}
% \begin{algorithmic}
% \caption{SAA Solution Procedure \citep{law2007simulation}}
% \label{SAA-pr}
%  \Statex \textbf{Step 1}: Set $\nu_0=5$ replications (with each $K=10$) and set $\nu=\nu_0$.
%  \Statex \textbf{Step 2:} Calculate $\bar{\psi}_K^\nu(X)$ and $h(\nu,p)$ from  $\psi_K^1,\cdots,\psi_K^\nu$ with a $95\%$ confidence interval.
%   \Statex \textbf{Step 3:} If $\frac{h(\nu,p)}{\bar{\psi}_k^\nu(X)} < \frac{\xi}{1+\xi}$, then use $\bar{\psi}_k^\nu(X)$ as a point estimate for $\Psi_k^*$ and stop. Otherwise, replace $\nu$ with $\nu+1$ and proceed for another replication of the simulation and go to Step 2.
%   \Statex \textbf{Step 4:} If $\nu>\nu^{\max}$ increase $K$ sample paths by $K^\text{step}$, set $\nu=\nu_0=5$ and repeat Steps 2 and 3 above.
% \end{algorithmic}
% \end{algorithm}
After solving the SAA Model using the SAA Solution Procedure and obtaining a solution for the problem setting under service time uncertainty, we compare these results with the results from our heuristic algorithms and present our computational insights
in Section~\ref{compstu}.

% \begin{figure}[H] %t h b
% 	\centering
% 	\includegraphics[scale=1.1]{SAA_4BAMA_16nov_200_comp.jpeg}
% 	\caption{Algorithm~\ref{alg-RMA}, Algorithm~\ref{alg-RB} and SAA Model Results for randomized scenarios for Example~1, where patient waiting time cost, $\alpha\in\{0.2,\cdots,1\}$, the idle time costs for $PA$ and $P$, $\beta_a=\beta_p=1$, the overtime costs for $PA$ and $P$, $o_a=o_p\in\{1.2,1.5,1.8\}$, and regular time in the day, $R=200$ minutes}
% 	\label{SAA_BA_comp4}
% 	\vspace{-0.2in}
% \end{figure}

\section{Computational Study} \label{compstu}
\vspace{-1mm}

In this section, we compare the performance of the SAA solution approach and a first come first appointment (FCFA) approach from literature to our heuristic algorithms on randomly generated problem sets using real-world cost parameters. The parameters and computational setup are discussed in more detail in Section~\ref{cdp}. Then in Section~\ref{compin}, we illustrate the results and discuss related insights.

\subsection{Computational Parameters and Setup}\label{cdp}
\vspace{-1mm}

% This paper will utilize the data collected by \cite{oh2013guidelines} to address the research questions that will be further discussed in Section \ref{prd}. The data is collected over 9 work days from a three-provider family medicine practice in Massachusetts. The dataset consists of the list of prescheduled appointments and the revised list at the end of the day, which reflects the no-shows, cancellations, reschedules and same-day appointments. The dataset also includes a complete chronology of patient flow, from nurse service times to patient waiting times in the exam room or the lobby. Therefore, the data introduced by \cite{oh2013guidelines}, accommodates the characteristics that are necessary for the analysis that will be carried out in this study.

% $\$800$/hour

According to \cite{cech}, overall healthcare idle time cost per clinic can exceed $\$$10,000 per year in lost revenue. However, losing just one patient results in roughly $\$$1,600 per year in lost revenue for a clinic. %Additionally, studies show that patient waiting time is a significant factor that affects overall patient experience \citep{cech,dunnill2004medical}. 
\cite{klassen2019appointment} use similar cost differentials to test their healthcare simulation model. \cite{daggy2010using} assume patient waiting time cost is $\$0.33$/minute and clinic overtime cost is $\$13.33$/minute. \cite{zhou2021coordinated} set the relative cost parameters for patient waiting time, $P$ idle time and $P$ overtime as 0.2, 1, 1.5, respectively. Taking the $P$ idle time cost as the baseline and assuming that it is equal to 1, \cite{lee2018outpatient} utilize a range of values for the patient waiting time cost 0.1, 0.2,$\cdots$,1, and $P$ overtime cost 1.2, 1.5, 1.8. Thus, literature shows that patient waiting time is associated with a smaller cost compared to healthcare provider idle time cost, which is often valued less than overtime cost. For these computational results, we preserve the cost structure suggested by the literature and set the cost associated with overtime greater than or equal to the cost of healthcare provider idle time, which is then greater than patient waiting time.    
Following \cite{klassen2019appointment}, we do not assign different costs for $PA$ and $P$ for the sake of simplicity. However, this can easily be incorporated into the experimental structure.  In a real-world scenario, $P$ is expected to have higher overtime and idle time costs compared to $PA$, but this difference is not expected to impact the results in Section \ref{compstu}.

% In a few of the experiments we exceed this range and set the cost of idle time equal to 0 in order to illustrate cases where $P$s are less concerned with the idle time, since it can be filled by carrying out other tasks. 

\looseness -1 For all results in this section, the number of scenarios per replication was set to $K=15$. The initial number of replications ($\nu_0$) and the maximum number of replications ($\nu^{\max}$) were set to 5 and 10, respectively, which corresponds to a total of 75 and 150 sample paths. %SAA calculations were computed for an error rate of $\xi=0.04$. In this section each run was computed for a set of different cost parameters akin to \cite{lee2018outpatient} in order to expand the parameter space that we test the algorithms. 
For each combination of cost parameters, SAA was able to obtain the point estimate for the objective value after only 5 replications 100\% of the time with an error rate of $\xi=0.04$. All of the results presented in Section~\ref{compin} use the same set of scenarios ($K$) and replications ($\nu$) across the SAA Solution and all Algorithm values. That is, the same sample paths used for the SAA solution are used for Algorithm~3, Algorithm~4, and the FCFA results.

\looseness -1 Our computational results assume normally distributed random service times in each stage (as in similar studies such as \cite{robinson2003scheduling}, \cite{koeleman2012optimal}), but uniform distributions could easily be used depending on the particular clinic. There are six patient types with the mean and standard deviations obtained from the literature \citep{huang2015appointment,oh2013guidelines} as shown in Table~\ref{table-newparMAIN}. For this study, the length of regular time in a planning horizon, $R$, is set to 300 minutes and there are two blocks scheduled in each planning horizon. These two blocks together resemble the morning shift of an outpatient clinic, making up approximately five hours of planned clinical time. Depending on the operational structure of the clinic, it is possible to schedule another two blocks as an afternoon session to obtain a 10-hour workday.

\begin{table}[h]
	\caption{Stochastic Service Times Problem Setting, where $m=6$, $r=\sum_{i=1}^{m} r_i= 16$, and $q=3$, for one block}
    \vspace{-3mm}
	\begin{center}
		%\large
        \resizebox{\textwidth}{!}{
		\begin{tabular}{|c|c|c|c|c|c|} 
			\hline
			\multirow{2}{*}{Patient Type}	& \multicolumn{2}{|c|}{Stage 1 ($\lambda$)} & \multicolumn{2}{|c|}{Stage 2 ($\mu$)} & \multirow{2}{*}{Demand ratio ($r$)} \\
   \cline{2-5}
			                & Mean (min) & Sd (min)  & Mean (min) & Sd (min)  &              \\
                   \hline
High Complexity (HC) & 17.8 & 10.7 & 19.5 & 8.2 & 2            \\
Low Complexity (LC)  & 8.5  & 5.1  & 16.6 & 9   & 4            \\
Medium Complexity (MC)        & 9.5  & 6.1  & 12.7 & 7   & 4            \\
Low (L)             & 6    & 3    &   -  &  -  & 3            \\
Medium (M)          & 10   & 6    &  -   &  -  & 2            \\
High (H)            & 18   & 12   &  -   &  -  & 1  \\         
			\hline
		\end{tabular}}
	\end{center}
%	\vspace{0.1in}

	\label{table-newparMAIN}
    \vspace{-2mm}
\end{table}   

When establishing the block schedule templates for Algorithms~\ref{alg-RMA} and~\ref{alg-RB}, they are created using the mean service times of each patient type, exactly as they would be created in an actual clinic. Of course the actual random service times are used in all result calculations, which means the algorithm results represent realistically achievable results without knowledge of future actual service times. Besides eliminating the need for a commercial solver, both heuristic algorithms also provide a much more efficient solution than SAA by reducing the run time from over 50~hours to less than 2~seconds per scenario. 
Table~\ref{Algs_comptime} displays the computation times of the algorithms and SAA solution on the same scenario. This demonstrates a significant computational advantage while emphasizing the implementable and practical aspect of the algorithms.

\vspace{-2mm}
\begin{table}[h]
\caption{Computation times of the SAA Solution, Algorithm~\ref{alg-RMA}, Algorithm~\ref{alg-RB}, and FCFA, for a randomized secnario where regular time in the planning horizon $R=300$ minutes.}
\centering
% \resizebox{\textwidth}{!}{
\begin{tabular}{|c|c|}
\hline
 \multicolumn{2}{|c|}{Computation Times Per Randomized Scenario}\\
\hline
 Approach & Computation Time \\
\hline
SAA Solution & 50+ hours \\
Algorithm~\ref{alg-RMA} & 0.70 seconds\\
Algorithm~\ref{alg-RB} & 1.16 seconds\\
FCFA  & 0.21 seconds \\ 
\hline
\end{tabular}
% }
% \vspace{0.01in}
\label{Algs_comptime}
\vspace{-5mm}
\end{table}

\looseness -1 In addition to the comparison with the SAA approach, we compare the performance of Algorithms~\ref{alg-RMA} and \ref{alg-RB} with another heuristic algorithm from literature. In settings where the $PA$ sees some of the patients before $P$ or acts as a substitute for $P$, \cite{white2017ice} assume a First Come, First Appointment (FCFA) scheduling rule (which can also be referred to as first come, first served or first in, first out) for the patients, which is commonly used in other studies \citep{fan2020appointment}. Similar to \cite{white2017ice}, we implement an FCFA-based algorithm where the patients in a block are assumed to arrive in a randomized order and are scheduled accordingly. %The same set of random service times in scenarios ($K$) and replications ($\nu$) are used in the FCFA results as well as the SAA solution and all algorithm results.

\subsection{Computational Results and Insights} \label{compin}
\vspace{-1mm}

\begin{figure}[h] %t h b
	\centering
	\includegraphics[scale=0.75]{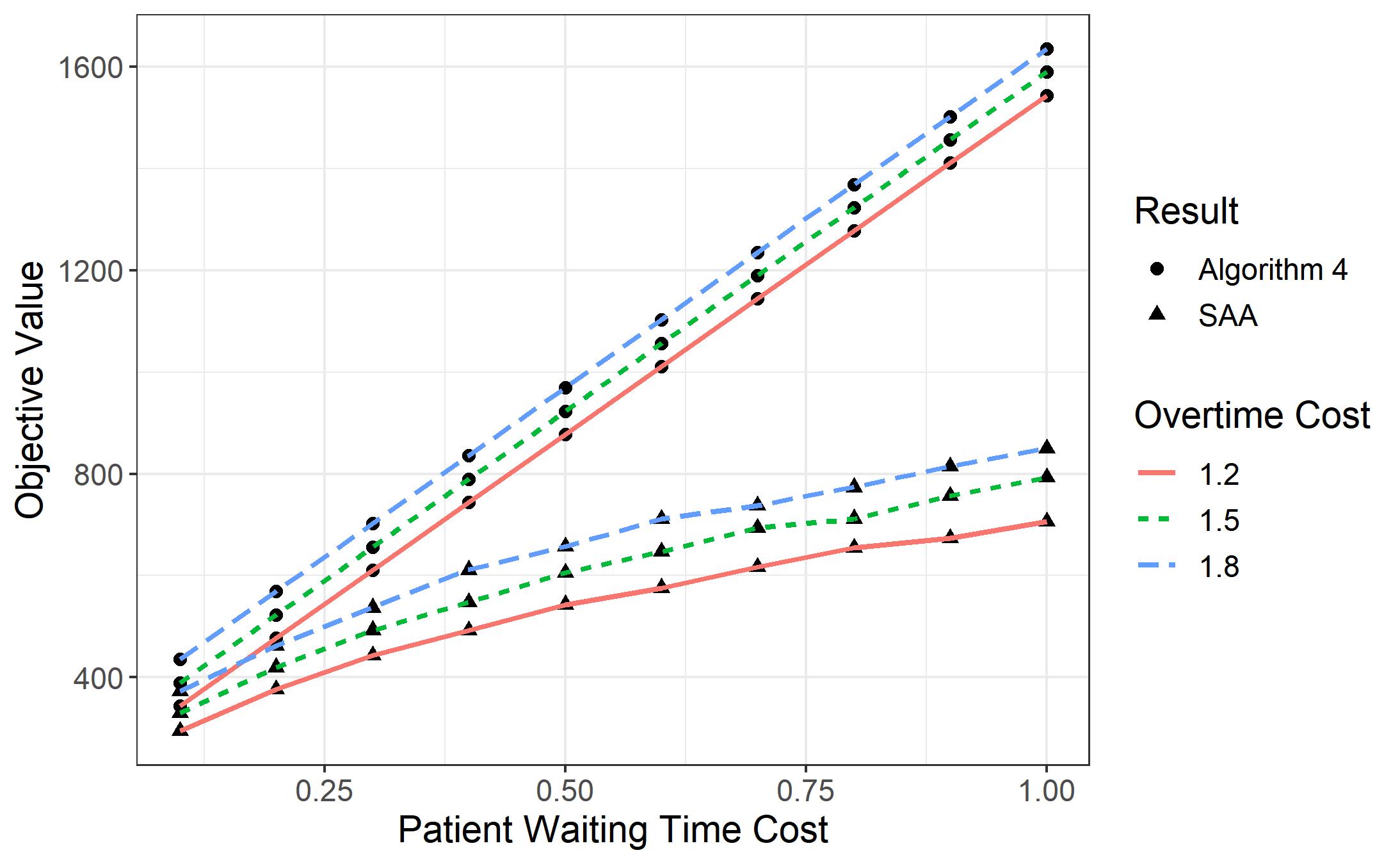}
 %{SAA_200_03Jan_comp.jpeg}
	\caption{Objective value results for Algorithm~\ref{alg-RB} and the solution using SAA for randomized scenarios, where patient waiting time cost, $\alpha\in\{0.1,\cdots,1\}$, the idle time costs for $PA$ and $P$, $\beta_a=\beta_p=1$, the overtime costs for $PA$ and $P$, $o_a=o_p\in\{1.2,1.5,1.8\}$, and regular time in the day, $R=300$ minutes}
	\label{SAA_BA_comp_July24}
	% \vspace{-0.2in}
\end{figure}

We begin by comparing the average objective function values for the solution using SAA with Algorithm~\ref{alg-RB} in Figure~\ref{SAA_BA_comp_July24}. This figure displays results for different overtime costs (1.2, 1.5, and 1.8) and varies the patient waiting time cost in the x-axis while keeping other parameters stable. Figure~\ref{SAA_BA_comp_July24} illustrates that both the solution using SAA and Algorithm~\ref{alg-RB} objective function values increase with respect to patient waiting time costs for all overtime costs. As the cost of patient waiting time decreases, the objective values of Algorithm~\ref{alg-RB} becomes closer to the solution using SAA. Thus, Algorithm~\ref{alg-RB} provides an implementable block schedule similar in cost to the values using SAA when clinics value patient waiting time less than provider idle time and overtime, which is the typical clinic preference. Algorithm~\ref{alg-RMA} objective value results (not graphed in Figure~\ref{SAA_BA_comp_July24}) would all be larger than the Algorithm~\ref{alg-RB} results.

% \textcolor{blue}{Table~\ref{Algs_SAA_July24} presents the objective function values of the solution using SAA, Algorithm~\ref{alg-RMA}, Algorithm~\ref{alg-RB} and average results of the FCFA based algorithm for a variety of cost parameters.  Algorithm~\ref{alg-RB} yields the smallest overall cost among other algorithm results, thus is the closest to the optimal solution using SAA. Algorithm~\ref{alg-RMA} performs better than the FCFA based algorithm for smaller values of the patient waiting time cost but the relationship is overturned for larger values.}

% \begin{table}[h]
% \centering
% \begin{tabular}{|c|c|c|c|c|}
% \hline
%  \multicolumn{5}{|c|}{ Objective Function Values for Planning Horizon Scheduling}\\
% \hline
% $\alpha$ & SAA & Algorithm~\ref{alg-RMA} & Algorithm~\ref{alg-RB} & FCFA (average)\\
% \hline
%  0.20 & 298.47 & 453.08  & 425.30 & 543.24  \\
% 0.40 & 351.18 & 621.31  & 563.82 & 657.94  \\
% 0.60 & 389.60 & 789.53  & 702.34 & 772.64  \\
% 0.80 & 435.02 & 957.76  & 840.86 & 887.34  \\
% 1.00 & 465.36 & 1125.99 & 979.38 & 1002.04 \\
% \hline
% \end{tabular}
% \caption{SAA and Algorithm results for randomized scenarios where idle time costs for $PA$ and $P$, $\beta_a=\beta_p=1$, the overtime costs for $PA$ and $P$, $o_a=o_p=1.5$, and regular time in the day \textcolor{blue}{$R=300$} minutes for a varying set of patient waiting time costs,~$\alpha$}
% \label{Algs_SAA_July24}
% \end{table}

Next Figure~\ref{fig:algs_comp} presents a comparison of Algorithm~3, Algorithm~4, and FCFA based on the objective function that is the weighted sum of the patient waiting time, idle time, and overtime costs of $PA$ and $P$ where the regular time in a planning horizon is 300 minutes and the idle time costs for $PA$ and $P$ are equal to 1. The figure displays the algorithm that yields the best (i.e. lowest) objective function value for a varying set of overtime costs $o_a, o_p$ and patient waiting time cost $\alpha$, extended to $\{0.9, \cdots, 2.1\}$ and $\{0.1,\cdots,0.8\}$, respectively. %Objective function values were calculated based on Table~\ref{Algs_overall_July24} results for Algorithms~\ref{alg-RMA} and \ref{alg-RB} and the average FCFA metrics. 
We observe that Algorithm~\ref{alg-RB} provides the lowest cost appointment template in 60\% of the cases analyzed. However, for smaller values of overtime costs and larger values of patient waiting time costs, the FCFA-based algorithm provides a lower overall cost. Note that in the portion of the cases where FCFA is favorable to Algorithm~\ref{alg-RB}, the cost parameters are closer to each other. This shows that when the cost of provider idle time, provider overtime, and patient waiting time are valued relatively equally, then random scheduling could be fine for an outpatient clinic. This assumption, however, contradicts many real-life settings (as found in other literature discussed at the beginning of Section~\ref{cdp}) where the clinics prioritize the efficiency (idle time and overtime) of the healthcare providers, which Algorithm~\ref{alg-RB} is better able to capture. 
To find cases where Algorithm~\ref{alg-RMA} is favorable, we would need to expand the cost parameters beyond Figure~\ref{fig:algs_comp}. We observed that for larger values of $o_a,o_p$ and $\beta_a,\beta_p$ beyond 10, and with $\alpha=0.1$, Algorithm~\ref{alg-RMA} yields the schedule with the lowest total cost. %Also note that the solution using SAA is always the lowest, but is a lot more complex and requires much more time to run making it impractical for most real-world settings.

\begin{figure}[th]%htbp
\vspace{22 mm}
	{\normalsize
		\setlength{\unitlength}{1.1mm}
		\begin{picture}(100,30)(0,15) 
		\thicklines
  	\put(25,55){\makebox(0,0){$\alpha$}}
  %      \put(25,60){\makebox(0,0){1.0}}		
%		\put(25,55){\makebox(0,0){0.9}}
		\put(25,50){\makebox(0,0){0.8}}
        \put(25,45){\makebox(0,0){0.7}}		
		\put(25,40){\makebox(0,0){0.6}}
		\put(25,35){\makebox(0,0){0.5}}
        \put(25,30){\makebox(0,0){0.4}}		
		\put(25,25){\makebox(0,0){0.3}}
		\put(25,20){\makebox(0,0){0.2}}
        \put(25,15){\makebox(0,0){0.1}}

	\put(112,7){\makebox(0,0){$o_a, o_p$}}
		\put(97,7){\makebox(0,0){2.1}}
		\put(82,7){\makebox(0,0){1.8}}
        \put(67,7){\makebox(0,0){1.5}}		
		\put(52,7){\makebox(0,0){1.2}}
		\put(37,7){\makebox(0,0){0.9}}

		\multiput(30,10)(27,0){1}{\vector(0,1){50}}
		\multiput(105,10)(27,0){1}{\line(0,1){43}}
		% \put(20,26){\makebox(0,0){0}}
		\multiput(30,10)(27,0){1}{\vector(1,0){90}}
        \multiput(30,53)(27,0){1}{\line(1,0){75}}

\multiput(45,22)(27,0){1}{\line(0,1){5}}
% \multiput(60,37)(27,0){1}{\line(0,1){5}}
\multiput(75,27)(27,0){1}{\line(0,1){5}}

\multiput(30,22)(27,0){1}{\line(1,0){15}}
\multiput(45,27)(27,0){1}{\line(1,0){15}}
\multiput(60,27)(27,0){1}{\line(1,0){15}}
\multiput(75,32)(27,0){1}{\line(1,0){15}}
\multiput(90,32)(27,0){1}{\line(1,0){15}}

		%===========Col2=============

%  \put(30,57){\makebox(15,5){FCFA}}		
%  \put(30,52){\makebox(15,5){FCFA}}		
  \put(30,47){\makebox(15,5){FCFA}}
  \put(30,42){\makebox(15,5){FCFA}}
  \put(30,37){\makebox(15,5){FCFA}}		
  \put(30,32){\makebox(15,5){FCFA}}
  \put(30,27){\makebox(15,5){FCFA}}		
  \put(30,22){\makebox(15,5){FCFA}}		
  \put(30,17){\makebox(15,5){Alg~\ref{alg-RB}}}
  \put(30,12){\makebox(15,5){Alg~\ref{alg-RB}}}
		
		%===========Col3=============

 % \put(45,57){\makebox(15,5){FCFA}}		
 % \put(45,52){\makebox(15,5){FCFA}}		
  \put(45,47){\makebox(15,5){FCFA}}
  \put(45,42){\makebox(15,5){FCFA}}
  \put(45,37){\makebox(15,5){FCFA}}		
  \put(45,32){\makebox(15,5){FCFA}}
  \put(45,27){\makebox(15,5){FCFA}}		
  \put(45,22){\makebox(15,5){Alg~\ref{alg-RB}}}		
  \put(45,17){\makebox(15,5){Alg~\ref{alg-RB}}}
  \put(45,12){\makebox(15,5){Alg~\ref{alg-RB}}}		
		
		%======Col4 ====

  %\put(60,57){\makebox(15,5){FCFA}}		
  %\put(60,52){\makebox(15,5){FCFA}}		
  \put(60,47){\makebox(15,5){FCFA}}
  \put(60,42){\makebox(15,5){FCFA}}
  \put(60,37){\makebox(15,5){FCFA}}		
  \put(60,32){\makebox(15,5){FCFA}}
  \put(60,27){\makebox(15,5){FCFA}}		
  \put(60,22){\makebox(15,5){Alg~\ref{alg-RB}}}		
  \put(60,17){\makebox(15,5){Alg~\ref{alg-RB}}}
  \put(60,12){\makebox(15,5){Alg~\ref{alg-RB}}}		
		
		%===========Col5=============
 
 % \put(75,57){\makebox(15,5){FCFA}}		
 % \put(75,52){\makebox(15,5){FCFA}}		
  \put(75,47){\makebox(15,5){FCFA}}
  \put(75,42){\makebox(15,5){FCFA}}
  \put(75,37){\makebox(15,5){FCFA}}		
  \put(75,32){\makebox(15,5){FCFA}}
  \put(75,27){\makebox(15,5){Alg~\ref{alg-RB}}}		
  \put(75,22){\makebox(15,5){Alg~\ref{alg-RB}}}		
  \put(75,17){\makebox(15,5){Alg~\ref{alg-RB}}}
  \put(75,12){\makebox(15,5){Alg~\ref{alg-RB}}}		
		
		%===========Col6=============
 
 % \put(90,57){\makebox(15,5){FCFA}}		
%  \put(90,52){\makebox(15,5){FCFA}}		
  \put(90,47){\makebox(15,5){FCFA}}
  \put(90,42){\makebox(15,5){FCFA}}
  \put(90,37){\makebox(15,5){FCFA}}		
  \put(90,32){\makebox(15,5){FCFA}}
  \put(90,27){\makebox(15,5){Alg~\ref{alg-RB}}}		
  \put(90,22){\makebox(15,5){Alg~\ref{alg-RB}}}		
  \put(90,17){\makebox(15,5){Alg~\ref{alg-RB}}}
  \put(90,12){\makebox(15,5){Alg~\ref{alg-RB}}}		
	
		\end{picture}}	
	\vspace{0.4in}	
	\caption{Lowest objective value between Algorithm~\ref{alg-RB}, Algorithm~\ref{alg-RMA}, and FCFA where idle time cost for $PA$ and $P$, $\beta_a= \beta_p =1$, and regular time in the day $R=300$ minutes for a varying set of patient waiting time costs $\alpha$ and overtime costs for $PA$ and $P$, $o_a=o_p$}
	\label{fig:algs_comp}	
    \vspace{-3mm}
\end{figure}
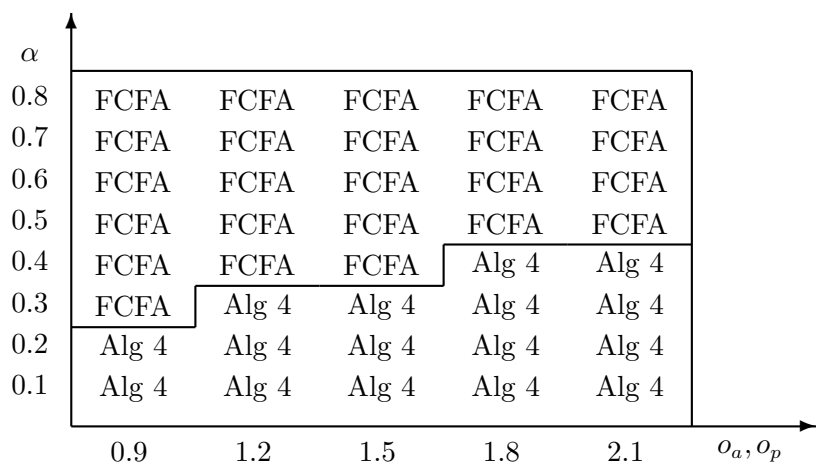 %March2025

\looseness -1 These computational results suggest that in many real-world scenarios where the overtime and idle time costs of the healthcare providers are associated with greater importance compared to the patient waiting time, Algorithm~\ref{alg-RB} using a block schedule is able to preserve this priority structure. Thus, repeating the same block schedule throughout the planning horizon smooths out the workload and provides a simple yet effective appointment template that is easy to follow for outpatient clinics.

% \vspace{-2mm}
\begin{table}[h]
\caption{SAA solution overtime, idle time, and wait time results (in minutes) for randomized scenarios where idle time costs for $PA$ and $P$, $\beta_a=\beta_p=1$, and regular time in the planning horizon $R=300$ minutes for a varying set of patient waiting time costs,~$\alpha$.}
\centering
\resizebox{\textwidth}{!}{
\begin{tabular}{|c|c|c|c|c|c|c|}
\hline
 \multicolumn{7}{|c|}{The Optimal Planning Horizon Solution using SAA}\\
\hline
\multirow{2}{*}{$o_a,o_p$} & \multirow{2}{*}{$\alpha$} & $PA$ Overtime & $P$ Overtime & $PA$ Idle Time & $P$ Idle Time & Waiting Time\\
& & (minutes) & (minutes) & (minutes) & (minutes) & (minutes)\\
 \hline
 \parbox[t]{2mm}{\multirow{3}{*}{\rotatebox[origin=c]{90}{1.2}}} &
0.20 & 50.57  & 114.66 & 26.69 & 50.57 & 663.18 \\
& 0.60 & 77.75  & 128.65 & 53.86 & 35.25 & 398.57 \\
& 0.80 & 101.71 & 147.03 & 77.83 & 52.73 & 282.32 \\
\hline
 \parbox[t]{2mm}{\multirow{3}{*}{\rotatebox[origin=c]{90}{1.5}}}
& 0.20 & 40.37 & 108.57 & 16.48 & 14.74 & 825.88 \\
& 0.60 & 72.16 & 124.08 & 48.27 & 30.76 & 456.45 \\
& 0.80 & 91.11 & 136.65 & 67.22 & 42.69 & 325.33 \\
\hline
\parbox[t]{2mm}{\multirow{3}{*}{\rotatebox[origin=c]{90}{1.8}}}
& 0.20 & 38.47 & 107.99 & 14.47 & 13.89 & 852.98 \\
& 0.60 & 61.59 & 117.00 & 37.70 & 23.51 & 547.51 \\
& 0.80 & 78.17 & 128.04 & 54.28 & 34.11 & 392.52\\
\hline
\end{tabular}}
% \vspace{0.01in}

\label{SAA_overall_July24}
% \vspace{-3mm}
\end{table}

\begin{table}[h]
\caption{Overtime, idle time, and wait time results (in minutes) for Algorithm~\ref{alg-RMA}, Algorithm~\ref{alg-RB}, and FCFA for randomized scenarios where regular time in the planning horizon $R=300$ minutes.}
\centering
\resizebox{\textwidth}{!}{
\begin{tabular}{|c|c|c|c|c|c|}
\hline
 \multicolumn{6}{|c|}{Algorithm Results for Planning Horizon Scheduling}\\
\hline
 \multirow{2}{*}{Approach} & $PA$ Overtime & $P$ Overtime & $PA$ Idle Time & $P$ Idle Time & Waiting Time\\
 & (minutes) & (minutes) & (minutes) & (minutes) & (minutes)\\
\hline
Algorithm~\ref{alg-RMA} & 41.87 & 109.71 & 17.98 & 7.30 & 1477.16\\
Algorithm~\ref{alg-RB} & 42.18 & 110.74 & 18.30 & 8.33 & 1333.69 \\
FCFA  & 58.52  & 136.77 & 34.63 & 43.19 & 1025.94 \\ 
\hline
\end{tabular}}
% \vspace{0.01in}
\label{Algs_overall_July24}
\vspace{-5mm}
\end{table}

Tables~\ref{SAA_overall_July24} and \ref{Algs_overall_July24} display the underlying time results (in minutes) for overtime, idle time, and patient waiting time. We observe that as the overtime costs increase both overtime and idle times decrease while patient waiting time increases. Conversely, when the waiting time cost ($\alpha$) increases, $P$ and $PA$ overtime and idle times increase while the patient waiting time decreases. Since the algorithms do not change their schedule based on cost parameters, the minutes of overtime, idle time, and wait time do not change with variation in these parameters. Thus, Table~\ref{Algs_overall_July24} only displays one time value per algorithm. For the algorithms, only the objective functions change with varying cost parameters (as seen in Figures~\ref{SAA_BA_comp_July24} and~\ref{fig:algs_comp}). 

% and thus the same idle time, overtime, and waiting time minutes occur for varying cost parameters 

%Table~\ref{Algs_overall_July24} presents a row for the best, worst, and average FCFA times because we used 5 different random arrival orders of patients (with the same random service times used in the other algorithms and SAA) to generate these results. The average solution is used for comparison purposes and in all objective function results.

From comparing the two tables, we see that Algorithm~\ref{alg-RB} provides a smaller $P$ idle time from the solution using SAA in all cost combinations explored and smaller $PA$ idle times for larger values of overtime costs, but results in more patient waiting time. While the patient waiting time appears to be a large increase from 825 minutes with SAA (when overtime is 1.5 and patient wait time is $\alpha = 0.2)$ to 1333 minutes with Algorithm~\ref{alg-RB}, these values are the total for all patients. This is only an increase of 16 minutes per patient, in reality.

Table~\ref{Algs_overall_July24} shows that Algorithm~\ref{alg-RMA} provides the smallest $P$ idle time among all approaches. However, with less than 1.5 minutes increase in both $PA$ and $P$ overtime and idle times, Algorithm~\ref{alg-RB} is able to reduce the patient waiting time incurred by Algorithm~\ref{alg-RMA} by approximately 143 minutes. This reduced patient waiting time for similar idle time and overtime is why Algorithm~\ref{alg-RB} dominates Algorithm~\ref{alg-RMA} in Figures~\ref{SAA_BA_comp_July24} and~\ref{fig:algs_comp}. This is also similar to the deterministic service time setting results discussed in Appendix~\ref{app2}.

\vspace{-2mm}
\section{Conclusion and Future Research Directions} \label{con}
\vspace{-1mm}

\looseness -1 This paper develops an outpatient appointment scheduling system for two-stage outpatient clinics under patient heterogeneity and service time uncertainty. We consider a two-stage outpatient clinic with one physician assistant ($PA$) in stage~1 and one physician ($P$) in stage~2. A block scheduling approach similar to that of \cite{lee2018outpatient} is developed for the two-stage scheduling setting which is based on the idea of production smoothing used in the Toyota Production System with the objective of balancing the workload in each workstation across a planning horizon. Our objective is to find appointment schedules that minimize a weighted sum of the healthcare providers’ idle time and overtime and patient waiting time.

First, for clinics facing negligible service time uncertainty (i.e., with deterministic service times), we develop an efficient scheduling algorithm to obtain a block schedule with zero total  system idle time while minimizing total patient waiting time. However, minimizing patient waiting time in a block schedule with zero system idle time is shown to be strongly $\mathcal{NP}$-Hard. We extend the idea of our efficient scheduling algorithm to multiple blocks and compare its performance with that of an exact method developed using mathematical programming. We find that our heuristic algorithm provides an appointment template with similar $P$ idle time and overtime as the mathematical program. %We also illustrate that our easily implementable algorithm block schedule preserves the priority structure where idle time and overtime of the healthcare providers are associated with greater costs compared to patient waiting time. 

Second, the problem setting is expanded to include stochastic patient service times when seeing healthcare providers. We then develop a stochastic programming model, solved using the Sample Average Approximation (SAA) approach. We show that our heuristic algorithm yields a block schedule with no-idle time for $P$, when the deviation from the mean service time is less than 46\% on average. %We obtain estimates for the Stochastic Programming Model using the SAA approach for the cost parameters chosen in compliance with the literature with a relative error of 4\%. 
Setting the complex SAA solution as the benchmark, we observe that our heuristic Algorithm~\ref{alg-RB} provides an appointment template with similar objective function costs when the patient wait time cost is low. We also show that Algorithm~\ref{alg-RB} will minimize $P$ idle time with only a 16 minute/patient increase in patient wait time. In addition, we compare the performance of our algorithms with a first come, first appointment (FCFA) scheduling rule commonly found in similar studies. We show that in scenarios where FCFA is favorable to Algorithm~\ref{alg-RB}, the idle time, overtime and patient wait time cost parameters are relatively similar to each other. However, this contradicts many real-life settings where clinics prioritize the efficiency (idle times and overtimes) of the healthcare providers, which Algorithms~\ref{alg-RMA} and \ref{alg-RB} are better able to capture. We contribute to the healthcare scheduling literature by providing insights into both deterministic and stochastic service time environments. 

% Finally, we examine the impact of patient no-shows. We propose several overbooking strategies to deal with no-shows in a block schedule and show that Algorithm~\ref{alg-B} performs well in a no-show setting.  

Further research directions can include introducing a two-stage outpatient clinical setting where the patient type or the need for further treatment by $P$ may change after the initial stage of $PA$ service is completed. Finally, an extended study of block schedule overbooking strategies with no-shows could provide interesting and practical insights. 

%Another promising extension of the current study is to develop appointment templates with two $PA$s in the first stage and two $P$s in the second stage of the outpatient clinic. Increasing the number of healthcare providers would increase the complexity of the problem and create additional challenges for appointment scheduling coordination. Finally, an extended study of block schedule overbooking strategies with no-shows could provide interesting and practical insights.

% \section*{Acknowledgement(s)}

% An unnumbered section, e.g.\ \verb"\section*{Acknowledgements}", may be used for thanks, etc.\ if required and included \emph{in the non-anonymous version} before any Notes or References.

% \section*{Disclosure statement}

% An unnumbered section, e.g.\ \verb"\section*{Disclosure statement}", may be used to declare any potential conflict of interest and included \emph{in the non-anonymous version} before any Notes or References, after any Acknowledgements and before any Funding information.

% \section*{Notes on contributor(s)}

% An unnumbered section, e.g.\ \verb"\section*{Notes on contributors}", may be included \emph{in the non-anonymous version} if required. A photograph may be added if requested.

% \section*{Nomenclature/Notation}

% An unnumbered section, e.g.\ \verb"\section*{Nomenclature}" (or \verb"\section*{Notation}"), may be included if required, before any Notes or References.

% \section*{Notes}

% An unnumbered `Notes' section may be included before the References (if using the \verb"endnotes" package, use the command \verb"\theendnotes" where the notes are to appear, instead of creating a \verb"\section*").

% \bibliographystyle{nonumber}

 % \bibliographystyle{ormsv080}
    % \bibliographystyle{apacite}
\bibliographystyle{chicago}

 \let\oldbibliography\thebibliography
 \renewcommand{\thebibliography}[1]{%
 	\oldbibliography{#1}%
 	\baselineskip11pt %Change this for line spacing within the same reference
 	\setlength{\itemsep}{7pt}% %Change this for spacing between two referneces
 }
    
    \bibliography{socialref2}

%% Here starts the e-companion (EC)
%%%%%%%%%%%%%%%%%%%%%%%%%%%%%%%%%%%%%%%%%%%%%%%%%%%%%%%%%%
 \ECSwitch

%\ECDisclaimer
%%%%%%%%%%%%%%%%%%%%%%%%%%%%%%%%%%%%%%%%%%%%%%%%%%%%%%%%%%

%%% Main head for the e-companion
% \ECHead{Appendix}

\begin{APPENDICES}

\section{Scheduling Excessive $Q$ Group Patients}\label{app0}

If $L_{a} > L_{p}$ then perform Algorithm~\ref{alg-bal} to remove excessive $Q$ Group Patients from each block so that workload $L_{a}$ is less than or equal to $L_{p}$. All the removed patients will form the last $(k+1)^{th}$ block. The idea behind Algorithm~\ref{alg-bal} is that it provides a list $V$ that consists of excessive $Q$ Group Patients, which forms the $(k+1)^{th}$ block while it balances the workloads $L_{a}$ and $L_{p}$ to be nearly equal for each of $k$ blocks. Table~\ref{table-Pr2} presents Example~2, where $m=4$, $q=2$, and $r=13$, which illustrates the formation of a block schedule with excessive $Q$ Group Patients. In Example~2, $L_{a} =\sum_{i=1}^{m} r_i\lambda_{i} = 180$ and $L_{p} =\sum_{i=1}^{m} r_i\mu_{i}= 130$, initially. Since $L_{a} > L_{p}$, we apply Algorithm~\ref{alg-bal} and obtain list $V$ which consists of one of type 1 and three of type 2 patients from each of $k$ blocks forming the last $(k+1)^{th}$ block ($\pi_3$ block in Example~2) as illustrated in Figure~\ref{figx:block3}. The new revised workloads for each of the $k$ blocks are $L_{a} = 125$ and  $L_{p} = 130$, as illustrated in Figure~\ref{figx:block1}.

\begin{algorithm}
\begin{algorithmic}
\caption{Balance-Workload Procedure}
\label{alg-bal}
 \Statex \textbf{Step 0}: Start with an empty list $V$.
 \Statex \hskip3.5em Set $u_i=r_i$ for patient type $i$, $i=1, 2, \ldots, q$.
 \Statex \hskip3.5em Note that $L_{a} =\sum_{i=1}^{m} r_i\lambda_{i} > L_{p} =\sum_{i=1}^{m} r_i\mu_{i}$.
 \Statex \textbf{Step 1:} Find index $i$ with the largest $u_i$ and the tie is broken for high $\lambda_{i}$.
 \Statex \hskip3.5em Remove one patient of type $i$ from each block and place it in
 \Statex \hskip3.5em the remove list $V$, and set $u_i=u_i-1$. Revise the new workloads, $L_{a}$, $L_{p}$.
  \Statex \hskip3.5em If $L_{a} > L_{p}$ then repeat Step 1; otherwise go to Step 2.
  \Statex \textbf{Step 2}: Form the last block with patients in the list $V$.
\end{algorithmic}
\end{algorithm}

\begin{table} [h]
	\caption{Example 2, where $m=4$, $r=\sum_{i=1}^{m} r_i= 13$, and $q=2$}
	\begin{center}
		%\large
		\begin{tabular}{|c|c|c|c|} 
			\hline
			Patient type $Ti$	& Stage 1 time,  $\lambda_{i}$ & Stage 2 time, $\mu_{i}$ & Demand ratio, $r_i$ \\
			\hline
			Type $T1$ & 10 & 0  &  4 \\
			Type $T2$ & 15 & 0  &  5 \\
			Type $T3$ & 20 & 25 & 1\\
			Type $T4$ & 15 & 35  & 3\\
			\hline
		\end{tabular}
	\end{center}
	% \vspace{0.01in}
	\label{table-Pr2}
\end{table}

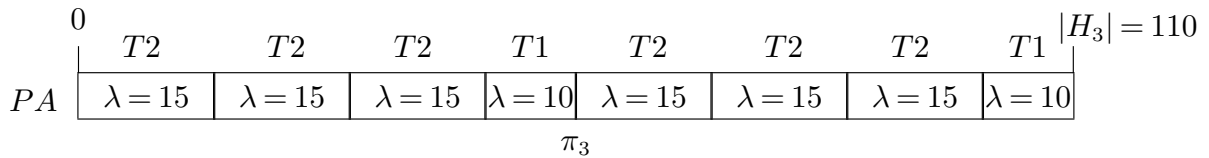
\begin{figure}[htbp]
\resizebox{\textwidth}{!}{
	{\normalsize
		\setlength{\unitlength}{1.11mm}
		\begin{picture}(135,13)(55,11) 
		%\thicklines
		
		\put(72,23){\makebox(0,0){$T2$}}
		\put(88,23){\makebox(0,0){$T2$}}
		\put(102,23){\makebox(0,0){$T2$}}
		\put(115,23){\makebox(0,0){$T1$}}
		
		\put(128,23){\makebox(0,0){$T2$}}
		\put(143,23){\makebox(0,0){$T2$}}
		\put(157,23){\makebox(0,0){$T2$}}
		\put(170,23){\makebox(0,0){$T1$}}

		\multiput(65,20)(27,0){1}{\line(0,1){3}}
		\put(65,26){\makebox(0,0){0}}
		
		\multiput(175,20)(27,0){1}{\line(0,1){3}}
		\put(181,25){\makebox(0,0){$|H_3|=110$}}

		%======Head1 ====	

		\put(65,15){\framebox(15,5){$\lambda=15$}}
		\put(80,15){\framebox(15,5){$\lambda=15$}}
		\put(95,15){\framebox(15,5){$\lambda=15$}}
		\put(110,15){\framebox(10,5){$\lambda=10$}}
	
		\put(120,15){\framebox(15,5){$\lambda=15$}}
		\put(135,15){\framebox(15,5){$\lambda=15$}}
		\put(150,15){\framebox(15,5){$\lambda=15$}}
	
		\put(165,15){\framebox(10,5){$\lambda=10$}}
		\put(60,17){\makebox(0,0){$PA$}}	
			\put(120,12){\makebox(0,0){$\pi_3$}}
		
		\end{picture}}}	
	%\vspace{-0.2in}	
	\caption{The Last Block Schedule $\pi_3$ for Example 2}
	\label{figx:block3}	
\end{figure}

Figure~\ref{figx:block4} shows the no-wait schedule $\Gamma= \{\pi_1\pi_2\pi_3\}$ of three blocks, where the no-wait constraint is maintained between the blocks in the scheduling terminology. Let the regular time available for a day be $R=300$ minutes. Thus, $PA$ experiences overtime of 65 minutes because their finish time of $\sigma$ is 365 minutes. The idle time for $PA$ is equal to 5 minutes, while there is no idle time or overtime for $P$. Patient waiting time at stage~2 is equal to 180 minutes and there is zero waiting time at stage~1.  

% Thus, $D_p(\sigma) = 0$ and $D_a(\sigma)=5$, $W_p(\sigma) = 180$, $W_a(\sigma)=0$, and the overtimes of $P$ and $PA$, $b_p(\sigma) = 0$, $b_a(\sigma)=65$.

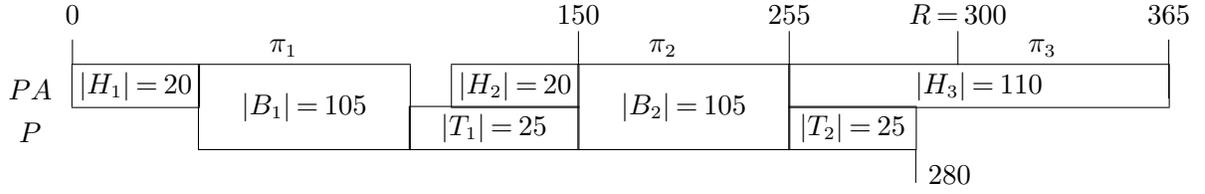
\begin{figure}[htbp]
\resizebox{\textwidth}{!}{
	{\normalsize
		\setlength{\unitlength}{1.14mm}
		\begin{picture}(145,15)(10,15) 
		%\thicklines

		\put(15,17){\makebox(0,0){$PA$}}
		
		\put(15,12){\makebox(0,0){$P$}}
		
			\multiput(20,20)(27,0){1}{\line(0,1){3}}
		\put(20,26){\makebox(0,0){0}}

		%======Head1 ====
		
		\put(20,15){\framebox(15,5){$|H_1|=20$}}

		%===========Body1=============
		
		\put(35,10){\framebox(25,10){$|B_1|=105$}}
		
		%===========Tail1=============
		
		\put(60,10){\framebox(20,5){$|T_1|=25$}}

		%======Head2 ====
		
		\put(65,15){\framebox(15,5){$|H_2|=20$}}

		%===========Body2=============
		
		\put(80,10){\framebox(25,10){$|B_2|=105$}}
		
		%===========Tail2=============
		
		\put(105,10){\framebox(15,5){$|T_2|=25$}}

		%======Head3 ====
		
		\put(105,15){\framebox(45,5){$|H_3|=110$}}

		\multiput(150,16)(27,0){1}{\line(0,1){7}}
		\put(150,26){\makebox(0,0){$365$}}
		
		\multiput(105,20)(27,5){1}{\line(0,1){4}}
		\put(105,26){\makebox(0,0){$255$}}
		
			\multiput(80,20)(27,5){1}{\line(0,1){4}}
		\put(80,26){\makebox(0,0){$150$}}
		
			\multiput(125,20)(27,5){1}{\line(0,1){4}}
		\put(125,26){\makebox(0,0){$R=300$}}
		
			\multiput(120,6)(27,5){1}{\line(0,1){4}}
		\put(124,7){\makebox(0,0){$280$}}

		\put(45,22){\makebox(0,0){$\pi_1$}}
		\put(90,22){\makebox(0,0){$\pi_2$}}
		\put(135,22){\makebox(0,0){$\pi_3$}}
		
		\end{picture}}}	
	\vspace{0.3in}	
	\caption{Example 2: The Schedule $\Gamma= \{\pi_1\pi_2\pi_3\}$ is no-wait schedule which is concatenation of three blocks.}
	\label{figx:block4}	
\end{figure}

% \section{Model for Two-stage Scheduling of a Block}\label{app1.0}

\section{Section~\ref{algsintro} Proofs}\label{app1}

\noindent\textbf{Proof of Theorem \ref{thm-block1}:} 

The proof is a reduction from 3-Partition, a known $\mathcal{NP}$-complete problem. First, the patient waiting time minimization problem with no-idle time is reformulated as a decision problem.

\noindent{\bf Decision problem:} Given $m$ types of patients with known demand ratios $r_i=1$ for all $i=1, \cdots, m$, service times  $\lambda_{i}$, $\mu_i$ of patient's type $i$ at stages 1 and 2, respectively, does there exist a no-idle time block schedule $\pi$, such that total waiting time of the patients are smaller than a given threshold $D=0$, $W_a(\pi)+W_p(\pi) \leq 0$?

\noindent{\bf 3-Partition (Garey and Johnson, 1979):} Given positive integers, $s$, $B$ and a set of integers $A=\{a_1,a_2,...,a_{3s}\}$ with $\sum_{i=1}^{3s}a_i=sB$ and $B/4 < a_i < B/2$ for $1 \leq i \leq 3s$, 3-Partition problem checks the existence of a partition of $A$ into three element sets $\{A_1,A_2,...,A_s\}$ such that $\sum_{a_j \in A_i} a_j = B, 1 \leq i \leq s$.

Given an instance of 3-Partition, we construct a specific instance of the decision problem. Consider the following problem setting:  $r_i=1$ for all $i=1, \cdots, m$. Patients are arranged in such a way that $\lambda_{q+1} \geq \lambda_{q+2} \geq \cdots \geq \lambda_m$ where $0<\lambda_i\leq \mu_i$ for $i=q+1,q+2,\cdots, m$, all other types have $\lambda_i>0$ and $\mu_i=0$ for $i=1,2,\cdots,q$. We define the service times of patients in accordance with the problem setting as follows:
   \[ \lambda_{q+1} = \lambda_{q+2} = \cdots = \lambda_m = B,\]
    \[ \mu_{q+1} = \mu_{q+2} = \cdots = \mu_m = 2B,\]
    \[ \lambda_1 =a_1, \lambda_2 =a_2, \cdots , \lambda_q = a_q,\]
    \[ \mu_1=\mu_2= \cdots = \mu_q =0.\]
 Here, it is assumed that $q=3s$ hence, there are $3s$ patients that have zero service times at stage 2. There are a total of $m-q$ patients that have nonzero service times at stage 2 ($0<\lambda \leq \mu$). We set $m=4s+1$ and $D=0$. 
 
The decision problem is in class $\mathcal{NP}$ for a given block schedule $\pi$. We show that the 3-Partition problem is reducible to the decision problem. That is, we show that there exists a solution (schedule) to the 3-Partition if and only if there exists a solution to the no-idle block scheduling decision problem with the total patient waiting time equal to 0. 

\begin{enumerate}
    \item If there exists a solution (schedule) to the 3-Partition, we will show that there exists a solution to the no-idle block schedule with the total patient waiting time equal to 0.
    
    Patients of type $i$ where $i=q+1, q+2, \cdots , m$, can be scheduled in a no-wait block schedule as in  Figure~\ref{fig:proof1}. By doing so, the partial schedule would have $m-q=s+1$ patients in a block and $s$ idle time slots of length $B$ each, for the $PA$. Since there exists a 3-Partition to the problem, it is possible to schedule patients of type $i$ where $i=1,2,\cdots,3s$, in these idle time-slots each slot fitting three patients without any idle at $PA$. Consequently, this would yield the no-idle block schedule as shown in Figure~\ref{fig:proof2}. Thus, the if part of the statement holds.
    
    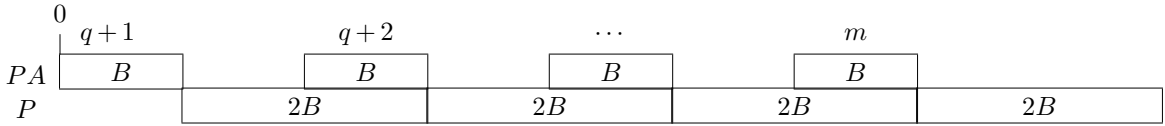
\begin{figure}[htbp]
    \resizebox{\textwidth}{!}{
		{\normalsize
			\setlength{\unitlength}{1mm}
			\begin{picture}(180,24)(12,4) 
				%\thicklines
				\put(27,23){\makebox(0,0){$q+1$}}
				
				%\put(43,23){\makebox(0,0){}}
				
				\put(65,23){\makebox(0,0){$q+2$}}
				
				\put(71,23){\makebox(0,0){}}
				
				\put(101,23){\makebox(0,0){$\cdots$}}
				
				\put(102,23){\makebox(0,0){}}
				
				\put(137,23){\makebox(0,0){$m$}}
				\multiput(20,20)(27,0){1}{\line(0,1){3}}
				\put(20,26){\makebox(0,0){0}}

				%\multiput(130,5)(27,0){1}{\line(0,1){3}}
				
				%======PA ====
				
				\put(20,15){\framebox(18,5){$B$}}
				
				\put(56,15){\framebox(18,5){$B$}}
				
				\put(92,15){\framebox(18,5){$B$}}
				
				\put(128,15){\framebox(18,5){$B$}}
				\put(15,17){\makebox(0,0){$PA$}}
				
				%======P==
				\put(38,10){\framebox(36,5){$2B$}}
				\put(74,10){\framebox(36,5){$2B$}}
				\put(110,10){\framebox(36,5){$2B$}}
				\put(146,10){\framebox(36,5){$2B$}}
				\put(15,12){\makebox(0,0){$P$}}
				
		\end{picture}}}	
			\vspace{-3 mm}	
		\caption{Block schedule construction}
		\label{fig:proof1}	
	\end{figure}

    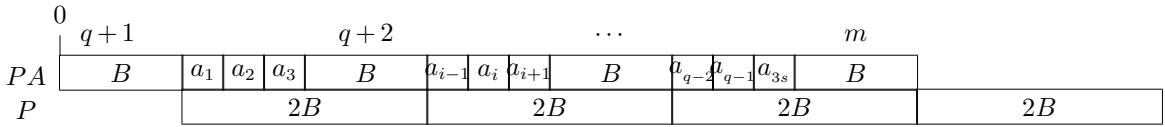
\begin{figure}[htbp]
    \resizebox{\textwidth}{!}{
		{\normalsize
			\setlength{\unitlength}{1mm}
			\begin{picture}(180,24)(12,4) 
				%\thicklines
				
				\put(27,23){\makebox(0,0){$q+1$}}
				
				%\put(43,23){\makebox(0,0){}}
				
				\put(65,23){\makebox(0,0){$q+2$}}
				
				\put(71,23){\makebox(0,0){}}
				
				\put(101,23){\makebox(0,0){$\cdots$}}
				
				\put(102,23){\makebox(0,0){}}
				
				\put(137,23){\makebox(0,0){$m$}}
				\multiput(20,20)(27,0){1}{\line(0,1){3}}
				\put(20,26){\makebox(0,0){0}}

				%\multiput(130,5)(27,0){1}{\line(0,1){3}}
				
				%======PA ====
				
				\put(20,15){\framebox(18,5){$B$}}
				\put(38,15){\framebox(6,5){$a_1$}}
				\put(44,15){\framebox(6,5){$a_2$}}
				\put(50,15){\framebox(6,5){$a_3$}}
				\put(56,15){\framebox(18,5){$B$}}
				\put(74,15){\framebox(6,5){$a_{i-1}$}}
				\put(80,15){\framebox(6,5){$a_{i}$}}
				\put(86,15){\framebox(6,5){$a_{i+1}$}}
				\put(92,15){\framebox(18,5){$B$}}
				\put(110,15){\framebox(6,5){$a_{_{q-2}}$}}
				\put(116,15){\framebox(6,5){$a_{_{q-1}}$}}
				\put(122,15){\framebox(6,5){$a_{_{3s}}$}}
				\put(128,15){\framebox(18,5){$B$}}
				\put(15,17){\makebox(0,0){$PA$}}
				
				%======P==
				\put(38,10){\framebox(36,5){$2B$}}
				\put(74,10){\framebox(36,5){$2B$}}
				\put(110,10){\framebox(36,5){$2B$}}
				\put(146,10){\framebox(36,5){$2B$}}
				\put(15,12){\makebox(0,0){$P$}}
				
		\end{picture}}}	
			\vspace{-3 mm}	
		\caption{No-idle time block schedule}
		\label{fig:proof2}	
	\end{figure}
    
    \item If there exists a solution to the no-idle block schedule with the total patient waiting time equal to 0, we will show that there exists a solution (schedule) to the 3-Partition. 
    
    Since the total patient waiting time is equal to 0, patient of type $i=q+1,q+2,\cdots,m$ must be scheduled as no-wait schedule as shown in Figure~\ref{fig:proof1}. This partial schedule creates $s$ idle time slots of length $B$ each, for the $PA$.
    
    Since there exists a no-idle block schedule and $\sum_{i=1}^{3s} a_i  = sB$, all remaining patient of type $i=1, 2,\cdots, q$, are scheduled in $s$ idle time-slots without any idle time. Let denote $A_l$ be the subset of patient type $i=1,2,\cdots,3s$, patients scheduled in time-slot $l$, where $l=1,2, \ldots, s$. 
	 Note that $\sum_{a_i \in A_l} a_i  = B$, $l=1,2, \ldots, s$ as there exists a no-idle block schedule.  We prove the following claim.
	 \noindent\textbf{Claim:} The size of each set $A_l$ is equal to 3, where $l=1,2,\cdots,s$.\\
	\noindent\textbf{Proof:} Note that from 3-Partition problem, we have $B/4<a_i<B/2$. Suppose the size of each set $A_l$ is equal to 4, then $B<\sum_{a_i \in A_l} a_i<2B$. As $\sum_{a_i \in A_l} a_i$ exceeds $B$, it violates the assumption of no-idle block schedule. Thus, the size of each set $A_l$ must be smaller than 4. Let's assume that the size of each set $A_l$ is equal to 2, then $B/2<\sum_{a_i \in A_l} a_i<B$. As $\sum_{a_i \in A_l} a_i$ is strictly smaller than $B$, it creates idle time in the schedule. Hence the size of each set $A_l$ must be greater than 2. The above facts imply that the size of each set $A_l$ must be equal to 3. Thus, the claim holds. $\blacksquare$
	
  As a consequence of the above claim, the no-idle block schedule forms the schedule structure similar to that shown in Figure~\ref{fig:proof2} which provides a partition of the patients of type $i=1,2,\cdots,3s$, into three element sets $\{A_1,A_2,\cdots,A_s\}$ such that $\sum_{a_i \in A_l} a_i  = B$, $l=1,2,\cdots,s$. Therefore, there exists a solution to the 3-Partition problem and the only if part of the statement holds. $\blacksquare$
    
\end{enumerate}
\medskip

For simplicity in exposition, we assume that $r_i=1$ for all $i=1,\cdots,m$. The results in Lemmas below are valid for any general $r_i \geq 1$ for all $i=1,\cdots,m$. 

\vspace{5mm}
\noindent\textbf{Proof of Lemma~\ref{L-block1}:} 
	
Let Algorithm~\ref{alg-MA} provide a sequence $\sigma$ of patients given an appointment at the clinic in a block, where $\sigma=\{ \sigma_{1}, \sigma_{2}, \ldots, \sigma_{m-q}, \sigma_{m-q+1}, \ldots, \sigma_{m}\}$ and $\sigma_{j}$ denote the patient type sequenced in $j^{th}$ position in $\sigma$. Consider two partial sequences, $\sigma^f$ and $\sigma^b$ of $\sigma$, where $\sigma= \sigma^f \cup \sigma^b$:\\ (i) Front sequence, $\sigma^f =\{\sigma_{1}, \sigma_{2}, \ldots, \sigma_{m-q}\}= \{q+1, q+2, \ldots, m\}$, \\(ii) Back sequence, $\sigma^b=\{\sigma_{m-q+1}, \sigma_{m-q+2}, \ldots, \sigma_{m}\}= \{1, 2, \ldots, q\}$.

We first show that the front sequence  $\sigma^f=\{\sigma_{1}, \sigma_{2}, \ldots, \sigma_{m-q}\}= \{q+1, q+2, \ldots, m\}$ forms no-idle schedule. Note that the patients in $\sigma^f$ see both $PA$ and $P$ and $0<\lambda_{\sigma_{j}}\leq\mu_{\sigma_{j}}$, for $j=1,2, \ldots, m-q$.

% In $\sigma^f$, the idle time of a physician and a physician's assistant for the $j^{th}$ patient can be expressed as follows respectively,

In schedule formed by $\sigma^f$, the idle time of a $P$ and a $PA$ for the $j^{th}$ patient can be expressed as follows respectively,
    \[d_a^j=e_a^j-f_a^{j-1},\]
    \[d_p^j=e_p^j-f_p^{j-1}.\]
All patients in $\sigma^f$ see the $PA$ first, thus there is no idle time for the $PA$, $e_a^j=f_a^{j-1}$,  $d_a^j=e_a^j-f_a^{j-1}=0$ for all $j=1,2, \ldots, m-q$.
The outpatient clinic setting described in Section~\ref{pd} assumes that there exists one $P$ and one $PA$ present in the system. Therefore, the start time of $j^{th}$ patient at stage 2 in $\sigma^f$ can be expressed as 
    $e_p^j=\max\{f_a^j,f_p^{j-1}\}$. 
The $P$ can see $j^{th}$ patient, either immediately after he/she leaves stage 1 or immediately after  $(j-1)^{th}$ patient leaves stage 2, whichever is the greatest. Thus, the $P$ idle time becomes,
    $d_p^j=\max\{f_a^j,f_p^{j-1}\}-f_p^{j-1}$.
We now show the following claim:

\noindent\textbf{Claim:} In  schedule formed by $\sigma^f$, $e_p^j=f_p^{j-1}$ holds $j=1,2, \ldots, m-q$.

\noindent\textbf{Proof:} The problem setting assumes that the $PA$ begins seeing the first patient at time 0 hence, the finish time of $j^{\text{th}}$ patient at stage 1 can be recursively calculated as:
    \[f_a^j=\sum_{l=1}^j \lambda_l.\]
Similarly, the $P$ begins seeing the first patient immediately after the patient leaves stage 1. Thus, the finish time of $j^{\text{th}}$ patient at stage 2 can be defined as follows:
    \[f_p^j \geq f_a^1 + \sum_{l=1}^j \mu_l.\]
This general expression becomes an equality when there is no idle time in the schedule. Without loss of generality, $f_p^{j-1} \geq f_a^1 + \sum_{l=1}^{j-1} \mu_l$ follows from the above inequality when $j>1$. 

By construction $0 < \lambda_{\sigma_j} \leq \mu_{\sigma_j}$, for $j=1,2, \ldots, m-q$, it follows that
% \begin{align*}
%     &\sum_{l=1}^j \mu_l \geq \sum_{l=1}^j \lambda_l\\
%     &\sum_{l=1}^{j-1} \mu_l \geq \sum_{l=1}^{j-1} \lambda_l\\
%     &\lambda_1 + \sum_{l=1}^{j-1} \mu_l \geq \lambda_j + \sum_{l=1}^{j-1} \lambda_l\\
% \end{align*}
    \[\sum_{l=1}^j \mu_l \geq \sum_{l=1}^j \lambda_l, \]

    \[\sum_{l=1}^{j-1} \mu_l \geq \sum_{l=1}^{j-1} \lambda_l,\]

    \[\lambda_1 + \sum_{l=1}^{j-1} \mu_l \geq \lambda_j + \sum_{l=1}^{j-1} \lambda_l.  \]

Since Algorithm~\ref{alg-MA} assumes $\lambda_{1} \geq \lambda_{2} \geq \ldots \geq \lambda_{m-q}$ hence, $\lambda_1 \geq \lambda_j$. Therefore, the above equation holds and can be rewritten as  $f_a^1 + \sum_{l=1}^{j-1} \mu_l  \geq \sum_{l=1}^j \lambda_l$. Thus, $f_p^{j-1} \geq f_a^j$ and $e_p^j=\max\{f_a^j,f_p^{j-1}\}=f_p^{j-1}$, the claim holds. $\blacksquare$

Following the claim, $d_p^j=f_p^{j-1}-f_p^{j-1}=0$ holds for $j=1,2, \ldots, m-q$. Thus, Algorithm~\ref{alg-MA} provides a no-idle schedule formed by $\sigma^f$. 

We now show that the block formed by $\sigma= \sigma^f \cup \sigma^b$ is a no-idle schedule. Note that the algorithm schedules patients in $\sigma^b=\{\sigma_{m-q+1}, \sigma_{m-q+2}, \ldots, \sigma_{m}\}= \{1, 2, \ldots, q\}$ at the end of schedule  formed by $\sigma^f$ as early as possible. Since $\mu_{1}= \mu_{2}=\ldots =\mu_{q}=0$ and there is no idle time at $PA$. Thus, the schedule formed by $\sigma$ is a no-idle schedule. $\blacksquare$ %\qedsymbol{$\blacksquare$}

\medskip

\noindent\textbf{Proof of Lemma \ref{l-block2}:} 

Without loss of generality, we assume that $r_i=1$ for all $i=1,\cdots,m$. Later, the results in this Lemma is generalized for $r_i \geq 1$ for all $i=1,\cdots, m$.
Recall the notation of front sequence, $\sigma^f =\{\sigma_{1}, \sigma_{2}, \ldots, \sigma_{m-q}\}= \{q+1, q+2, \ldots, m\}$, and back sequence, $\sigma^b=\{\sigma_{m-q+1}, \sigma_{m-q+2}, \ldots, \sigma_{m}\}= \{1, 2, \ldots, q\}$. We derive the expression from Lemma~\ref{l-block2} separately for the partial schedules formed by these sequences. 
In the block schedule $\pi$ given by Algorithm~\ref{alg-MA}, $W_a(\pi)=0$ as $\pi$ is formed by both $\sigma^f$ and $\sigma^b$.
The patient waiting time in $\pi$ formed by $\sigma^b$ is equal to zero at the second stage, as these patients only see $PA$ and do not see $P$. Following from Lemma~\ref{L-block1}, $PA$ begins seeing the first patient at time 0 thus, for any arbitrary no idle block schedule $\pi$, the waiting times of the patients are recursively calculated for $\pi$ formed by $\sigma^f$ as in Table~\ref{recu}.
 $$W_p(\pi)=(m-q-1)\mu_{q+1}+(m-q-2)\mu_{q+2}+ \cdots +\mu_{m-1} - (m-q-1)\lambda_{q+2}-(m-q-2)\lambda_{q+3}- \cdots -\lambda_m.$$

\begin{table}[]
\caption{Patient wait times in the No-idle schedule $\pi$, where $r_i=1$ for all $i=1, \cdots, m$}
    \centering
    \resizebox{\textwidth}{!}{\begin{tabular}{|c|c|c|c|c|c|c|}
    \hline
         $j$ & $e_a^j$ & $f_a^j$ & $e_p^j$ & $f_p^j$ & $w_a^j$ & $w_p^j$ \\
         \hline
         $q+1$ & 0 & $\lambda_{q+1}$ & $\lambda_{q+1}$ & $\lambda_{q+1}$ + $\mu_{q+1}$ & 0 & 0\\
         $q+2$ & $\lambda_{q+1}$ & $\lambda_{q+1}$ + $\lambda_{q+2}$ & $\lambda_{q+1}$ + $\mu_{q+1}$ & $\lambda_{q+1}$ + $\mu_{q+1}$ + $\mu_{q+2}$ & 0 & $\mu_{q+1}- \lambda_{q+2}$  \\
         $q+3$ & $\lambda_{q+1} + \lambda_{q+2}$ & $\lambda_{q+1} + \lambda_{q+2} + \lambda_{q+3}$ & $\lambda_{q+1} + \mu_{q+1} + \mu_{q+2}$ & $\lambda_{q+1} + \mu_{q+1} + \mu_{q+2} + \mu_{q+3}$ & 0 & $\mu_{q+1} + \mu_{q+2} - \lambda_{q+2}- \lambda_{q+3}$ \\
         \vdots & \vdots & \vdots & \vdots & \vdots & \vdots & \vdots \\
         $i$ & $\sum_{l=q+1}^{i-1} \lambda_l$ & $\sum_{l=q+1}^i \lambda_l$ & $\lambda_{q+1} + \sum_{l=q+1}^{i-1} \mu_l$ & $\lambda_{q+1} + \sum_{l=q+1}^i \mu_l$ & 0 & $\sum_{l=q+1}^{i-1} \mu_l - \sum_{l=q+2}^i \lambda_l$\\
         \vdots & \vdots & \vdots & \vdots & \vdots & \vdots & \vdots \\
         $m$ & $\sum_{l=q+1}^{m-1} \lambda_l$ & $\sum_{l=q+1}^{m} \lambda_l$ & $\lambda_{q+1} + \sum_{l=q+1}^{m-1} \mu_l$ & $\lambda_{q+1} + \sum_{l=q+1}^{m} \mu_l$ & 0 & $\sum_{l=q+1}^{m-1} \mu_l - \sum_{l=q+2}^{m} \lambda_l$\\
         \hline
    \end{tabular}}
    % \vspace{0.01in}
    
    \label{recu}
\end{table}
 
Let $\gamma= m-q$ and rearranging, we have
$$W_p(\pi)=(\gamma-1)\mu_{q+1}+(\gamma-2)\mu_{q+2}+ \cdots +\mu_{m-1} - (\gamma-1)\lambda_{q+2}-(\gamma-2)\lambda_{q+3}-\cdots-\lambda_m.$$ 

 By adding, we have $$W_a(\pi)+W_p(\pi) = \sum_{j=1}^{\gamma-1} (\gamma-j)(\mu_{q+j} - \lambda_{q+j+1}), \,\, \hbox{where} \,\,\gamma= m-q.$$

 We now generalize the above formula for $r_i \geq 1$ for all $i=1,\cdots, m$.
Note that Algorithm~\ref{alg-MA} provides a sequence, $\sigma$ for block $\pi$, where $\sigma=\{\sigma_{1}, \sigma_{2}, \ldots, \sigma_{r-v}, \sigma_{r-v+1}, \ldots, \sigma_{r}\}$. Recall the two partial sequences for a block, $\sigma^f$ and $\sigma^b$ of $\sigma$, where $\sigma= \sigma^f \cup \sigma^b$, $\sigma^f =\{\sigma_{1}, \sigma_{2}, \ldots, \sigma_{r-v}\}$, and $\sigma^b=\{\sigma_{r-v+1}, \sigma_{r-v+2}, \ldots, \sigma_{r}\}$.

By setting $q=v$ and $m=r$,  we have $W_a(\pi)+W_p(\pi) = \sum_{j=1}^{\gamma-1} (\gamma-j)(\mu_{v+j} - \lambda_{v+j+1}), \,\, \hbox{where} \,\,\gamma= r-v$ for the general case. $\blacksquare$	
	
\vspace{5mm}
\noindent \textbf{Proof of Lemma \ref{modif}:} 

We prove Lemma~\ref{modif} in three parts, first, we show that schedule $\hat{s} \in S$ yields a no-idle front block schedule, then we show that Algorithm~\ref{alg-MA} yields the smallest patient waiting time among all possible schedules obtained (in set $S$) for the front schedule with $Q^+$ Group Patients, where patients are ordered: $\lambda_{q+1} \geq \lambda_{q+2} \geq \cdots \geq \lambda_m$.

\noindent\textbf{Part 1:} Algorithm~\ref{alg-MA} provides a schedule in $S_1$ which is subset of $S$. In $S_1$, the ties of $PA$ service times ($\lambda_{i}$) are broken in favor of the patient with the smallest $(\mu_i-\lambda_i)$ in the front schedule. Note that $S$ is the set of all possible schedules for the front schedule with $Q^+$ Group Patients are ordered: $\lambda_{q+1} \geq \lambda_{q+2} \geq \cdots \geq \lambda_m$.  
We note that any schedule $\hat{s} \in S$ is no idle schedule. The proof is similar to that of Lemma~\ref{L-block1}. 

\noindent \textbf{Part 2:} The front schedule provided by Algorithm~\ref{alg-MA} arranges the patients such that $\lambda_{q+1} \geq \lambda_{q+2} \geq \ldots \geq \lambda_m$, for $i=q+1,q+2, \ldots, m$. Patients that have equal service times at stage 1 ($\lambda_i$) can be arranged in any order in the front schedule which could yield multiple solutions. However, in such a case Algorithm~\ref{alg-MA} breaks the ties in favor of the patient with the smallest service time in stage 2, which ensures obtaining a unique front schedule $s' \in S_1$.

\noindent\textbf{Part 3:} We use interchange arrangement to prove that schedule $s' \in S_1$ has the minimum total patient waiting time among the schedules in $S$. Let $\bar{s} \in \{S-S_1\}$ be a no-idle block schedule that is not given by Algorithm~\ref{alg-MA}. Suppose, we assume that $\bar{s}$ has the minimum total patient waiting time among all no-idle block schedules. In the schedule $\bar{s} \in \{S-S_1\}$, patients are scheduled in the order of $\lambda_{q+1} \geq \lambda_{q+2} \geq \cdots \geq \lambda_{q+j} \geq \lambda_{q+j+1} \geq \cdots \geq \lambda_m$ for the front schedule and $\mu_i \geq \lambda_i$ holds for every patient $i$. The schedule $\bar{s} \in \{S-S_1\}$ is illustrated as in Figure~\ref{fig:mod_alg1}.

\begin{figure}[htbp]
\resizebox{\textwidth}{!}{
	{\normalsize
		\setlength{\unitlength}{1mm}
		\begin{picture}(175,24)(-28,4) 
			%\thicklines
			
			\put(-20,20){\line(0,1){5}}
			\put(-20,27){\makebox(0,0){0}}
			
			\put(-25,17){\makebox(0,0){$PA$}}
	      	\put(-25,12){\makebox(0,0){$P$}}
	      	 	
	      	\put(20,22){\vector(1,0){18}}
	      	\put(3,23){$L_1$}
	      	\put(20,22){\vector(-1,0){40}}
	      	\put(38,20){\line(0,1){5}}

			\put(20,4){\vector(1,0){19}}
			\put(13,5.5){$L_2$}
			\put(20,4){\vector(-1,0){28}}
			\put(40,10){\line(0,-1){8}}
			\put(-8,10){\line(0,-1){8}}

			%\multiput(130,5)(27,0){1}{\line(0,1){3}}

			%======PA & P Fornt====
			
			\put(-20,15){\framebox(12,5){$\lambda_{q+1}$}}
			\put(-8,15){\framebox(9,5){$\lambda_{q+2}$}}
	
			\put(-15,15){\makebox(50,5){$\cdots \cdots$}}
		
			\put(-8,10){\framebox(13,5){$\mu_{q+1}$}}
			\put(5,10){\framebox(13,5){$\mu_{q+2}$}}
			
				\put(3,9){\makebox(50,5){$\cdots \cdots$}}

			%======PA ====
			
			\put(20,15){\framebox(18,5){$\lambda_{q+j-1}$}}
			\put(38,15){\framebox(12,5){$\lambda_{q+j}$}}
			\put(50,15){\framebox(12,5){\tiny$\lambda_{q+j+1}$}}
			
			\put(62,15){\framebox(10,5){\tiny$\lambda_{q+j+2}$}}
			\put(68,15){\makebox(50,5){$\cdots  \cdots \cdots  \cdots$}}
		
			\put(112,15){\framebox(7,5){$\lambda_m$}}
	
			\put(119,15){\framebox(6,5){$\lambda_1$}}
			\put(125,15){\framebox(6,5){$\lambda_2$}}
			\put(132,15){\makebox(7,5){$\cdots  \cdots$}}
			\put(140,15){\framebox(5,5){$\lambda_q$}}

			%======P==
			\put(40,10){\framebox(20,5){$\mu_{q+j-1}$}}
			%\put(38,10){\framebox(22,5){$\mu_{q+j-1}$}}
			\put(60,10){\framebox(20,5){$\mu_{q+j}$}}
			\put(80,10){\framebox(15,5){$\mu_{q+j+1}$}}
			\put(95,10){\framebox(13,5){\small$\mu_{q+j+2}$}}
			\put(116,9){\makebox(6,5){$\cdots \cdots  \cdots $}}
			\put(128,9.5){\framebox(8,5){$\mu_m$}}

	\end{picture}}}	
	%	\vspace{0.1in}	
	\caption{No-idle block schedule $\bar{s} \in S-S_1$}
	\label{fig:mod_alg1}	
\end{figure}
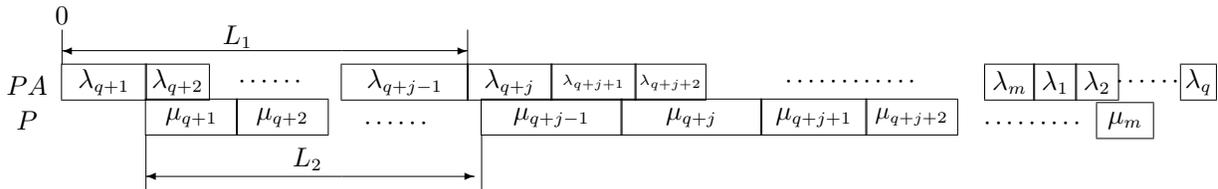

 %    \begin{figure}[htbp]
	% 	{\normalsize
	% 		\setlength{\unitlength}{1mm}
	% 		\begin{picture}(180,24)(12,4) 
	% 			\thicklines
				
	% 			\multiput(20,20)(27,0){1}{\line(0,1){3}}
	% 			\put(20,26){\makebox(0,0){0}}

	% 			%\multiput(130,5)(27,0){1}{\line(0,1){3}}
				
	% 			%======PA ====
				
	% 			\put(20,15){\framebox(18,5){$\lambda_{q+1}$}}
	% 			\put(38,15){\framebox(12,5){$\lambda_{q+2}$}}
	% 			\put(50,15){\framebox(12,5){$\lambda_{q+3}$}}
				
	% 			\put(62,15){\framebox(10,5){$\lambda_{q+4}$}}
	% 			\put(72,15){\makebox(50,5){$\cdots$}}
	% 			\put(122,15){\framebox(7,5){$\lambda_m$}}
	% 			\put(129,15){\framebox(6,5){$\lambda_1$}}
	% 			\put(135,15){\framebox(6,5){$\lambda_2$}}
	% 			\put(141,15){\makebox(7,5){$\cdots$}}
	% 			\put(152,15){\framebox(5,5){$\lambda_q$}}
	% 			\put(15,17){\makebox(0,0){$PA$}}
				
	% 			%======P==
	% 			\put(38,10){\framebox(22,5){$\mu_{q+1}$}}
	% 			\put(60,10){\framebox(20,5){$\mu_{q+2}$}}
	% 			\put(80,10){\framebox(15,5){$\mu_{q+3}$}}
	% 			\put(95,10){\framebox(13,5){$\mu_{q+4}$}}
	% 			\put(121,10){\makebox(6,5){$\cdots$}}
	% 			\put(142,10){\framebox(8,5){$\mu_m$}}
	% 			\put(15,12){\makebox(0,0){$P$}}
				
	% 	\end{picture}}	
	% 	%	\vspace{0.1in}	
	% 	\caption{No-idle block schedule $S$}
	% 	\label{fig:mod_alg1}	
	% \end{figure}
There exist two adjacent patients in $\bar{s} \in \{S-S_1\}$ such that $\lambda_{q+j}=\lambda_{q+j+1}$ and $\mu_{q+j} > \mu_{q+j+1}$. Switching the order of these patients in the schedule would yield another no-idle block schedule $s'$, which is illustrated in Figure~\ref{fig:mod_alg2}. This change only affects the waiting times of the patients $q+j$ and $q+j+1$, as the appointment start and finish times of the patients $q+j+2,\cdots,m$ remain the same. The waiting times of patients $q+j$ and $q+j+1$ in the schedule $\bar{s}$ can be calculated as below. Let $L_1 =\sum_{u=q+2}^{q+j-1} \lambda_{u}$ and $L_2 =\sum_{u=q+1}^{q+j-2} \mu_{u}$. Note that $L_2 -L_1 \geq 0$.
\[w_p^{q+j}=L_2+\mu_{q+j-1}-L_1-\lambda_{q+j},\]
\[w_p^{q+j+1}=L_2+\mu_{q+j-1}+\mu_{q+j}-L_1-\lambda_{q+j}-\lambda_{q+j+1},\]
\[w_p^{q+j}+w_p^{q+j+1}=2L_2+2\mu_{q+j-1}+\mu_{q+j}-2L_1 -2\lambda_{q+j}-\lambda_{q+j+1},\]
\[w_p^{q+j}+w_p^{q+j+1}=2(L_2-L_1)+2\mu_{q+j-1}+\mu_{q+j}-3\lambda_{q+j}.\]
By switching the order of the patients $q+j$ and $q+j+1$, we obtain schedule $s'$ where the waiting times of these patients are updated as follows:
\[w_{p'}^{q+j}=L_2+\mu_{q+j-1}+\mu_{q+j+1}-L_1-\lambda_{q+j}-\lambda_{q+j+1},\]
\[w_{p'}^{q+j+1}=L_2+\mu_{q+j-1}-L_1-\lambda_{q+j+1},\]
\[w_{p'}^{q+j}+w_{p'}^{q+j+1}=2L_2+2\mu_{q+j-1}+\mu_{q+j+1}-2L_1-
\lambda_{q+j}-2\lambda_{q+j+1},\]
\[w_{p'}^{q+j}+w_{p'}^{q+j+1}=2(L_2-L_1)+2\mu_{q+j-1}+\mu_{q+j+1}-3\lambda_{q+j}.\]

The difference between the total waiting times of the patients in schedule $\bar{s}$ and $s'$ is equal to \[\Delta(w_p)=(w_{p'}^{q+j}+w_{p'}^{q+j+1})-(w_{p}^{q+j}+w_{p}^{q+j+1}),\]
\[\Delta(w_p)=(2\mu_{q+j-1}+\mu_{q+j+1}-3\lambda_{q+j})-(2\mu_{q+j-1}+\mu_{q+j}-3\lambda_{q+j}),\]
\[\Delta(w_p)=\mu_{q+j+1}-\mu_{q+j}.\]

Since $\mu_{j} > \mu_{j+1}$, it follows that $\Delta(w_p)<0$. Thus, the schedule $s'$ yields a smaller total patient waiting time than schedule $\bar{s}$. This contradicts the assumption that there exists a schedule $\bar{s} \in \{S-S_1\}$ that yields the minimum total patient waiting time among all no-idle block schedules in $S$. 
For any schedule $\bar{s} \in \{S-S_1\}$, a finite number of such switching the order of two adjacent patients transforms $\bar{s}$ into schedule $s' \in S_1$. Each switching improves the total patient waiting time. Thus, $s' \in S_1$ is optimal among set $S$. Therefore, we conclude that Algorithm~\ref{alg-MA} provides a schedule in $S_1$ with the smallest total patient waiting time among all possible no-idle schedules in $S$. $\blacksquare$

\begin{figure}[htbp]
\resizebox{\textwidth}{!}{
	{\normalsize
		\setlength{\unitlength}{1mm}
		\begin{picture}(175,24)(-28,4) 
			%\thicklines
			
			\put(-20,20){\line(0,1){5}}
			\put(-20,27){\makebox(0,0){0}}
			
			\put(-25,17){\makebox(0,0){$PA$}}
			\put(-25,12){\makebox(0,0){$P$}}
			
				\put(20,22){\vector(1,0){18}}
			\put(3,23){$L_1$}
			\put(20,22){\vector(-1,0){40}}
			\put(38,20){\line(0,1){5}}

			\put(20,4){\vector(1,0){19}}
			\put(13,5.5){$L_2$}
			\put(20,4){\vector(-1,0){28}}
			\put(40,10){\line(0,-1){8}}
			\put(-8,10){\line(0,-1){8}}

			%\multiput(130,5)(27,0){1}{\line(0,1){3}}

			%======PA & P Fornt====
			
			\put(-20,15){\framebox(12,5){$\lambda_{q+1}$}}
			\put(-8,15){\framebox(9,5){$\lambda_{q+2}$}}
			
			\put(-15,15){\makebox(50,5){$\cdots \cdots$}}
			
			\put(-8,10){\framebox(13,5){$\mu_{q+1}$}}
			\put(5,10){\framebox(13,5){$\mu_{q+2}$}}
			
			\put(3,9){\makebox(50,5){$\cdots \cdots$}}

			%======PA ====
			
			\put(20,15){\framebox(18,5){$\lambda_{q+j-1}$}}
		%	\put(38,15){\framebox(12,5){$\lambda_{q+j}$}}
		%	\put(50,15){\framebox(12,5){\tiny$\lambda_{q+j+1}$}}
			\put(62,15){\framebox(10,5){\tiny$\lambda_{q+j+2}$}}
			\put(68,15){\makebox(50,5){$\cdots  \cdots \cdots  \cdots$}}
			
			\put(38,15){\framebox(12,5){\tiny$\lambda_{q+j+1}$}}
			\put(50,15){\framebox(12,5){$\lambda_{q+j}$}}

			\put(112,15){\framebox(7,5){$\lambda_m$}}
			
			\put(119,15){\framebox(6,5){$\lambda_1$}}
			\put(125,15){\framebox(6,5){$\lambda_2$}}
			\put(132,15){\makebox(7,5){$\cdots  \cdots$}}
			\put(140,15){\framebox(5,5){$\lambda_q$}}

			%======P==
			\put(40,10){\framebox(20,5){$\mu_{q+j-1}$}}
		%	\put(60,10){\framebox(20,5){$\mu_{q+j}$}}
		%	\put(80,10){\framebox(15,5){$\mu_{q+j+1}$}}
			\put(95,10){\framebox(13,5){$\mu_{q+j+2}$}}
			\put(116,9){\makebox(6,5){$\cdots \cdots  \cdots $}}
			\put(128,9.5){\framebox(8,5){$\mu_m$}}
			
				\put(60,10){\framebox(15,5){$\mu_{q+j+1}$}}
				\put(75,10){\framebox(20,5){$\mu_{q+j}$}}

	\end{picture}}}	
	%	\vspace{0.1in}	
	\caption{No-idle block schedule $s' \in S-S_1$}
	\label{fig:mod_alg2}	
\end{figure}
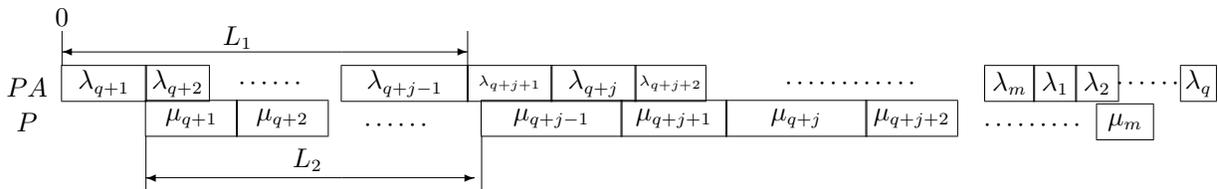

% \begin{figure}[htbp]
% 		{\normalsize
% 			\setlength{\unitlength}{1mm}
% 			\begin{picture}(180,24)(12,4) 
% 				\thicklines
				
% 				\multiput(20,20)(27,0){1}{\line(0,1){3}}
% 				\put(20,26){\makebox(0,0){0}}

% 				%\multiput(130,5)(27,0){1}{\line(0,1){3}}
				
% 				%======PA ====
				
% 				\put(20,15){\framebox(18,5){$\lambda_{q+1}$}}
% 				\put(38,15){\framebox(12,5){$\lambda_{q+3}$}}
% 				\put(50,15){\framebox(12,5){$\lambda_{q+2}$}}
				
% 				\put(62,15){\framebox(10,5){$\lambda_{q+4}$}}
% 				\put(72,15){\makebox(50,5){$\cdots$}}
% 				\put(122,15){\framebox(7,5){$\lambda_m$}}
% 				\put(129,15){\framebox(6,5){$\lambda_1$}}
% 				\put(135,15){\framebox(6,5){$\lambda_2$}}
% 				\put(141,15){\makebox(7,5){$\cdots$}}
% 				\put(152,15){\framebox(5,5){$\lambda_q$}}
% 				\put(15,17){\makebox(0,0){$PA$}}
				
% 				%======P==
% 				\put(38,10){\framebox(22,5){$\mu_{q+1}$}}
% 				\put(60,10){\framebox(15,5){$\mu_{q+3}$}}
% 				\put(75,10){\framebox(20,5){$\mu_{q+2}$}}
% 				\put(95,10){\framebox(13,5){$\mu_{q+4}$}}
% 				\put(121,10){\makebox(6,5){$\cdots$}}
% 				\put(142,10){\framebox(8,5){$\mu_m$}}
% 				\put(15,12){\makebox(0,0){$P$}}
				
% 		\end{picture}}	
% 		%	\vspace{0.1in}	
% 		\caption{No-idle block schedule $S'$}
% 		\label{fig:mod_alg2}	
% 	\end{figure}

\newpage
\noindent \textbf{Proof of Lemma~\ref{Perf-bound-block1}:} 

Let Algorithm~\ref{alg-MA} provides a sequence, $\sigma$ for block $\pi$, where $\sigma=\{\sigma_{1}, \sigma_{2}, \ldots, \sigma_{r-v}, \sigma_{r-v+1}, \ldots, \sigma_{r}\}$. Recall the two partial sequences for a block, $\sigma^f$ and $\sigma^b$ of $\sigma$, where $\sigma= \sigma^f \cup \sigma^b$, $\sigma^f =\{\sigma_{1}, \sigma_{2}, \ldots, \sigma_{r-v}\}$, and $\sigma^b=\{\sigma_{r-v+1}, \sigma_{r-v+2}, \ldots, \sigma_{r}\}$.

It is sufficient to develop waiting time expression for partial schedule formed by $\sigma^f$ as patient waiting time occurs only at $\sigma^f$. Note that $W_a(\sigma^f)=0$,  $W_p(\sigma^b)=0$ and $W_a(\sigma^b)=0$. Thus, we develop below the expression to quantify, $W_p(\sigma^f)$.

Wait time of patient schedule for $P$ at $\sigma_{1}$: $0$,

Wait time of patient schedule for $P$ at $\sigma_{2}$: $\mu_{\sigma_{1}}- \lambda_{\sigma_{2}}$. Since $\mu_{\sigma_{1}} \geq \lambda_{\sigma_{1}} \geq \lambda_{\sigma_{2}}$, the wait time is nonnegative, i.e., $\mu_{\sigma_{1}}- \lambda_{\sigma_{2}} \geq 0$.

Wait time of patient schedule for $P$ at $\sigma_{3}$: $\mu_{\sigma_{1}} + \mu_{\sigma_{2}}- \lambda_{\sigma_{2}}- \lambda_{\sigma_{3}}$. Since $\mu_{\sigma_{1}} \geq \lambda_{\sigma_{1}} \geq \lambda_{\sigma_{2}}$ and $\mu_{\sigma_{2}} \geq \lambda_{\sigma_{2}} \geq \lambda_{\sigma_{3}}$, the wait time is nonnegative, i.e., $\mu_{\sigma_{1}} + \mu_{\sigma_{2}}- \lambda_{\sigma_{2}}- \lambda_{\sigma_{3}} \geq 0$, so on finally,

Wait time of patient schedule for $P$ at $\sigma_{r-w}$: $\mu_{\sigma_{1}} + \mu_{\sigma_{2}} + \dots + \mu_{\sigma_{r-v-1}}- \lambda_{\sigma_{2}}- \lambda_{\sigma_{3}} - \dots - \lambda_{\sigma_{r-v}} \geq 0$.

Adding all wait times, we obtain $W_p(\sigma^f) = (r-v-1)[\mu_{\sigma_{1}}- \lambda_{\sigma_{2}}] +  (r-v-2)[\mu_{\sigma_{2}}- \lambda_{\sigma_{3}}] + \dots +
[\mu_{\sigma_{r-v-1}}- \lambda_{\sigma_{r-v}}] \geq 0$. By substituting
each $[\mu_{\sigma_{i}}- \lambda_{\sigma_{i+1}}]$ with
 larger value $[\gamma_1 - \gamma_2]$ in $W_p(\sigma^f)$, we obtain 
$W_p(\sigma^f) \leq [(r-v-1) + (r-v-2) + \dots + 1] [\gamma_1 - \gamma_2]$, That is, we have the following bound:

$$ W_p(\sigma^f) \leq  \frac{(r-v)(r-v-1)}{2}[\gamma_1 - \gamma_2].$$	

We now show this bound is tight. For a bound to be tight there should exist two positive coefficients $c_1$ and $c_2$ such that the following is satisfied:

$$ c_1 \frac{(r-v)(r-v-1)}{2}[\gamma_1 - \gamma_2] \leq W_p(\sigma^f) \leq  c_2 \frac{(r-v)(r-v-1)}{2}[\gamma_1 - \gamma_2].$$

As the total waiting time is characterized by $W_p(\sigma^f) = (r-v-1)[\mu_{\sigma_{1}}- \lambda_{\sigma_{2}}] +  (r-v-2)[\mu_{\sigma_{2}}- \lambda_{\sigma_{3}}] + \dots +
[\mu_{\sigma_{r-v-1}}- \lambda_{\sigma_{r-v}}]$, without loss of generality we can rewrite the expression as follows for some $a_1,a_2$, where $2\leq a_1 \leq r-v$ and $1\leq a_2 \leq r-v-1$:
\[ W_p(\sigma^f) = (r-v-1)[\mu_{\sigma_{1}}- \lambda_{\sigma_{2}}] +  (r-v-2)[\mu_{\sigma_{2}}- \lambda_{\sigma_{3}}] + \dots\]
\[
+ (r-v-a_1+1)[\gamma_1-\lambda_{a_1}] +\dots + (r-v-a_2)[\mu_{a_2}-\gamma_2]+\dots +
[\mu_{\sigma_{r-v-1}}- \lambda_{\sigma_{r-v}}].\] 

That is
\[ W_p(\sigma^f) = [ \gamma_1 - \gamma_2] + (r-v-1)[ \mu_{\sigma_1} - \lambda_{\sigma_2}] + (r-v-2)[ \mu_{\sigma_2} - \lambda_{\sigma_3}] + \dots \]
\[ + (r-v-a_1)\gamma_1-(r-v-a_1+1)\lambda_{a_1} +\dots + (r-v-a_2)\mu_{a_2}-(r-v-a_2-1)\gamma_2+\dots +
[\mu_{\sigma_{r-v-1}}- \lambda_{\sigma_{r-v}}].\]

Hence it follows that $W_p(\sigma^f) \geq \gamma_1 - \gamma_2$. Let $c_1=\frac{2}{(r-v)(r-v-1)}$ and $c_2=1$, then it follows that:

$$ c_1 \frac{(r-v)(r-v-1)}{2} [\gamma_1 - \gamma_2] \leq W_p(\sigma^f) \leq  c_2 \frac{(r-v)(r-v-1)}{2} [\gamma_1 - \gamma_2],$$

$$ [\gamma_1 - \gamma_2] \leq W_p(\sigma^f) \leq \frac{(r-v)(r-v-1)}{2}[\gamma_1 - \gamma_2].$$

Thus, the bound is tight. $\blacksquare$

\section{Section~\ref{msmp} Proofs} \label{app1.1a}

\noindent \textbf{Proof of Lemma~\ref{Perf-bound-block-muti}:}

Note that Algorithm~\ref{alg-RMA} is an extension of Algorithm~\ref{alg-MA} for multiple blocks. Note that Algorithm~\ref{alg-MA} provides a patient sequence $\sigma$ for a block $\pi$, where $\sigma=\{\sigma_{1}, \sigma_{2}, \ldots, \sigma_{r-v}, \sigma_{r-v+1}, \ldots, \sigma_{r}\}$. Recall the two partial sequences for a block, $\sigma^f$ and $\sigma^b$ of $\sigma$, where $\sigma= \sigma^f \cup \sigma^b$, $\sigma^f =\{\sigma_{1}, \sigma_{2}, \ldots, \sigma_{r-v}\}$, and $\sigma^b=\{\sigma_{r-v+1}, \sigma_{r-v+2}, \ldots, \sigma_{r}\}$. The patient waiting time occurs only at $\sigma^f$. Note that for a block, $W_a(\sigma^f)=0$,  $W_p(\sigma^b)=0$ and $W_a(\sigma^b)=0$. Moreover, in Lemma~\ref{Perf-bound-block1}, we show for a block that $W_p(\sigma^f) \leq  \frac{(r-v)(r-v-1)}{2}[\gamma_1 - \gamma_2]$. Consider a multiple block schedule, $\hat{\Gamma}$ as shown in Figure~\ref{figx:block-multi}: the schedule consists of $k$ jobs (corresponding $k$ blocks), each job having a bead $H$, a body $B$, and a tail $T$, where  $B=\sum_{j=2}^{r} \lambda_{\sigma_{j}}$, $B+T=\sum_{j=1}^{r-v} \mu_{\sigma_{j}}$, and $H=\lambda_{\sigma_{1}}$ as defined in Section~\ref{bs}. In Figure~\ref{figx:block-multi}, jobs are scheduled no-wait in manner, where $\theta  =\max \{0, \, (B+T) - (B+H)\} =  \max \{0, \, [\,\sum_{j=1}^{r-v} \mu_{\sigma_{j}} - \sum_{j=1}^{r} \lambda_{\sigma_{j}}]\}$. 
It has to be noted that the patient wait times remain the same in each block in schedule $\hat{\Gamma}$ shown in Figure~\ref{figx:block-multi} since all repetitively following blocks (jobs) are scheduled as no-wait in manner. Thus, the bound, $W_p(\sigma^f) \leq  \frac{(r-v)(r-v-1)}{2}[\gamma_1 - \gamma_2]$ is also valid for all repetitively scheduled blocks in Figure~\ref{figx:block-multi}.

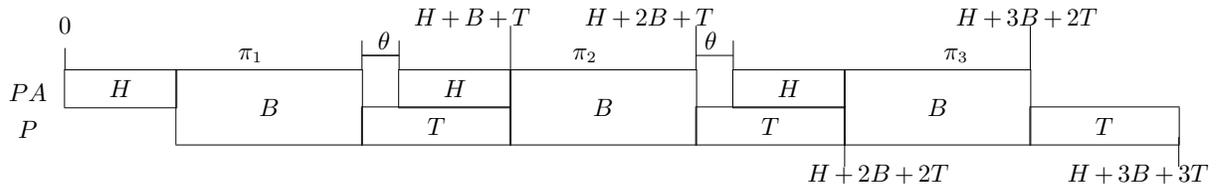
\begin{figure}[htbp]
\resizebox{\textwidth}{!}{
	{\normalsize
		\setlength{\unitlength}{1.15mm}
		\begin{picture}(165,15)(10,14) 
			%\thicklines	
			\put(15,17){\makebox(0,0){$PA$}}
			
			\put(15,12){\makebox(0,0){$P$}}
			
			\multiput(20,20)(27,0){1}{\line(0,1){3}}
			\put(20,26){\makebox(0,0){0}}

			%======Head1 ====
			
			\put(20,15){\framebox(15,5){$H$}}

			%===========Body1=============
			
			\put(35,10){\framebox(25,10){$B$}}
			
			%===========Tail1=============
			
			\put(60,10){\framebox(20,5){$T$}}

			%======Head2 ====
			
			\put(65,15){\framebox(15,5){$H$}}

			%===========Body2=============
			
			\put(80,10){\framebox(25,10){$B$}}
			
			%===========Tail2=============
			
			\put(105,10){\framebox(20,5){$T$}}

			%======Head3 ====
			
			\put(110,15){\framebox(15,5){$H$}}

			%===========Body3=============
			
			\put(125,10){\framebox(25,10){$B$}}
			
			%===========Tail3=============
			
			\put(150,10){\framebox(20,5){$T$}}

			%======Othres ====
			
			%\multiput(64.5,20)(27,5){1}{\line(0,1){4}}
			%\multiput(60.5,20)(27,5){1}{\line(0,1){4}}

            \multiput(65,20)(27,5){1}{\line(0,1){4}}
			\multiput(60,20)(27,5){1}{\line(0,1){4}}   
			\put(63,24){\makebox(0,0){$\theta$}}
				\multiput(60,22)(27,5){1}{\line(1,0){5}}
				\multiput(60,22)(27,5){1}{\line(1,0){5}}
			\multiput(110,20)(27,5){1}{\line(0,1){4}}
			\put(107,24){\makebox(0,0){$\theta$}}
			\multiput(105,22)(27,5){1}{\line(1,0){5}}
			%	\put(60,20)(27,5){\line(1,0){4}}

			\multiput(105,20)(27,5){1}{\line(0,1){6}}
			\put(90,26){$H+2B+T$}
			
			\multiput(80,20)(27,5){1}{\line(0,1){6}}
			\put(67,26){$H+B+T$}
			
			\multiput(150,20)(27,5){1}{\line(0,1){6}}
			\put(140,26){$H+3B+2T$}
			
			\multiput(125,7)(27,5){1}{\line(0,1){4}}
			\put(120,5){$H+2B+2T$}
			
			\multiput(170,7)(27,5){1}{\line(0,1){4}}
			\put(155,5){$H+3B+3T$}

			\put(45,22){\makebox(0,0){$\pi_1$}}
			\put(90,22){\makebox(0,0){$\pi_2$}}
			\put(140,22){\makebox(0,0){$\pi_3$}}
			
	\end{picture}}}	
	\vspace{0.2in}	
	\caption{Block schedule $\hat{\Gamma}= \{\pi_1\pi_2\pi_3\}$ is a concatenation of $k$ blocks as shown here, where $k=3$ and $\sigma$ is patient sequence in $\pi_i$, $i=1,2, \ldots, k$.}
	\label{figx:block-multi}	
\end{figure}

\medskip

Algorithm~\ref{alg-RMA} obtains schedule $\Gamma$ by adjusting $\hat{\Gamma}$ as follows: It schedules patients in block $\pi_2$ by advancing their appointment times by $\theta$ with respect to the schedule $\hat{\Gamma}$ in Figure~\ref{figx:block-multi}, patients in block $\pi_3$ by advancing their appointment time by $2\theta$, so on finally until patients in block $\pi_k$ by advancing their appointment times by $(k-1)\theta$. Consequently, in each block, only $(r-v)$ patients of $Q^+$ type wait time in each block increase as follows: in block $\pi_2$ $Q^+$ patients wait times increase by $(r-v)\theta$, in block $\pi_3$ $Q^+$ patients wait times increase by $2(r-v)\theta$, so on finally until in block $\pi_k$ $Q^+$ patients wait times increase by $ (k-1)(r-v)\theta$.

Thus, the total patient waiting time $W$ of the schedule $\Gamma$ given by Algorithm~\ref{alg-RMA} is bounded as follows:

$$W \leq k W_p(\sigma^f) + (r-v)\theta + 2(r-v)\theta + \ldots  + (k-1)(r-v)\theta.$$

That is,

$$W \leq \frac{k(r-v)(r-v-1)}{2}[\gamma_1 - \gamma_2] + {\frac{k(k-1)(r-v)}{2} }\theta. \,\,\,\, \blacksquare$$

\section{Section~\ref{ucd} Proofs} \label{app1.1}

\noindent \textbf{Proof of Lemma~\ref{L-block1_st}:}
	
Let Algorithm~\ref{alg-MA} provide a sequence, $\sigma$, where
$\sigma=\{\sigma_{1}, \sigma_{2}, \ldots, \sigma_{r-v}, \sigma_{r-v+1}, \ldots, \sigma_{r}\}$. Consider two partial sequences, $\sigma^f$ and $\sigma^b$ of $\sigma$, where $\sigma= \sigma^f \cup \sigma^b$. We prove Lemma~\ref{L-block1_st} by induction. Note that $\pi$ is the block schedule formed by sequence $\sigma$.

\medskip

Let $\lambda_{\sigma_i}$ and $\mu_{\sigma_i}$ denote the processing time of patients scheduled in $i^{th}$ position, $\sigma_i$ in schedule $\pi$. Note that the first patient's processing time $\lambda_{\sigma_1}$ in $\pi$ may be reduced to the  $\lambda_{\sigma_1} - {\frac{w}{2}}\lambda_{\sigma_1}$ in the worst case. Thus, it is sufficient to set the appointment times of patients scheduled in positions,  $\sigma_1$ and $\sigma_2$ to $\tau^{\sigma_1}_a=0$ and $\tau^{\sigma_2}_a=\lambda_{\sigma_1} - {\frac{w}{2}}\lambda_{\sigma_1}$, respectively. As the processing time of the first patient can at least be $\lambda_{\sigma_1} - { \frac{w}{2}}\lambda_{\sigma_1}$ this appointment time assignment ensures zero $PA$ idle time. Thus, the base case is shown. Similarly, the first two patients' processing time $\lambda_{\sigma_1}+ \lambda_{\sigma_2}$ in $\pi$ may be reduced to the  $\lambda_{\sigma_1} +\lambda_{\sigma_2} - {\frac{w}{2}}[\lambda_{\sigma_1} +\lambda_{\sigma_2}]$ in the worst case. Thus, it is sufficient to set the appointment time of patient scheduled in position,  $\sigma_3$ to $\tau^{\sigma_3}_a= \lambda_{\sigma_1} +\lambda_{\sigma_2} - {\frac{w}{2}}[\lambda_{\sigma_1} +\lambda_{\sigma_2}]$. 

Assume that the appointment times of first $i$ patients are defined as $\tau^{\sigma_{i}}_a = 
(1-\frac{w}{2}) \sum_{k=1}^{i-1} \lambda_{\sigma_k}$, for some $i \in \{ 2, 3, \ldots, r \}$ for the inductive step and $\tau^{\sigma_1}_a=0$, as shown in the base case. Therefore, no $PA$ idle time is incurred by the first $i$ patients. Then, the processing time of the first $i$ patients, can at worst be reduced to $\lambda_{\sigma_i} - {\frac{w}{2}}\lambda_{\sigma_i} + (1-\frac{w}{2}) \sum_{k=1}^{i-1} \lambda_{\sigma_k}$ which is equivalent to $(1-\frac{w}{2}) \sum_{k=1}^{i} \lambda_{\sigma_k}$. Setting the appointment time of the patient in position $\sigma_{i+1}$ to $\tau^{\sigma_{i+1}}_a = (1-\frac{w}{2}) \sum_{k=1}^{i} \lambda_{\sigma_k}$ would yield no idle time at $PA$ for the first $i+1$ patients. As this holds for any $i+1 \in \{ 2, 3, \ldots, r \}$, in $\pi$ with the above appointment times of the patients, there is no idle time at $PA$ for any realized sample path.

% By continuing this argument, we find that $\tau^{\sigma_1}_a=0$ and $\tau^{\sigma_{i}}_a = 
% (1-{w \over 2}) \sum_{k=1}^{i-1} \lambda_{\sigma_k}$, $i=2, 3, \ldots, r$. It is evident that in $\sigma$ with above appointment times of the patients, there is no idle time at $PA$ for any realized sample path.

We now show that $\sigma^f$ does not have any idle time at $P$ for any realized sample path under the condition in the lemma that  $w \leq \min_{1 \leq i \leq r-v-1} \left[\frac{2\sum_{k=1}^{i}(\mu_{\sigma_k} - \lambda_{\sigma_{k+1}})} {\sum_{k=1}^{i}(\mu_{\sigma_k} + \lambda_{\sigma_{k+1}})}\right]$. Note that $\mu_{\sigma_1} \geq \lambda_{\sigma_1} \geq \lambda_{\sigma_2}$. Thus, the patient wait time scheduled at position $\sigma_{1}$ before the treatment by the $P$ in the deterministic case is $(\mu_{\sigma_1} - \lambda_{\sigma_2}) \geq 0$. To maintain no idle at $P$, we must have 
$(\mu_{\sigma_1} - \lambda_{\sigma_2}) \geq {\frac{w}{2}} (\mu_{\sigma_1} + \lambda_{\sigma_2})$.
That is $w \leq \frac{2(\mu_{\sigma_1} - \lambda_{\sigma_2})}{\mu_{\sigma_1} + \lambda_{\sigma_2}}$. Thus, the base case is shown. Similarly, the  wait time of  patient  in $\sigma_{2}$ before the treatment by the $P$ in the deterministic case is $(\mu_{\sigma_1} + \mu_{\sigma_2}- \lambda_{\sigma_2} - \lambda_{\sigma_3}) \geq 0$ since $\mu_{\sigma_1} \geq \lambda_{\sigma_1} \geq \lambda_{\sigma_2}$ and $\mu_{\sigma_2} \geq \lambda_{\sigma_2} \geq \lambda_{\sigma_3}$. To maintain no idle at $P$, we must have  
$(\mu_{\sigma_1} + \mu_{\sigma_2}- \lambda_{\sigma_2} - \lambda_{\sigma_3}) \geq {\frac{w}{2}} (\mu_{\sigma_1} + \mu_{\sigma_2} + \lambda_{\sigma_2} + \lambda_{\sigma_3})$, i.e., $w \leq \frac{2(\mu_{\sigma_1} - \lambda_{\sigma_2}) + 2(\mu_{\sigma_2} - \lambda_{\sigma_3})}{(\mu_{\sigma_1} + \lambda_{\sigma_2}) + (\mu_{\sigma_2} + \lambda_{\sigma_3})}$.

Suppose $w \leq \frac{2\sum_{k=1}^{i}(\mu_{\sigma_k} - \lambda_{\sigma_{k+1}})} {\sum_{k=1}^{i}(\mu_{\sigma_k} + \lambda_{\sigma_{k+1}})}$ for some $i \in \{ 1,2, \cdots , r \}$ for the inductive step. Then, the wait time of the patient in $\sigma_{i+1}$ before the treatment by the $P$ in the deterministic case is $(\mu_{\sigma_{i+1}} + \sum_{k=1}^{i} \mu_{\sigma_k} - \lambda_{\sigma_{i+1}} - \sum_{k=1}^{i+1} \lambda_{\sigma_k}) \geq 0$ since $\mu_{\sigma_1} \geq \lambda_{\sigma_1} \geq \lambda_{\sigma_2} \geq \cdots \geq \lambda_{\sigma_{i+1}}$ and $\mu_{\sigma_2} \geq \lambda_{\sigma_2} \geq \cdots \geq \lambda_{\sigma_{i+1}}$. To maintain no idle at $P$, we must have $\frac{w}{2}(\mu_{\sigma_{i+1}}+\lambda_{\sigma_{i+2}}) + \frac{w}{2}  \sum_{k=1}^{i}(\mu_{\sigma_k} + \lambda_{\sigma_{k+1}}) \leq (\mu_{\sigma_{i+1}}-\lambda_{\sigma_{i+2}}) + \sum_{k=1}^{i}(\mu_{\sigma_k} - \lambda_{\sigma_{k+1}})$, which can equivalently be expressed as $w \leq \frac{2\sum_{k=1}^{i+1}(\mu_{\sigma_k} - \lambda_{\sigma_{k+1}})}{\sum_{k=1}^{i+1}(\mu_{\sigma_k} + \lambda_{\sigma_{k+1}})}$. As this holds for any $i+1 \in \{ 1, 2, \ldots, r \}$, the claim holds. Since $\sigma^b$ does not involve $P$, the result follows. $\blacksquare$

% In generalizing this expression for any patient in $\sigma_{i}$, we get 
% $w \leq \min_{1 \leq i \leq r-v-1} \left[{2\sum_{k=1}^{i}(\mu_{\sigma_k} - \lambda_{\sigma_{k+1}}) \over \sum_{k=1}^{i}(\mu_{\sigma_k} + \lambda_{\sigma_{k+1}})} \right]$. Since $\sigma^b$ does not involve $P$, the result follows. $\blacksquare$

\section{Deterministic Service Times Computational Results} \label{app2}

In this section, we test our Planning Horizon Model, introduced in Section~\ref{msmp}, on a variety of randomly generated problem sets using real-world cost parameters which we discuss in Section~\ref{cdp}. We compare our heuristic solutions to the optimal value and discuss the related insights.

% Similarly, Tables~\ref{resres1.1}, \ref{resres1.3} and \ref{resres1.2} present the average computational results for the Planning Horizon Model, Algorithm~\ref{alg-RMA} and Algorithm~\ref{alg-RB} for a varying set of overtime costs for 300 minutes of regular time in a day ($R$). 

\begin{table}[]
 \caption{Model Results (in minutes) for 100 randomized service times, where the patient waiting time cost, $\alpha=1$, the idle time costs for $PA$ and $P$, $\beta_a=\beta_p=5$, and regular time in the day, $R=200$ minutes}
    \centering
    \resizebox{\textwidth}{!}{
    \begin{tabular}{|c|c|c|c|c|c|c|}
     \hline
         & \multicolumn{6}{c|}{Planning Horizon Model}\\
         \hline
         $o_a,o_p$ & obj & $PA$ Overtime & $P$ Overtime & $PA$ Idle Time & $P$ Idle Time & Waiting Time\\
         \hline
         1  & 148.03  & 38.19 & 80.35 & 0.03 & 0.00    & 29.34 \\
         5  & 620.43  & 38.18 & 79.79 & 0.02 & 0.01 & 30.43 \\
         10 & 1210.08 & 38.18 & 79.72 & 0.02 & 0.00    & 30.98 \\
         30 & 3567.91 & 38.18 & 79.71 & 0.02 & 0.00    & 31.11 \\
         \hline
    \end{tabular}
    }
    % \vspace{0.01in}
   
    \label{resres2.1}
\end{table}

\begin{table}[]
    \caption{Algorithm~\ref{alg-RMA} Results (in minutes) for 100 randomized service times, where the patient waiting time cost, $\alpha=1$, the idle time costs for $PA$ and $P$, $\beta_a=\beta_p=5$, and regular time in the day, $R=200$ minutes}
    \centering
    \resizebox{\textwidth}{!}{
    \begin{tabular}{|c|c|c|c|c|c|c|}
     \hline
         & \multicolumn{6}{c|}{Algorithm~\ref{alg-RMA}, Two-Stage Scheduling for Planning Horizon}\\
         \hline
         $o_a,o_p$ & obj & $PA$ Overtime & $P$ Overtime & $PA$ Idle Time & $P$ Idle Time & Waiting Time \\
         \hline
         1  & 408.67  & 51.53 & 85.91 & 14.23 & 0.00 & 200.08 \\
         5  & 958.43  & 51.53 & 85.91 & 14.23 & 0.00 & 200.08 \\
         10 & 1645.63 & 51.53 & 85.91 & 14.23 & 0.00 & 200.08 \\
         30 & 4394.43 & 51.53 & 85.91 & 14.23 & 0.00 & 200.08\\
         \hline
    \end{tabular}
    }
    % \vspace{0.01in}

    \label{resres2.3}
\end{table}

\begin{table}[]
    \caption{Algorithm~\ref{alg-RB} Results (in minutes) for randomized service times, where the patient waiting time cost, $\alpha=1$, the idle time costs for $PA$ and $P$, $\beta_a=\beta_p=5$, and regular time in the day, $R=200$ minutes}
    \centering
    \resizebox{\textwidth}{!}{
    \begin{tabular}{|c|c|c|c|c|c|c|}
     \hline
         & \multicolumn{6}{c|}{Algorithm~\ref{alg-RB}, Improved Two-Stage Scheduling for Planning Horizon}\\
         \hline
         $o_a,o_p$ & obj & $PA$ Overtime & $P$ Overtime & $PA$ Idle Time & $P$ Idle Time & Waiting Time \\
         \hline
         1  & 270.75  & 51.53 & 85.91 & 14.23 & 0.00 & 62.16 \\
         5  & 820.51  & 51.53 & 85.91 & 14.23 & 0.00 & 62.16 \\
         10 & 1507.71 & 51.53 & 85.91 & 14.23 & 0.00 & 62.16 \\
         30 & 4256.51 & 51.53 & 85.91 & 14.23 & 0.00 & 62.16\\
         \hline
    \end{tabular}
    }
    % \vspace{0.01in}

    \label{resres2.2}
\end{table}

 Each randomly generated problem set was tested using the Planning Horizon Model, Algorithm~\ref{alg-RMA} and Algorithm~\ref{alg-RB}, with results in Tables~\ref{resres2.1}, \ref{resres2.3}, and \ref{resres2.2} respectively. We varied the overtime costs ($o_a,o_p$) from 1 to 30 similar to \cite{klassen2019appointment}, while keeping other parameters stable. We set the length of a planning horizon ($R$) as 200 minutes.
 Note that the Planning Horison Model objective function values change a lot and the minutes of overtime, idle time and wait time change a little with respect to different cost parameters across random sample paths for each iteration. The heuristic algorithm objective function values also change a lot with respect to different cost parameters across the same random sample paths. However, the heuristic algorithm minutes of overtime, idle time and wait time decisions do not change with respect to different cost parameters (and using the same random sample paths) since the algorithm only looks at patient types to fill the schedule. Thus, the patient waiting time, overtime, and idle times for $PA$ and $P$ remain the same for changing overtime costs in Tables~\ref{resres2.3} and~\ref{resres2.2}. 
 Algorithms~\ref{alg-RMA}, \ref{alg-RB}, and the Planning Horizon Model yield very similar results for $P$ idle time. The algorithms provide no-idle time schedules for $P$ and the Planning Horizon Model yields a near zero $P$ idle time.

 Scheduling all $Q^+$ Group patients in a non-increasing order of their service times, Algorithm~\ref{alg-RMA} yields a block schedule with no $P$ idle time, but includes high patient waiting time due to greater $P$ service times compared to those of $PA$. By alternating between $Q^+$ and $Q$ Group patients in the schedule, Algorithm~\ref{alg-RB} manages to reduce this patient waiting time by nearly 70\%, while maintaining the same values for $P$ and $PA$ idle time and overtime, as illustrated in Tables~\ref{resres2.3} and \ref{resres2.2}. This reduction in patient waiting time while every other metric is similar illustrates that Algorithm~\ref{alg-RB} provides an appointment template with a smaller total cost compared to Algorithm~\ref{alg-RMA}.

Even though Algorithm~\ref{alg-RB} forms a no-idle time schedule within each block, repeating each block to form the planning horizon appointment template can create additional idle time for $PA$, which then results in more overtime. However, the Planning Horizon Model is able to generate different block sequences for different blocks of the day and therefore optimally reduce the idle time incurred between blocks. On the other hand, Algorithm~\ref{alg-RB} provides a solution that is easily implementable (i.e. same sequence for all blocks), with comparable $P$ idle time and overtime to the optimal model solution and only increases $PA$ idle time and overtime by roughly 14 minutes and patient waiting time by 32 minutes, that corresponds to 6.675$\%$ and 16$\%$ difference respectively. In a real-world scenario, small idle times are often used to complete other tasks and can be less concerning if predefined tasks are available to fill that time \citep{klassen2019appointment}. Idle time can also prevent a disruption in the overall daily schedule by establishing a buffer between appointments under service time uncertainty as studied in Section~\ref{ucd}.     

\begin{table}[h]
\caption{Different Block Patterns}
    \centering
    % \resizebox{\textwidth}{!}{
    \begin{tabular}{|c|ccccccccc|}
    \hline
         Pattern & \multicolumn{9}{|c|}{Patient Types}\\
         \hline
         A1 & T4 & \textcolor{blue}{\textbf{\textit{T2}}} & T3 & \textcolor{blue}{\textbf{\textit{T1}}} & T4 & \textcolor{blue}{\textbf{\textit{T1}}} & \textcolor{blue}{\textbf{\textit{T1}}} & T4 & \textcolor{blue}{\textbf{\textit{T2}}}\\
         A2 & T4 & \textcolor{blue}{\textbf{\textit{T1}}} & \textcolor{blue}{\textbf{\textit{T1}}} & T4 & \textcolor{blue}{\textbf{\textit{T2}}} & T3 & \textcolor{blue}{\textbf{\textit{T1}}} & T4 & \textcolor{blue}{\textbf{\textit{T2}}}\\
         A3 & T4 & \textcolor{blue}{\textbf{\textit{T1}}} & \textcolor{blue}{\textbf{\textit{T1}}} & T4 & \textcolor{blue}{\textbf{\textit{T1}}} & \textcolor{blue}{\textbf{\textit{T2}}} & T4 & \textcolor{blue}{\textbf{\textit{T2}}} & T3\\
         B1 & T3 & \textcolor{blue}{\textbf{\textit{T1}}} & T4 & \textcolor{blue}{\textbf{\textit{T2}}} & T4 & \textcolor{blue}{\textbf{\textit{T1}}} & \textcolor{blue}{\textbf{\textit{T2}}} & T4 & \textcolor{blue}{\textbf{\textit{T1}}}\\
         B2 & T4 & \textcolor{blue}{\textbf{\textit{T2}}} & \textcolor{blue}{\textbf{\textit{T1}}} & T4 & \textcolor{blue}{\textbf{\textit{T1}}} & \textcolor{blue}{\textbf{\textit{T1}}} & T4 & \textcolor{blue}{\textbf{\textit{T2}}} & T3\\
         B3 & T4 & \textcolor{blue}{\textbf{\textit{T1}}} & \textcolor{blue}{\textbf{\textit{T2}}} & T4 & \textcolor{blue}{\textbf{\textit{T2}}} & T3 & \textcolor{blue}{\textbf{\textit{T1}}} & T4 & \textcolor{blue}{\textbf{\textit{T1}}}\\
         B4 & T3 & \textcolor{blue}{\textbf{\textit{T1}}} & T4 & \textcolor{blue}{\textbf{\textit{T1}}} & \textcolor{blue}{\textbf{\textit{T1}}} & T4 & \textcolor{blue}{\textbf{\textit{T2}}} & T4 & \textcolor{blue}{\textbf{\textit{T2}}}\\
         B5 & \textcolor{blue}{\textbf{\textit{T2}}} & T4 & \textcolor{blue}{\textbf{\textit{T1}}} & \textcolor{blue}{\textbf{\textit{T1}}} & T4 & \textcolor{blue}{\textbf{\textit{T2}}} & T3 & \textcolor{blue}{\textbf{\textit{T1}}} & T4\\
         Algorithm~\ref{alg-RMA} & T3 & T4 & T4 & T4 & \textcolor{blue}{\textbf{\textit{T1}}} & \textcolor{blue}{\textbf{\textit{T1}}} & \textcolor{blue}{\textbf{\textit{T1}}} & \textcolor{blue}{\textbf{\textit{T2}}} & \textcolor{blue}{\textbf{\textit{T2}}}\\
         Algorithm~\ref{alg-RB} & T3 & \textcolor{blue}{\textbf{\textit{T1}}} & T4 & \textcolor{blue}{\textbf{\textit{T1}}} & \textcolor{blue}{\textbf{\textit{T1}}} & T4 & \textcolor{blue}{\textbf{\textit{T2}}} & T4 & \textcolor{blue}{\textbf{\textit{T2}}}\\
         \hline
    \end{tabular}
    % }
    % \vspace{0.01in}
    
    \label{res1.3}
\end{table}

% In practice, assigning costs for different performance metrics can be challenging for an outpatient clinic, this can render an already time consuming process of solving a mathematical model to design a daily appointment template even more impractical.

 Table~\ref{res1.3} presents different block patterns provided by the Planning Horizon Model (A1-B5), and Algorithms~\ref{alg-RMA} and \ref{alg-RB}. Contrary to heuristic Algorithms~\ref{alg-RMA} and \ref{alg-RB}, the Planning Horizon Model can yield different block schedules within a planning horizon or between parameter sets. We observe that for all 100 randomized parameter sets, there exist three different patterns of patient types in the first block (A1,A2,A3) and five different patterns in the second block (B1,B2,B3,B4,B5). One can observe that most of the optimal Planning Horizon Model patterns alternate between $Q$ (shown in bold and italics) and $Q^+$ Group Patients at consecutive positions in the block. This observation shows that the optimal solution provided by the Planning Horizon Model and the appointment template provided by Algorithm~\ref{alg-RB} adopt a similar approach. In fact, the pattern obtained by Algorithm~\ref{alg-RB} is identical to the pattern (B4) that is utilized by the Planning Horizon Model for the second block when the patient waiting time cost and overtime costs for $P$ and $PA$ are all equal to 1 and the idle time costs for $P$ and $PA$ are set equal to 0.

 \section{Service Time Uncertainty Models and Procedures} \label{app3}

\noindent \textbf{Stochastic Programming Model for a Single Block}
\begin{align}
    \min\; & \Psi_s=\sum_{s=1}^N \rho_s \Pi_s(w_a^{j,s},\bar{w}_p^{j,s},d_a^{j,s},d_p^{j,s}) \label{m17}\\
    \text{subject to }  & \Pi_s=\alpha \sum_{j=1}^r w^{j,s}_a + \alpha \sum_{j=1}^r \bar{w}^{j,s}_p+\beta_a \sum_{j=1}^r d^{j,s}_a+\beta_p \sum_{j=1}^r d^{j,s}_p\hspace{3mm}s=1,2,\cdots,N\label{m29}\\
    & \sum_{j=1}^r x_{lj}=1\hspace{3mm} l=1,2,\cdots,r\label{m18}\\
    & \sum_{l=1}^{r}  x_{lj}=1\hspace{3mm}  j=1,2,\cdots,r\label{m19}\\
    & \sum_{z\in\{1,2,\cdots,q\}} x_{z1}=0\label{m36}\\
    & e_a^{j+1,s}\geq e_a^{j,s} + \sum_{l=1}^r x_{lj} \lambda_{l,s} \hspace{3mm}  j=1,2,\cdots,r-1;\hspace{3mm}s=1,2,\cdots,N\label{m20}\\
    & e_p^{j,s} \geq e_a^{j,s} + \sum_{l=1}^r x_{lj} \lambda_{l,s} \hspace{3mm} j=1,2,\cdots,r;\hspace{3mm}s=1,2,\cdots,N\label{m21}\\
    & e_p^{j+1,s} \geq e_p^{j,s} + \sum_{l=1}^r x_{lj} \mu_{l,s} \hspace{3mm} j=1,2,\cdots,r-1;\hspace{3mm}s=1,2,\cdots,N\label{m22}\\
    & d_a^{j+1,s} \geq e_a^{j+1,s}-e_a^{j,s}-\sum_{l=1}^r x_{lj} \lambda_l \hspace{3mm} j=1,2,\cdots,r-1;\hspace{3mm}s=1,2,\cdots,N\label{m23}\\
    & d_p^{j+1,s} \geq e_p^{j+1,s}-e_p^{j,s}-\sum_{l=1}^r x_{lj}\mu_l \hspace{3mm}  j=1,2,\cdots,r-1;\hspace{3mm}s=1,2,\cdots,N\label{m24}\\
    & e_a^{1,s}=0\label{m25}\\
    % & w_a^{j,s} \geq e_a^{j,s}-e_a^{j-1,s}-\sum_{l=1}^r x_{l,j-1}\lambda_{l,s}\hspace{3mm} j=2,\cdots,r;\hspace{3mm}s=1,2,\cdots,N\\
    & \tau_a^j \leq e_a^{j,s}\hspace{3mm} j=2,\cdots,r;\hspace{3mm}s=1,2,\cdots,N\\
    & w_a^{j,s} \geq e_a^{j,s}-\tau_a^j \hspace{3mm} j=2,\cdots,r;\hspace{3mm}s=1,2,\cdots,N \\
    & \bar{w}_p^{j,s} \geq e_p^{j,s}-e_a^{j,s}-\sum_{l=1}^r x_{lj} \lambda_{l,s} -M(1-y_j)\hspace{3mm} j=1,2,\cdots,r;\hspace{3mm}s=1,2,\cdots,N\label{m26}\\
    & \frac{y_j}{M} \leq \sum_{l=1}^r \mu_{l,s} x_{lj} \hspace{3mm} j=1,2,\cdots,r;\hspace{3mm}s=1,2,\cdots,N\label{m27}\\
    & My_j \geq \sum_{l=1}^r \mu_{l,s} x_{lj}\hspace{3mm} j=1,2,\cdots,r;\hspace{3mm}s=1,2,\cdots,N\label{m28}\\
    % & \Pi_s \leq (1+\delta)\Pi^*_s\hspace{3mm}s=1,2,\cdots,N\label{m30}\\
    & x_{lj}\in \{0,1\}\hspace{3mm} l=1,2,\cdots,r;\hspace{3mm} j=1,2,\cdots,r\label{m31}\\
    & y_j\in \{0,1\}\hspace{3mm} j=1,2,\cdots,r\label{m32}\\
    & d_a^{j,s}, d_p^{j,s}, e_a^{j,s}, e_p^{j,s}, w_a^{j,s},\bar{w}_p^{j,s}, \tau_a^j \geq 0\hspace{3mm}   j=1,2,\cdots,r;\hspace{3mm}s=1,2,\cdots,N\label{m33}\\
    & M\text{: sufficiently large number} \nonumber
\end{align}

 \newpage

\textbf{Stochastic Programming Model for the Planning Horizon}
\vspace{-1mm}
\begin{align}
    \min\; & \Psi^p_s=\sum_{s=1}^N \rho_s\Pi^p_s(w_a^{t,s},\bar{w}^{t,s}_p, d^{t,s}_a, d^{t,s}_p, b_a^{s},b_p^{s})\\
    % \min\; &  U_p(\sigma)+ I_a(\sigma) + I_p(\sigma)+O_a(\sigma)+O_p(\sigma) \label{mm1}\\
    \text{subject to } & \sum_{t=(c-1)r+1}^{cr} x_{lct}=1\hspace{3mm} l=1,2,\cdots,r\hspace{3mm} c=1,2, \cdots,k \label{mmm2}\\
    & \sum_{c=1}^{k}\sum_{l=1}^{r}  x_{lct}=1\hspace{3mm}   t=1,2,\cdots,n \label{mmm3}\\
    & \Pi^p_s= \alpha \sum_{t=1}^n w^{t,s}_a + \alpha \sum_{t=1}^n \bar{w}^{t,s}_p+\beta_a \sum_{t=1}^n d^{t,s}_a+\beta_p \sum_{t=1}^n d^{t,s}_p + o_a b_a^{s}+ o_p b_p^{s}\hspace{3mm}s=1,2,\cdots,N\\
    & \sum_{z\in\{1,2,\cdots,q\}} x_{z11}=0\label{mmm4}\\
    & e_a^{t+1,s}\geq e_a^{t,s} + \sum_{c=1}^{k}\sum_{l=1}^{r} x_{lct} \lambda_l^{s} \hspace{3mm}  t=1,2,\cdots,n-1;\hspace{3mm}s=1,2,\cdots,N\label{mmm5}\\
    & e_p^{t,s} \geq e_a^{t,s} + \sum_{c=1}^{k}\sum_{l=1}^{r} x_{lct} \lambda_l^s \hspace{3mm} t=1,2,\cdots,n;\hspace{3mm}s=1,2,\cdots,N\label{mmm6}\\
    & e_p^{t+1,s} \geq e_p^{t,s} + \sum_{c=1}^{k}\sum_{l=1}^{r} x_{lct} \mu_l^s \hspace{3mm} t=1,2,\cdots,n-1;\hspace{3mm}s=1,2,\cdots,N\label{mmm7}\\
    & d_a^{t+1,s} \geq e_a^{t+1,s}-e_a^{t,s}- \sum_{c=1}^{k} \sum_{l=1}^r x_{lct} \lambda_l^s \hspace{3mm} t=1,2,\cdots,n-1;\hspace{3mm}s=1,2,\cdots,N\label{mmm8}\\
    & d_p^{t+1,s} \geq e_p^{t+1,s}-e_p^{t,s}- \sum_{c=1}^{k} \sum_{l=1}^r x_{lct}\mu_l^s \hspace{3mm}  t=1,2,\cdots,n-1;\hspace{3mm}s=1,2,\cdots,N\label{mmm9}\\
    & e_a^{1,s}=0;\hspace{3mm}s=1,2,\cdots,N\label{mmm10}\\
    % & w_a^{t,s} \geq e_a^{t,s}-e_a^{t-1,s}-\sum_{l=1}^r x_{l,t-1}\lambda_{l}^s\hspace{3mm} t=2,\cdots,r;\hspace{3mm}s=1,2,\cdots,N\\
    & \tau_a^t \leq e_a^{t,s}\hspace{3mm} t=2,\cdots,n;\hspace{3mm}s=1,2,\cdots,N\\
    & w_a^{t,s} \geq e_a^{t,s}-\tau_a^t \hspace{3mm} t=2,\cdots,n;\hspace{3mm}s=1,2,\cdots,N \\
    & \bar{w}_p^{t,s} \geq e_p^{t,s}-e_a^{t,s}-\sum_{c=1}^k\sum_{l=1}^r x_{lct} \lambda_l^s -M(1-y_t)\hspace{3mm} t=1,2,\cdots,n;\hspace{3mm}s=1,2,\cdots,N\label{mmm12}\\
    & \frac{y_t}{M} \leq \sum_{c=1}^k \sum_{l=1}^r \mu_l^s x_{lct} \hspace{3mm} t=1,2,\cdots,n;\hspace{3mm}s=1,2,\cdots,N\label{mmm13}\\
    & My_t \geq \sum_{c=1}^k \sum_{l=1}^r \mu_l^s x_{lct}\hspace{3mm} t=1,2,\cdots,n;\hspace{3mm}s=1,2,\cdots,N\label{mmm14}\\
    & b_a^{s} \geq e_a^{1,s} + \sum_{c=1}^k \sum_{t=(c-1)r+1}^{cr} \sum_{l=1}^r \lambda_l^{s} x_{lct} + \sum_{t=1}^n d_a^{t,s} - R\hspace{3mm}s=1,2,\cdots,N\label{mmm15}\\
    & b_p^{s} \geq e_p^{1,s} + \sum_{c=1}^k \sum_{t=(c-1)r+1}^{cr} \sum_{l=1}^r \mu_l^{s} x_{lct} + \sum_{t=1}^n d_p^{t,s} - R\hspace{3mm}s=1,2,\cdots,N\label{mmm16}\\
    % & W_p(\sigma)= \sum_{t=1}^n \bar{w}^t_p\label{mmm23}\\
    % & U_p(\sigma)=\alpha W_p(\sigma)\label{mmm24}\\
    % & D_a(\sigma)= \sum_{t=1}^n d^t_a\label{mmm25}\\
    % & I_a(\sigma)=\beta_a D_a(\sigma)\label{mmm26}\\
    % & D_p(\sigma)= \sum_{t=1}^n d^t_p\label{mmm27}\\
    % & I_p(\sigma)=\beta_p D_p(\sigma)\label{mmm28}\\
    % & O_a(\sigma)=o_a b_a\label{mmm29}\\
    % &  O_p(\sigma)=o_p b_p\label{mmm30}\\
    & x_{lct}\in \{0,1\}\hspace{3mm} l=1,2,\cdots,r;\hspace{3mm} c=1,2,\cdots,k;\hspace{3mm} t=1,2,\cdots,n\label{mmm17}\\
    & y_t\in \{0,1\}\hspace{3mm} t=1,2,\cdots,n\label{mmm18}\\
    & d_a^{t,s}, d_p^{t,s}, e_a^{t,s}, e_p^{t,s}, w_a^{t,s}, \bar{w}_p^{t,s}, \tau_a^t, b_a^{s}, b_p^{s} \geq 0\hspace{3mm}   t=1,2,\cdots,n\hspace{3mm}s=1,2,\cdots,N\label{mmm19}
    %& M\text{: sufficiently large number}\label{mm20} \nonumber
\end{align}

\newpage

\noindent \textbf{SAA Model for a Single Block}

\noindent
Because the number of sample paths typically is quite large, finding an optimal objective function value, $\Psi^*_s$ and its solution for the stochastic single block problem can quickly become computationally intractable. We apply the Sample Average Approximation (SAA) Method (Kleywegt et al. 2001) to obtain an approximate solution close to the optimum for the stochastic single period problem. We briefly describe the SAA method: (i) First formulate the single block problem to obtain a solution for a sample of size $K$ which has the objective function of $\psi_K(X)=\frac{1}{K} \sum_{s=1}^K  \Pi_s(w_a^{j,s},\bar{w}_p^{j,s},d_a^{j,s},d_p^{j,s})$ as given below, (ii) then solve the single block problem formulated below repetitively $\nu$ times using a smaller $K$ to obtain an approximate solution where $\nu$ is the number of replications. Detailed steps are provided below. 

 \vspace{-1mm}
\begin{align}
    \min\; & \psi_K(X)=\frac{1}{K} \sum_{s=1}^K  \Pi_s(w_a^{j,s},\bar{w}_p^{j,s},d_a^{j,s},d_p^{j,s}) \label{m170}\\
    \text{subject to } & \Pi_s=\alpha \sum_{t=1}^n w^{t,s}_a + \alpha \sum_{j=1}^r \bar{w}^{j,s}_p+\beta_a \sum_{j=1}^r d^{j,s}_a+\beta_p \sum_{j=1}^r d^{j,s}_p\hspace{3mm}s=1,2,\cdots,K\label{m290}\\
    & \sum_{j=1}^r x_{lj}=1\hspace{3mm} l=1,2,\cdots,r\label{m180}\\
    & \sum_{l=1}^{r}  x_{lj}=1\hspace{3mm}  j=1,2,\cdots,r\label{m190}\\
    & \sum_{z\in\{1,2,\cdots,q\}} x_{z1}=0\label{m360}\\
    & e_a^{j+1,s}\geq e_a^{j,s} + \sum_{l=1}^r x_{lj} \lambda_{l,s} \hspace{3mm}  j=1,2,\cdots,r-1;\hspace{3mm}s=1,2,\cdots,K\label{m200}\\
    & e_p^{j,s} \geq e_a^{j,s} + \sum_{l=1}^r x_{lj} \lambda_{l,s} \hspace{3mm} j=1,2,\cdots,r;\hspace{3mm}s=1,2,\cdots,K\label{m210}\\
    & e_p^{j+1,s} \geq e_p^{j,s} + \sum_{l=1}^r x_{lj} \mu_{l,s} \hspace{3mm} j=1,2,\cdots,r-1;\hspace{3mm}s=1,2,\cdots,K\label{m220}\\
    & d_a^{j+1,s} \geq e_a^{j+1,s}-e_a^{j,s}-\sum_{l=1}^r x_{lj} \lambda_l \hspace{3mm} j=1,2,\cdots,r-1;\hspace{3mm}s=1,2,\cdots,K\label{m230}\\
    & d_p^{j+1,s} \geq e_p^{j+1,s}-e_p^{j,s}-\sum_{l=1}^r x_{lj}\mu_l \hspace{3mm}  j=1,2,\cdots,r-1;\hspace{3mm}s=1,2,\cdots,K\label{m240}\\
    & e_a^{1,s}=0\label{m250}\\
    % & w_a^{j,s} \geq e_a^{j,s}-e_a^{j-1,s}-\sum_{l=1}^r x_{l,j-1}\lambda_{l,s}\hspace{3mm} j=2,\cdots,r;\hspace{3mm}s=1,2,\cdots,N\\
    & \tau_a^j \leq e_a^{j,s}\hspace{3mm} j=2,\cdots,r;\hspace{3mm}s=1,2,\cdots,N\\
    & w_a^{j,s} \geq e_a^{j,s}-\tau_a^j \hspace{3mm} j=2,\cdots,r;\hspace{3mm}s=1,2,\cdots,N \\
    & \bar{w}_p^{j,s} \geq e_p^{j,s}-e_a^{j,s}-\sum_{l=1}^r x_{lj} \lambda_{l,s} -M(1-y_j)\hspace{3mm} j=1,2,\cdots,r;\hspace{3mm}s=1,2,\cdots,K\label{m260}\\
    & \frac{y_j}{M} \leq \sum_{l=1}^r \mu_{l,s} x_{lj} \hspace{3mm} j=1,2,\cdots,r;\hspace{3mm}s=1,2,\cdots,K\label{m270}\\
    & My_j \geq \sum_{l=1}^r \mu_{l,s} x_{lj}\hspace{3mm} j=1,2,\cdots,r;\hspace{3mm}s=1,2,\cdots,K\label{m280}\\
    & x_{lj}\in \{0,1\}\hspace{3mm} l=1,2,\cdots,r;\hspace{3mm} j=1,2,\cdots,r\label{m310}\\
    & y_j\in \{0,1\}\hspace{3mm} j=1,2,\cdots,r\label{m320}\\
    & d_a^{j,s}, d_p^{j,s}, e_a^{j,s}, e_p^{j,s}, w_a^{j,s},  \bar{w}_p^{j,s}, \tau_a^j \geq 0\hspace{3mm}   j=1,2,\cdots,r;\hspace{3mm}s=1,2,\cdots,K\label{m330}\\
    & M\text{: sufficiently large number}\label{m340} \nonumber
\end{align}

The SAA Model assumes equal probability for the occurrence of each scenario in the objective function, where the constraints are identical to the Stochastic Programming Model in Section~\ref{sto}. A SAA model can be developed for the planning horizon similar to that of the single block case. We now illustrate the solution procedure for the single block case.

\medskip
\noindent \textbf{SAA Solution Procedure for Two-stage Block Scheduling Model with Stochastic Service Times}

Here, $X=(x_{lj})$ denotes an optimal solution to one replication of the SAA model. More specifically, $X^u=(x^u_{lj})$ denotes an optimal solution to replication $u$ of SAA model and the $x^u_{lj}$ solution corresponds to the block schedule denoted as $\sigma_u$, where $u=1,2, \ldots, \nu$. Since the original Stochastic Programming Model intends to find one optimal $X^*=(x^*_{lj})$ (corresponding block schedule, $\sigma^*$) optimizing over all sample paths, we run the replication of the SAA model along this same idea and take the $\bar{\psi}_K^\nu(X)$, the objective function value of the solution $X$.

Preserving the constraints of the Stochastic Programming Model, the SAA Model assumes equal occurrence probability for each scenario and aims to minimize the average of a selected number of scenarios ($K$), that is the sample size, replicating each scenario set $K$, $\nu$ times. The number of replications or the sample size can be increased until achieving a desired precision level. Solving the SAA Model for each sample yields the objective value of $\psi_K^u(X)$ for each replication $u$ where $u=1,\cdots,\nu$. By taking the average of all of the objective values provided by the SAA Model, we obtain an estimate of the (mean) optimal objective value of the Stochastic Programming Model as follows:
\vspace{-0.05in}
  \[\bar{\psi}_K^\nu(X)=\frac{1}{\nu}(\psi_K^1(X^1)+\cdots+\psi_K^\nu(X^\nu))\]
Using this mean estimate, the sample variance can be calculated as ${S_K}^\nu=\frac{1}{\nu}\sum_{u=1}^\nu [\psi_K^u(X^u)-\bar{\psi}_K^\nu(X)]^2$. 
Next, the $100(1-p)\%$ confidence interval for the expected optimal objective value for the problem ($E[\Psi_K]$) can be constructed as follows:
\[\bar{\psi}_K^\nu(X) \pm t_{\nu-1,1-p/2}\sqrt{\frac{S_K^\nu}{\nu-1}}\]
where the half-length of the confidence interval is equal to $h(\nu,p)= t_{\nu-1,1-p/2}\sqrt{\frac{S_K^\nu}{\nu-1}}$.

Following \cite{stauffer2021elasticity}, if the estimate $\bar{\psi}_K^\nu(X)$ satisfies the following
\[\frac{\bar{\psi}_K^\nu(X)-E[\Psi_K]}{E[\Psi_K]}=\xi\]
then it is possible to claim that $\bar{\psi}_K^\nu(X)$ has a relative error of $\xi$. Increasing the replications until $\frac{h(\nu,p)}{\bar{\psi}_K^\nu(X)} \leq \xi$ is satisfied would allow $\bar{\psi}_K^\nu(X)$ to have a relative error at most $\frac{\xi}{1-\xi}$ with probability of approximately $1-p$ \citep{law2007simulation}. Thus, in order to obtain an estimate of $E[\Psi_K]$ with a relative error of $\xi$ with $100(1-p)\%$ confidence interval, we follow the steps of the SAA Solution Procedure \citep{law2007simulation} as follows.

\vspace{0.2in}

\noindent\textbf{Step 1}: Set $\nu_0=5$ replications (with each $K=15$) and set $\nu=\nu_0$.\\
\textbf{Step 2:} Calculate $\bar{\psi}_K^\nu(X)$ and $h(\nu,p)$ from  $\psi_K^1,\cdots,\psi_K^\nu$ with a $95\%$ confidence interval.\\
\textbf{Step 3:} If $\frac{h(\nu,p)}{\bar{\psi}_k^\nu(X)} < \frac{\xi}{1+\xi}$, then use $\bar{\psi}_k^\nu(X)$ as a point estimate for $\Psi_k^*$ and stop. Otherwise, replace $\nu$ with $\nu+1$ and proceed for another replication of the simulation and go to Step 2.\\
\textbf{Step 4:} If $\nu>\nu^{\max}$ increase $K$ sample paths by $K^\text{step}$, set $\nu=\nu_0=5$ and repeat Steps 2 and 3 above.

\vspace{0.2in}

\medskip

\noindent{\bf Obtaining an estimate of the (mean) optimal objective value of the Stochastic Programming Model using SAA model with sample size $K$:}

\medskip

\noindent{\bf Notation for Replication, $u$:} \\The solution $X^u=(x^u_{lj})$ (corresponding  block schedule, $\pi$) with objective function $\psi_K^u(X^u)$, where $u=1,2, \ldots, \nu$.

The steps for calculating the estimate of the optimal objective function are provided in the pseudocode below. We explain these steps in detail next.

\vspace{0.2 in}

\noindent \textbf{Step 1:} Initially set the solution that yields the minimum average to the first solution. $\bar{X}=X_1$\\
\textbf{Step 2:} For $u=2$ to $\nu$ : Calculate $\frac{1}{u}\sum_{v=1}^u \psi^v_K(X^u)$ and $ \frac{1}{u}\sum_{v=1}^u \psi^v_K(\bar{X})$. \\
\textbf{Step 3:} If $\frac{1}{u}\sum_{v=1}^u \psi^v_K(X^u) < \frac{1}{u}\sum_{v=1}^u \psi^v_K(\bar{X})$, set $\bar{X}=X^u$ and $\bar{\psi}^u_K=\frac{1}{u}\sum_{v=1}^u \psi^v_K(\bar{X})$.

% \begin{algorithm}
%     \begin{algorithmic}
%         \caption{\textcolor{blue}{Obtaining the Optimal Objective Value for $\nu$ replications}}
%         \State $\bar{X}=X_1$
%         \For{$u = 2$ to $\nu$} 
%             \If{$\frac{1}{u}\sum_{v=1}^u \psi^v_K(X^u) < \frac{1}{u}\sum_{v=1}^u \psi^v_K(\bar{X})$}
%                 \State $\bar{X}=X^u$
%                 \State $\bar{\psi}^u_K=\frac{1}{u}\sum_{v=1}^u \psi^v_K(\bar{X})$
%             \EndIf
%          \EndFor
%     \end{algorithmic}
% \end{algorithm}

\medskip

\noindent{\bf Mean of the first two Replications:}

First Replication, $u=1$: $X^1=(x^1_{lj})$ with objective function $\psi_K^1(X^1)$ in the SAA model.

Second Replication, $u=2$: $X^2=(x^2_{lj})$ with objective function $\psi_K^2(X^2)$ in the SAA model.

Now we have two block schedules corresponding to $X^1=(x^1_{lj})$ and $X^2=(x^2_{lj})$.

\medskip

We evaluate the solution $X^1=(x^1_{lj})$ from the first replication of the SAA model in the second replication sample paths of SAA and find the corresponding objective function, $\psi_K^2(X^1)$ and then find the average of the two $\frac{1}{2}(\psi_K^1(X^1)+\psi_K^2(X^1))$.

\medskip

Similarly, we evaluate solution $X^2=(x^2_{lj})$ from the second replication of SAA model in the first replication sample paths of SAA and find the corresponding objective function, $\psi_K^1(X^2)$ and then find the average of the two $\frac{1}{2}(\psi_K^1(X^2)+\psi_K^2(X^2))$.

\medskip

We then find the minimum of the two averages and register the solution $\bar{X}$ corresponding to the minimum average as below. Note that either $\bar{X} =X^1$ or $\bar{X} =X^2$ whichever provides the minimum average:
\[\bar{\psi}_K^2= \min\{\frac{1}{2}(\psi_K^1(X^1)+\psi_K^2(X^1)), \frac{1}{2}(\psi_K^1(X^2)+\psi_K^2(X^2))\}\]

We now have the best solution $\bar{X}$ for the first two replications, i.e., 
$\bar{\psi}_K^2= \frac{1}{2}(\psi_K^1(\bar{X})+\psi_K^2(\bar{X})).$

 \medskip

\noindent{\bf Mean of first three Replications:}

Note that we have the best solution $\bar{X}$ for the first two replications.

Let the solution for the third replication, $u=3$: $X^3=(x^3_{lj})$ with objective function $\psi_K^3(X^3)$ in the SAA model.

% \medskip

Now we have two block schedule solutions corresponding to $\bar{X}$ and $X^3=(x^3_{lj})$.

\medskip

We evaluate solution $\bar{X}$ from the first and second replications of the SAA model in the third replication sample paths of SAA and find the corresponding objective function, $\psi_K^3(\bar{X})$ and then find the average of three $\frac{1}{3}(\psi_K^1(\bar{X})+\psi_K^2(\bar{X}) +\psi_K^3(\bar{X}))$.

\medskip

Similarly, we evaluate solution $X^3=(x^3_{lj})$ in the first and second replication sample paths of SAA model and find the corresponding total objective functions, $\psi_K^1(X^3)+ \psi_K^2(X^3)$ and then find the average of three as follows: $\frac{1}{3}(\psi_K^1(X^3)+\psi_K^2(X^3)+\psi_K^3(X^3))$.

\medskip

We then find the minimum of the two averages and register the solution $\bar{X}$ corresponding to the best average as found below. 

\begin{center}
    $\bar{\psi}_K^3= \min\{\frac{1}{3}(\psi_K^1(\bar{X})+\psi_K^2(\bar{X}) +\psi_K^3(\bar{X})), \frac{1}{3}(\psi_K^1(X^3)+\psi_K^2(X^3)+\psi_K^3(X^3))\}$
\end{center}

% \[\bar{\psi}_K^3= \min\{\frac{1}{3}(\psi_K^1(\bar{X})+\psi_K^2(\bar{X}) +\psi_K^3(\bar{X})), \frac{1}{3}(\psi_K^1(X^3)+\psi_K^2(X^3)+\psi_K^3(X^3))\}\]

We now have the best solution $\bar{X}$ for the first three replications  (Note that either $\bar{X} =\bar{X}$ (the previous $\bar{X}$) or $\bar{X} =X^3$ whichever provides the minimum average:), i.e., 

\begin{center}
    $\bar{\psi}_K^3= \frac{1}{3}(\psi_K^1(\bar{X})+\psi_K^2(\bar{X})+\psi_K^3(\bar{X}))$
\end{center}

% $$\bar{\psi}_K^3= \frac{1}{3}(\psi_K^1(\bar{X})+\psi_K^2(\bar{X})+\psi_K^3(\bar{X})).$$

% \medskip

\noindent{\bf Mean of $\nu$ Replications:}

Assume that we follow this procedure and find that we now have the best solution $\bar{X}$ for $(\nu-1)$ replications.

Let the solution for Replication $\nu$: $X^\nu=(x^\nu_{lj})$ with objective function $\psi_K^\nu(X^\nu)$ in the SAA model.

% \bigskip

Now we have two block schedules corresponding to $\bar{X}$ and $X^\nu=(x^\nu_{lj})$.

% \medskip

We evaluate solution $\bar{X}$ in replication $\nu$ sample paths of SAA and find corresponding objective function, $\psi_K^\nu(\bar{X})$ and then find the average, $\frac{1}{\nu}(\psi_K^1(\bar{X})+\psi_K^2(\bar{X})+ \cdots + \psi_K^{\nu-1}(\bar{X})+\psi_K^\nu(\bar{X}))$.

% \medskip

Similarly, we evaluate solution $X^\nu=(x^\nu_{lj})$ from the replication $\nu$ of SAA model in the previous $(\nu-1)$ replication sample paths of SAA models and find corresponding total objective function values, $\psi_K^1(X^\nu)+\psi_K^2(X^\nu)+ \cdots + \psi_K^{\nu-1}(X^{\nu}))$ and then find the average, $\frac{1}{\nu}(\psi_K^1(X^\nu)+\psi_K^2(X^\nu)+ \cdots + \psi_K^{\nu-1}(X^{\nu})+\psi_K^\nu(X^\nu))$.

% \medskip

We then find the minimum of the two above averages and register the solution $\bar{X}$ corresponding to the minimum average as found above. Note that either $\bar{X} = \bar{X}$ (as previous $\bar{X}$ value) or $\bar{X} =X^\nu$ whichever provides the minimum average. Then, $\bar{\psi}_K^\nu$ is set equal to this minimum average.

We now have the best solution $\bar{X}$ for $\nu$ replications with the average of $\bar{\psi}_K^\nu$.

\newpage
\section{Extension to the General Multi-stage Outpatient Clinical Setting} \label{general_setting}

\cite{huang2015appointment} consider a multi-stage outpatient clinical setting where parallel resources—such as two physician assistants ($PA$s) and two physicians ($P$s)—are used. Typically, in such a multi-stage system, each physician is assisted by one physician assistant. Therefore, a natural strategy in this case is to decompose the system into two separate subsystems, each consisting of one $PA$ and one $P$. Pre-designed single appointment blocks can then be assigned to each subsystem in such a way that the workload is balanced (i.e., each subsystem is assigned a nearly equal number of blocks). For example, if the system's capacity target is to treat four blocks, then assigning two blocks to each of two subsystems ensures balanced workload distribution. Each subsystem can then be scheduled independently according to the methodologies developed in previous sections. In general, for a system with $p$ $PA$s and $p$ $P$s, the system can be decomposed into $p$ separate subsystems, each consisting of one $PA$ and one $P$.

\vspace{5mm}
\section{Two-stage Block Scheduling Model with No-show Uncertainty}\label{robu}
\vspace{-1mm}

In this section, we study the impact of patient no-shows on block schedules. In practice, clinics typically consider overbooking, i.e., assigning more than one patient for some time slots, to improve the utilization rate of the $P$ and $PA$ when no-shows occur. 

Assuming zero patient no-show probabilities, a single block schedule provided by Algorithm~\ref{alg-B} is as follows: % $S:|$HC HC L MC MC MC MC LC L LC L LC LC M M H$|$.

\noindent $S$: $|HC \;\; HC \;\; L \;\; MC\;\;MC\;\;MC\;\;MC\;\;LC\;\;L\;\;LC\;\;L\;\;LC\;\;LC\;\;M\;\;M\;\;H|$.
 
 \noindent Note that $S$ is a no-idle block schedule from our heuristics. Following the literature \citep{lee2018outpatient}, we consider two overbooking strategies to obtain schedules.

Let $p^{+}=0.2$ and $p=0.3$ be the no-show probabilities of $Q^{+}$ and $Q$ group patient types, respectively. In the appointment schedule with overbooking, we need to overbook and schedule $e^+=10\times 0.2 = 2$ additional patients (That is, the expected number of $Q^{+}$ no-shows for a single block, $e^+$ is 2). Similarly, we need to overbook and schedule $e=0.3 \times 6=1.8$ additional $Q$ patients (That is, the expected number of $Q$ patient no-shows, $e$, is 1.8 which is rounded to an integer number, 2). Thus, the total number of overbookings needed is $e+e^+=4$ in our analysis.
We propose two overbooking block schedules. The first uses level front loading while the second uses a fully front loaded strategy.

%Note that $e^+=2$ and $e=2$ in our example. Now we propose the first overbooking schedule that front-loads the overbooked patients. That is, the number $e^+ + e =2$ of both types are scheduled in the first block where $Q^{+}$ and $Q$ appear.

% \noindent{\bf Level Overbook Schedule, $S_{L}$:} In this schedule, (i) each schedule position can be overbooked by the same type of patient, (ii) each schedule position can be overbooked at most by one patient, (iii) $(e^+ + e)$ patients are overbooked at the earliest available time, (iv) $(e^+ + e)$ patients are distributed across the block as equally as possible as shown below. \\ $S_{L}$:
% $|T_3(1)T_4T_4T_4T_2(1)T_2(1)T_1T_1T_1||T_3(1)T_4T_4T_4T_2(1)T_2T_1T_1T_1|$.\\
% $T_{3}(1)$ means the patient position of $T_{3}$ is overbooked by one more $T_{3}$ patient. $T_{2}(1)$ means the patient position of $T_{2}$ is overbooked by one more $T_{2}$ patient.

\medskip

\noindent{\bf Level Front Load Overbook Schedule, $S_{LF}$:} In this schedule, (i) each schedule position can be overbooked by the same type of patient,  (ii) each schedule position can be overbooked at most by one patient, (iii) $(e^+ + e)$ patients are overbooked at the earliest available time in the first block, and (iv) no patients are overbooked in the subsequent blocks.

\noindent $S_{LF}$:
 $|HC(1) \,\; HC(1) \,\; L(1) \,\; MC\,\;MC\,\;MC\,\;MC\,\;LC\,\;L(1)\,\;LC\,\;L\,\;LC\,\;LC\;\;M\,\;M\,\;H|$.

%\medskip

\noindent{\bf Fully Front Load Overbook Schedule, $S_{FF}$:} In this schedule, (i) each schedule position can be overbooked by the same type of patient,  (ii) $(e^+ + e)$ patients are overbooked each schedule position in any number at the earliest available time in the first block, and (iii) no patients are overbooked in the subsequent blocks. \\ $S_{FF}$: $|HC(2) \;\; HC \;\; L(2) \;\; MC\;\;MC\;\;MC\;\;MC\;\;LC\;\;L\;\;LC\;\;L\;\;LC\;\;LC\;\;M\;\;M\;\;H|$.\\
$HC(2)$ means the patient position of $HC$ is overbooked by two more $HC$ patients. $L(2)$ means the patient position of $L$ is overbooked by two more $L$ patients.

%\medskip

In the above overbooking schedules $S_{LF}$ and $S_{FF}$, 20 patients are scheduled. Thus, the number of sample paths is $2^{20}=1,048,576$ since there are two states of each patient in the schedule: (i) the patient does not show up with the probability, $p^{+}=0.2$ (or $p=0.3$) and (ii) the patient shows up with probability, $1-p^{+}=0.8$ (or $1-p=0.7$), depending on their patient group type.

%We assume positive no-show probabilities for both $Q^{+}$ and $Q$ group patients within a single block schedule generated according to parameters introduced in Section~\ref{cdp}. We also assume deterministic service times in Section~\ref{overcomp}, so the only difference is the different overbooking strategies.
% and then analyze stochastic service times in Section~\ref{over-sto}.

%\subsection{Computational Results for No-shows and Overbooking under Deterministic Service Times}\label{overcomp}
\subsection{Computational Results for No-shows and Overbooking}\label{overcomp}
\vspace{-1mm}

This section compares computational results for the three schedules, $S$ (original block schedule without overbooking), $S_{LF}$, and $S_{FF}$ by enumerating all sample paths and estimating the expected objective function values. We also compare patient waiting times, $P$, $PA$ idle times and overtimes among the three schedules. 
The service times used in this study are assumed to be deterministic and from the parameters introduced in Section~\ref{cdp} based on Table~\ref{table-newparMAIN}. As the analysis is carried out in a single block (using Algorithm~\ref{alg-B}), the regular time in a day is set to $R=150$ minutes. %For this analysis $S_{LF}$ and $S_{FF}$ strategies were implemented on the schedule provided by Algorithm~\ref{alg-B}. We analyze and compare these strategies using robust optimization method.

In Figure~\ref{no_show_2plots} results, we enumerate all possible no-show sample paths, calculate their respective probabilities, and compare different strategies based on the expected objective function values per patient for various cost parameters. Since Algorithm~\ref{alg-B} incurs zero patient waiting time in the deterministic setting, its objective function value is constant (6.77) for increasing patient waiting time cost. Without overbooking, the Algorithm~\ref{alg-B} block schedule results in decreasing patient wait times and decreasing $PA$ and $P$ overtime as the number of no-shows increases. However, the $PA$ and $P$ idle time is increasing. The left graph of Figure~\ref{no_show_2plots} illustrates similar cost parameters to Section~\ref{compstu}. Here we see that Algorithm~\ref{alg-B} provides a lower overall cost at all but the lowest patient waiting time costs. If we raise the relative overtime costs to $o_a=o_p=1.5$ and $o_a=o_p=1.8$, the Algorithm~\ref{alg-B} total cost per patient only increases to 7.33 and 7.89, respectively. This illustrates that Algorithm~\ref{alg-B} is relatively robust to no-show behavior.

\begin{figure}[thb]
    \centering %[width=0.5\linewidth]
    \includegraphics[scale=0.80]{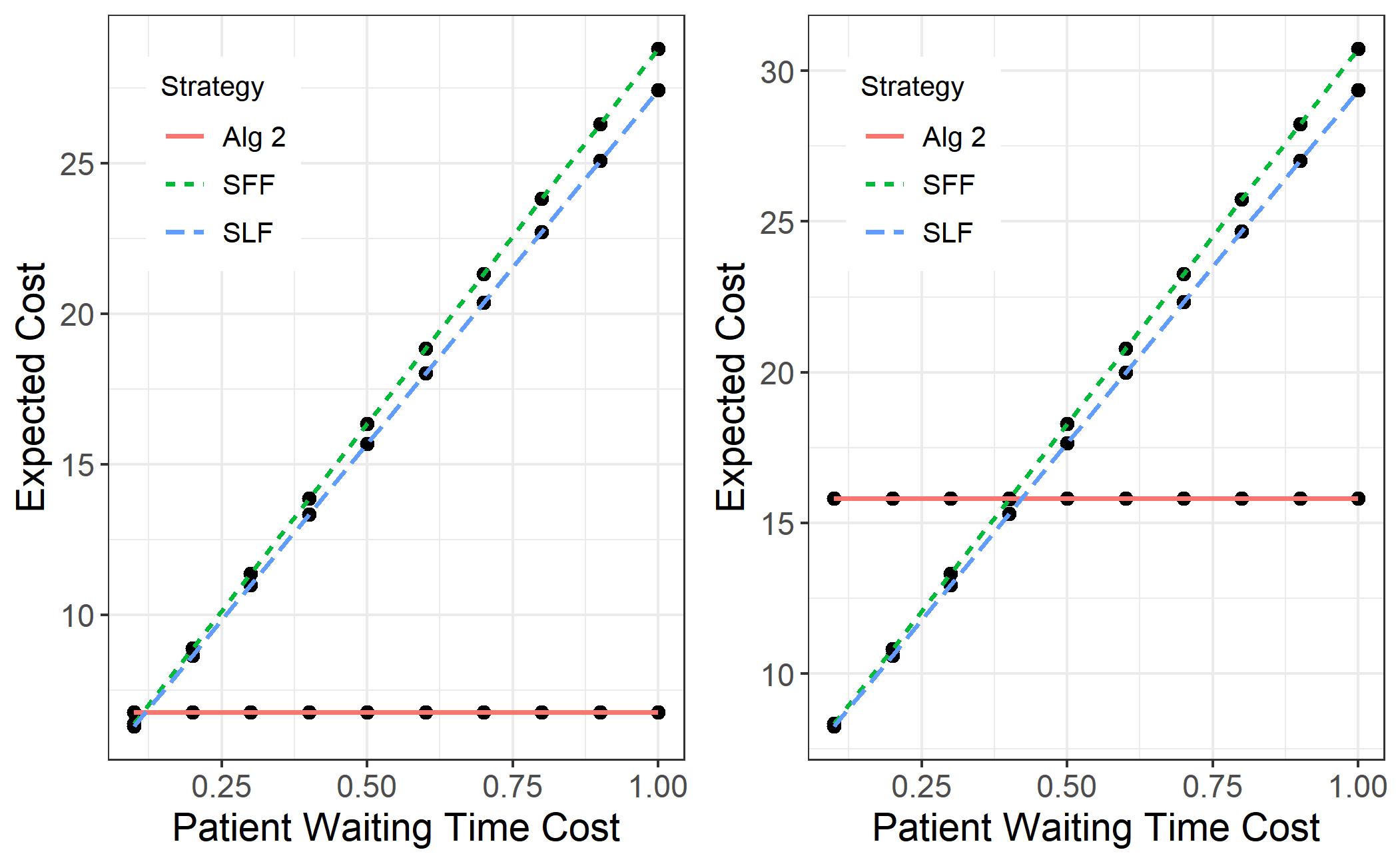}
    \caption{The expected objective function values per (scheduled) patient for Algorithm~\ref{alg-B}, $S_{LF}$ and $S_{FF}$
    where patient waiting time cost, $\alpha\in\{0.1,\cdots,1\}$, the overtime costs for $PA$ and $P$, $o_a=o_p=1.2$, and regular time in the day, $R=150$.
    The idle time costs for $PA$ and $P$, $\beta_a=\beta_p=1$ for the graph on the left of the figure and $\beta_a=\beta_p=3$ for the graph on the right of the figure.}
    \label{no_show_2plots}
\end{figure}
\vspace{-1mm}

In the right graph of Figure~\ref{no_show_2plots}, we increase the idle time cost from 1 to 3, which could occur in clinics that value $PA$ and $P$ idle time much more than any other cost. In this case, the overbooking strategies provide a better solution than Algorithm~\ref{alg-B} for a wider range of patient waiting time costs. This occurs because the overbooking strategies have shorter idle times but much longer patient wait times per sample path than Algorithm~\ref{alg-B}. The overtime amount is slightly less for Algorithm~\ref{alg-B}, but the real trade-off is between idle time and wait time. Thus, which strategy a clinic prefers depends on how that clinic values those two cost parameters. 
Among the overbooking strategies considered, we observe that $S_{LF}$ provides a lower per-patient cost than $S_{FF}$. This occurs because we enumerate all possible no-show positions and for all sample paths with no-shows later in the day, $S_{LF}$ performs better with a few overbooked time slots later in the day.

\end{APPENDICES}

% Appendix here
% Options are (1) APPENDIX (with or without general title) or
%             (2) APPENDICES (if it has more than one unrelated sections)
% Outcomment the appropriate case if necessary
%
% \begin{APPENDIX}{<Title of the Appendix>}
% \end{APPENDIX}
%
%   or
%
% \begin{APPENDICES}
% \section{<Title of Section A>}
% \section{<Title of Section B>}
% etc
% \end{APPENDICES}

% Acknowledgments here
% \ACKNOWLEDGMENT{The authors gratefully acknowledge the existence of
% the Journal of Irreproducible Results and the support of the Society
% for the Preservation of Inane Research.}

% References here (outcomment the appropriate case)

% CASE 1: BiBTeX used to constantly update the references
%   (while the paper is being written).
% \bibliographystyle{informs2014} % outcomment this and next line in Case 1
%\bibliography{<your bib file(s)>} % if more than one, comma separated

% CASE 2: BiBTeX used to generate mypaper.bbl (to be further fine tuned)
%\input{mypaper.bbl} % outcomment this line in Case 2

%If you don't use BiBTex, you can manually itemize references as shown below.

%%%%%%%%%%%%%%%%%
\end{document}